\documentclass[a4paper,11pt]{article}
\usepackage{graphicx}
\usepackage{epsfig}
\usepackage{amsmath}
\usepackage{amssymb}
\usepackage{mathtools}
\usepackage[top=3cm, bottom=3.5cm, left=2cm, right=2cm]{geometry}
\usepackage{microtype}
\usepackage{lmodern}
\usepackage[utf8]{inputenc}
\usepackage{bera}

\DeclarePairedDelimiter{\ceil}{\lceil}{\rceil}

\newcommand{\ds}{\displaystyle }

\newcommand{\bee}{\begin{equation}}
\newcommand{\eed}{\end{equation}}
\author{{\Large Thomas M. Michelitsch$^a$, Federico Polito$^b$ and
Alejandro P. Riascos$^c$ }\\ \\ 
\footnotesize{$^a$ Sorbonne Universit\'e, Institut Jean le Rond d’Alembert, 
CNRS UMR 7190} \\
\footnotesize{4 place Jussieu, 75252 Paris cedex 05, France} \\ 
\footnotesize{E-mail: michel@lmm.jussieu.fr}\\[1ex]
\footnotesize{$^b$ Department of Mathematics ``Giuseppe Peano'', University of Torino, Italy}  \\ 
\footnotesize{E-mail: federico.polito@unito.it}\\[1ex]
\footnotesize{$^c$Instituto de Fisica, Universidad Nacional Aut\'onoma de M\'exico} \\ 
\footnotesize{Apartado Postal 20-364, 01000 Ciudad de M\'exico, M\'exico}  \\
\footnotesize{E-mail: aperezr@fisica.unam.mx}
}
\title{ON DISCRETE TIME PRABHAKAR-GENERALIZED FRACTIONAL
POISSON PROCESSES AND RELATED STOCHASTIC DYNAMICS}

\begin{document}

\maketitle

\begin{abstract}
Recently the so-called Prabhakar generalization of the fractional Poisson counting process attracted much interest
for his flexibility to adapt to real world situations.
In this renewal process the waiting times between events are IID {\it continuous} random variables. In the present paper we analyze discrete-time counterparts: 
Renewal processes with {\it integer} IID interarrival times which converge in  
well-scaled continuous-time limits to the Prabhakar-generalized fractional Poisson process. 
These processes exhibit
non-Markovian features and long-time memory effects.
We recover for special choices of parameters the discrete-time 
versions of classical cases, such as the fractional Bernoulli process and the standard Bernoulli process as discrete-time approximations of the fractional Poisson and the standard Poisson process, respectively.
We derive difference equations of generalized
fractional type that govern these discrete time-processes where in well-scaled
continuous-time limits known evolution equations of generalized fractional Prabhakar type are recovered.
We also develop in Montroll-Weiss fashion the ``Prabhakar Discrete-time random walk (DTRW)'' 
as a random walk on a graph 
time-changed with a discrete-time version of Prabhakar renewal process.
We derive the generalized fractional discrete-time Kolmogorov-Feller difference
equations governing the resulting stochastic motion. Prabhakar-discrete-time processes open a promising field capturing several aspects in the dynamics of complex systems.
\\ \\
{\it Keywords:\\[1ex]
Discrete-time renewal processes; Generalized Kolmogorov-Feller equations;\\
Prabhakar fractional calculus; Non-Markov random walks on graphs}
\end{abstract}

\section{\small INTRODUCTION}
\label{GFPP}
Classically, random occurrence of events in time is modeled 
as standard Poisson processes, i.e.\ with exponentially distributed 
waiting-times and the memoryless Markovian property. Applying this model to random motions in continuous spaces or graphs lead to continuous-time Markov chains.
However, one has recognized that many phenomena in anomalous transport and diffusion including the dynamics in certain \emph{complex} systems 
exhibit power law distributed waiting-times with non-Markovian long-time memory features that are not compatible with classical exponential patterns \cite{KutnerMasoliver1990,GorenfloMainardi2006,Gorenflo2010, MetzlerKlafter2000,MetzlerKlafter2004}. %
\\[1ex]
A powerful approach to tackle these phenomena is obtained by admitting fat-tailed power-law waiting-time densities where Mittag-Leffler functions come into play as natural generalizations of the classical exponentials. A prototypical example is the \emph{fractional Poisson process} (FPP), a counting process with unit-size jumps and IID Mittag-Leffler distributed waiting-times.  Fractional diffusion equations governing the fractional Poisson process and many other related properties are already present in the specialized literature --- see e.g.\
\cite{HilferAnton1995,RepinSaichev2000,Laskin2003,
Laskin2009,MainardiGorenfloScalas2004,GorenfloMainardi2013, Zaslavsky2002,BeghinOrsinger2009,MeerschaertEtal2011,Meerschert-et-al2018}. 
\\[1ex]
Broadly speaking, 
Markov chains permitting arbitrary waiting times define so called 
semi-Markov processes \cite{PachonPolitoRicciuti2018}.
This area was introduced independently 
by L\'evy \cite{Levy1956} Smith \cite{Smith1955} and 
Tak\'acs \cite{Takacs1958} and fundamentals of 
this theory were derived by Pyke \cite{Pyke1961} and 
Feller \cite{Feller1964} among others. 
In these classical models, semi-Markov processes have as special cases continuous-time renewal processes, i.e. 
the waiting times are IID absolutely continuous random variables.
On the other hand, many applications
require intrinsic discrete-time scales, and thus semi-Markov processes where the waiting times are discrete integer random variables open an interesting field which merits deeper analysis.
Discrete-time renewal processes are relatively little touched in the literature compared 
to their continuous-time counterparts. A discrete variant of the above-mentioned Mittag-Leffler distribution was derived by Pillai and Jayakumar \cite{PilaiJayakumar1995}.
An application in terms of a discrete-time random walk (DTRW) diffusive transport
model was developed recently \cite{Angstmann-et-al2017}.
A general approach for discrete-time semi-Markov process and 
time-fractional difference equations was developed in a recent contribution by Pachon, Polito and Ricciuti \cite{PachonPolitoRicciuti2018}. The aim of the present paper is to develop new pertinent discrete-time counting processes that contain for certain parameter choices classical counterparts such as fractional Bernoulli and standard Bernoulli, as well as in well-scaled continuous-time limits their classical continuous-time counterparts such as fractional Poisson and standard Poisson. A further goal of this paper is to analyze the resulting stochastic dynamics
on graphs.
\\[1ex]
 Before we outline the detailed organisation of the present paper, let us briefly give a general overview on the principal lines. The paper has two parts. The first part is concerned with the development of a discrete version of the Prabhakar renewal process where we prove by a `well-scaled continuous-time limit procedure' the connection with the continuous-time Prabhakar renewal process. This part includes Sections \ref{Prabhakar-ren-process}-\ref{DT} and Appendices \ref{generfuncts}-\ref{Sib} where we recall the needed mathematical machinery. 
\\[1ex]
In the second part of the paper we introduce the Prabhakar Discrete-Time-Random Walk (Prabhakar DTRW) on (undirected) graphs as a normal walk on the graph time-changed with a Prabhakar discrete-time renewal process (Sections \ref{CTRM-Prabhakar}, \ref{influence-ini}).
For this stochastic motion we derive the evolution equation of generalized fractional difference type. The limits of Poisson and fractional Poisson are also considered which recover classical results known in the literature. Supplementary materials to this part can be found in the Appendices \ref{Sibuyarenewal}-\ref{Sibuya-DTRW}. 
\\[1ex]
The detailed organisation of our paper is as follows. As a point of departure, we introduce a class of discrete-time renewal processes which represent approximations of the 
continuous-time Prabhakar process. The Prabhakar renewal process was first introduced by Cahoy and Polito \cite{PolitoCahoy2013} and the continuous-time random walk (CTRW) model based on this process was developed by Michelitsch and Riascos \cite{TMM-APR-PhysicaA2020,MichelitschRiascosGFPP2019}.
Generalized diffusion equations with time fractional Prabhakar derivative and tempered time fractional Prabhakar derivative have been derived by Sandev et al \cite{Sandev-et-al2018}.
We describe the Prabhakar renewal process in Section \ref{Prabhakar-ren-process}. For a thorough review of properties and definitions of Prabhakar-related fractional calculus
we refer to the recent review article of Giusti et al \cite{Giusti2019}.
\\[1ex]
Section \ref{discrete-time-Prab} is devoted to derive discrete-time versions 
of the Prabhakar renewal process using a composition of two `simple' processes. Then, in Section \ref{Continuous-time-GFPP}
we show that under suitable scaling conditions
the continuous-time Prabhakar process is recovered.
\\[1ex]
In Section \ref{generalization} we develop a general procedure to generate discrete-time approximations of the 
Prabhakar process. We construct these processes in such a way that the waiting time distributions are vanishing at $t=0$. 
The choice of this condition turns out to be
crucial to obtain state-probabilities allowing to define proper Cauchy initial value problems in stochastic motions.
\\[1ex]
As a prototypical example, we analyze in Section \ref{DT} the most 
simple version of discrete-time Prabhakar process 
with above mentioned good initial conditions.
We call this version of Prabhakar discrete-time counting process the `Prabhakar Discrete-Time Process' (PDTP). We derive the state-probabilities 
(probabilities for $n$ arrivals in a given time interval), 
i.e. the discrete-time counterpart of the Prabhakar-generalized fractional Poisson distribution
which was deduced in the references \cite{PolitoCahoy2013,TMM-APR-PhysicaA2020,MichelitschRiascosGFPP2019}.
The PDTP is defined as a generalization of the `fractional Bernoulli process' introduced
in \cite{PachonPolitoRicciuti2018} which is contained for a certain choice of parameters as well as the
standard Bernoulli counting process. We prove these connections by means of explicit formulas.
We show explicitly that the discrete-time waiting time and state distributions of a PDTP converge in well-scaled continuous-time limits to their known continuous-time Prabhakar function type counterparts.
These results contain for a certain choice of parameters the well-known classical cases of fractional Poisson and 
standard Poisson distributions, respectively. We show that the well-scaled continuous-time limits yield the state probabilities
of Laskin's fractional Poisson \cite{Laskin2003} and standard Poisson distributions, respectively. 
\\[1ex]
In Section \ref{stae-diffeqs} we derive for the PDTP the discrete-time versions of the generalized fractional Kolmogorov-Feller equations 
that are solved by the PDTP state-probabilities. These equations constitute discrete-time convolutions of generalized fractional type reflecting long-time
memory effects and non-Markovian features (unless in the classical standard Bernoulli with Poisson continuous-time limit case). We show that discrete-time fractional Bernoulli and standard Bernoulli processes are contained for certain choice parameters and that the same is true for their continuous-time limits: They recover the classical Kolmogorov-Feller equations of fractional Poisson and standard Poisson, respectively.
\\[1ex]
Section \ref{avaerga-number} is devoted to the analysis of the expected number of arrivals and their asymptotic features. This part is motivated by the important role 
of this quantity for a wide class of diffusion problems and stochastic motions 
in networks and lattices. 
\\[1ex]
As an application we develop in Section \ref{CTRM-Prabhakar} in Montroll-Weiss fashion the 'Discrete-Time-Random Walk' (DTRW) on undirected networks and analyze a normal random walk subordinated to the PDTP. We call this walk the `Prabhakar DTRW'.
The developed DTRW approach is a general model to subordinate random walks on graphs to discrete-time
counting processes. Although we focus on undirected graphs the DTRW approach can be extended
to general walks such as on directed graphs or strictly increasing walks on the integer line.
Such an example is briefly outlined in Appendix \ref{Sibuyarenewal}, namely a strictly increasing walk subordinated to the Sibuya counting process.
\\[1ex]
Further we derive for the Prabhakar DTRW discrete-time Kolmogorov-Feller generalized fractional difference equations that govern the resulting stochastic motion on undirected graphs and 
demonstrate by explicit formulas the contained classical cases of fractional Bernoulli and standard Bernoulli
and their fractional Poisson and Poisson continuous-time limits, respectively.
The applications of this section are motivated by the huge upswing of network 
science which has become a rapidly
growing interdisciplinary field 
\cite{NohRieger2004,Newman2010,RiascosWang-Mi-Mi2019,TMM-APR-ISTE2019,RiascosMateos2014,RiascosBoyerHerringerMateos2019,RiascosMateos2020} (and see the references therein).
\\[1ex]
Section \ref{influence-ini} is devoted to analyze the influence of the initial condition in the discrete-time waiting time density on DTRW features. We explore resulting effects introducing the feature of uncertainty and randomness in the DTRW initial conditions. 
\\[1ex]
All proofs in the paper are accompanied by detailed derivations and supplementary materials in the Appendices.

\section{\small PRABHAKAR CONTINUOUS-TIME RENEWAL PROCESS}
\label{Prabhakar-ren-process}

Among several generalizations of the fractional Poisson process which were proposed in the literature, the so called Prabhakar type generalization which
we refer to as `{\it Generalized Fractional 
Poisson Process (GFPP)}' or also `{\it Prabhakar process}' seems to be one of the most pertinent candidates. The GFPP was first
introduced by Cahoy and Polito \cite{PolitoCahoy2013} and applied to stochastic motions in networks and lattices by Michelitsch and Riascos 
\cite{TMM-APR-PhysicaA2020,MichelitschRiascosGFPP2019,michel-riascos-springer2020}. 

The Prabhakar function which is a three-parameter generalization of the Mittag-Leffler function was introduced in 1971 by Prabhakar \cite{Prabhakar1971}
and has attracted much attention recently due to its great flexibility to adapt to real-world situations. Meanwhile, the Prabhakar function has been 
identified as a matter of great interest worthy of thorough investigation. 
For a comprehensive review of properties and physical applications with generalized fractional calculus emerging from Prabhakar functions we refer to the recent review article by Giusti et al. \cite{GiustiPolitoMainardi-etal2020} and consult also \cite{GaraGorenfloPolitoTomovski2014}.

The interesting feature of the Prabhakar process is that it contains the fractional Poisson process as 
well as the Erlang- and standard Poisson processes as special cases.
The related Prabhakar-generalized fractional derivative operators may be considered 
as among the most sophisticated tools to cover certain aspects of complexity in physical systems \cite{dosSantos2019,MainardiGarrappa2015}. 
The continuous-time Prabhakar renewal process is characterized by waiting time density with Laplace transform \cite{PolitoCahoy2013}
\begin{equation}
\label{GFPPcontinuous} 
{\tilde \chi}_{\alpha,\nu }(s)=\frac{\xi_0^{\nu }}{(\xi_0+s^{\alpha})^{\nu }} ,\hspace{0.75cm} \xi_0>0, \hspace{1cm} 0 <\alpha \leq 1, \hspace{0.25cm} \nu  >0 .
\end{equation}
Laplace inversion yields 
the waiting-time PDF of the GFPP \cite{PolitoCahoy2013,TMM-APR-PhysicaA2020}
\begin{equation}
\label{Prabdensity}
  \chi_{\alpha,\nu }(t)= \xi_0^{\nu } t^{\alpha\nu -1} \sum_{m=0}^{\infty}
 \frac{(\nu )_m}{m!}\frac{(-\xi_0 t^{\alpha})^m}{\Gamma(\alpha m+ \nu \alpha)} = 
 \xi_0^{\nu } t^{\nu \alpha-1} E_{\alpha,\nu \alpha}^{\nu }(-\xi_0 t^{\alpha}) ,\hspace{1cm} t \in \mathbb{R}^{+},
 \end{equation}
which we refer to as Prabhakar-Mittag-Leffler density. The choice of this name is since this 
expression appears as a generalization of the 
Mittag-Leffler density (and recovering the Mittag-Leffler density for $\nu=1$). Expression (\ref{Prabdensity}) contains the Prabhakar-Mittag-Leffler 
function (also referred to as Prabhakar function) \cite{Prabhakar1971} defined by
\begin{equation} 
 \label{genmittag-Leff}
 E_{a,b}^c(z) = \sum_{m=0}^{\infty}
 \frac{(c)_m}{m!}\frac{z^m}{\Gamma(am + b)} ,\hspace{0.5cm} \Re\{a\} >0 ,
 \hspace{0.5cm} \Re\{b\} > 0 ,\hspace{0.5cm} c,\, z \in \mathbb{C},
\end{equation}
where $(c)_m$ indicates the Pochhammer-symbol 
\begin{equation}
\label{Pochhammer}
(c)_m =  
\frac{\Gamma(c+m)}{\Gamma(c)} = \left\{\begin{array}{l} 1 ,\hspace{1cm} m=0, \\ \\ 
                                                            
     c(c+1)\ldots (c+m-1) ,\hspace{1cm} m=1,2,\ldots \end{array}\right. 
\end{equation}
Further aspects on generalizations of Mittag-Leffler functions such as the Prabhakar function are outlined in \cite{Giusti2019,ShulkaPrajabati2007}, 
and for an analysis of properties and applications we refer to the references
\cite{GiustiPolitoMainardi-etal2020,GaraGorenfloPolitoTomovski2014,HauboldMathaiSaxena2011,BeghinOrsinger2010,Mathai2010}. 
The GFPP recovers for $\nu =1$ with $0<\alpha<1$ the Laskin fractional Poisson process  \cite{Laskin2003}, 
for $\nu >0$ with $\alpha=1$ the (generalized) Erlang process, and for
$\alpha=1$, $\nu =1$ the standard Poisson process and their related distributions. 
For details and derivations consult \cite{PolitoCahoy2013,TMM-APR-PhysicaA2020,MichelitschRiascosGFPP2019,michel-riascos-springer2020}. 

\section{\small DISCRETE-TIME VARIANTS OF THE GFPP}
\label{discrete-time-Prab}
This section is devoted to the construction of discrete-time variants of the Prabhakar renewal process 
by means of a composition of two `simple' discrete-time processes. 
First of all we recall the concept of `{\it discrete-time renewal process}'
where also the term ‘{\it discrete-time renewal chain}’ is used in the literature \cite{PachonPolitoRicciuti2018,Barbu2008}.
\\[1ex]
We introduce the strictly increasing random walk $X= (X_n)_{n\geq 1}$, such that
\begin{equation}
\label{compound-var1}
X_n =\sum_{j=1}^nZ_j ,\hspace{1cm} Z_j \in \mathbb{N} ,\hspace{1cm} X_0=0,
\end{equation}
where the steps are non-zero IID integer random variables
$Z_j=k \in \mathbb{N}$ a.s., following each the same distribution $\mathbb{P} (Z_j=k)=w(k)$. With the choice of $w(0)=0$ the walk
(\ref{compound-var1}) becomes strictly increasing.
A random walk $X$ defined in (\ref{compound-var1}) is the natural discrete-time counterpart to a (strictly increasing) subordinator \cite{MeerschaertEtal2011,PachonPolitoRicciuti2018} (and see the references therein).
\\[1ex]
In a discrete-time renewal process
the random integers $X_n$ of (\ref{compound-var1}) 
indicate the times when events occur;  we refer them to as `arrival times' or `renewal times' where $n \in \mathbb{N}_0$ 
counts the events. We also use the terms `renewals' and `arrivals'. 
Let us now introduce the generating function of the waiting-time distribution $\mathbb{P}(Z=k)= w(k)$ as
\begin{equation}\label{genfuone}
\mathbb{E} u^Z ={\bar w}(u)= \sum_{k=0}^{\infty} u^kw(k) = w(1)u+w(2)u^2+\ldots ,\hspace{1cm} |u| \leq 1,
\end{equation} 
where ${\bar w}(u)|_{u=1}=1$ reflects normalization of the $w$-distribution.
Generally generating functions are highly elegant and powerful tools which we will use extensively in the present paper. 
For some definitions and properties we refer to Appendix \ref{generfuncts}.
\\[1ex]
Consider now two discrete-time renewal processes, I and II, having waiting-time distributions $w_I$, $w_{II}$, defined
in the above general way having both zero initial conditions $w_{I}(0)=w_{II}(0)=0$. 
Then we generate a new discrete-time renewal process resulting from a composition of these two `elementary' processes. Specifically, its waiting-time distribution $\mathcal{W}(k)$ has generating function such that
\begin{equation}
\label{genfucombi} 
{\bar {\cal W}}(u) =\sum_{n=1}^{\infty}w_{II}(n)({\bar w}_I(u))^n = {\bar w}_{II}({\bar w}_I(u)),
\end{equation}
with
\begin{equation}
\label{expectval}
({\bar w}_I(u))^n=\mathbb{E}_{w_I} u^{X} = \sum_{k_1=1}^{\infty}\sum_{k_2=1}^{\infty}\ldots \sum_{k_n=1}^{\infty}w_I(k_1)w_I(k_2)\ldots w_I(k_n)u^{\sum_{j=1}^n k_j}  ,\hspace{1cm} X_0=0,
\end{equation}
where $X=(X_n)_{n\geq 1}$ is the partial sum \eqref{compound-var1} in which the random jumps are $w_I$-distributed and is characterized by the generating function (\ref{expectval}). The event counter $n$ is then considered random in (\ref{genfucombi}) 
with distribution $w_{II}$.  
We observe that ${\bar {\cal W}}(u=1)=1$ reflects normalization of the new ${\cal W}$-distribution which furthermore fulfills the desired initial condition ${\bar {\cal W}}(u=0)={\cal W}(t=0)=0$. This new waiting-time distribution then is characterized by the probabilities
\begin{equation}
\label{Wdistri}
\mathbb{P}(Z=t) = {\cal W}(t)= \frac{1}{t!} \frac{d^t}{du^t}{\bar {\cal W}}(u)\Big|_{u=0}= \sum_{n=1}^{\infty}w_{II}(n)(w_I\star)^n(t) ,\hspace{0.5cm} t\in \mathbb{N}_0 ,
\end{equation}
where $(w_I\star)^n(t)$ stands for convolution power
(See Appendix \ref{generfuncts} for details.).
Now, let us assume that process I is a {\it Sibuya counting process} 
with waiting-times following `{\it Sibuya$(\alpha)$}' 
(See Appendix \ref{Sib} for definitions and some properties). For the process II 
we choose the waiting time distribution 
\begin{equation}
\label{WIIfracB}
\begin{array}{l} \ds
\ds
 w_B^{(\nu )}(k)\big|_{k=0} =0 ,\\ \\ \ds
w_B^{(\nu )}(k)= p^{\nu }q^{k-1}
(-1)^{k-1}\left(\begin{array}{l} -\nu  \\ k-1 \end{array}\right) =
\frac{(\nu )_{k-1}}{(k-1)!}
 p^{\nu }q^{k-1} , \hspace{1cm} k\in \mathbb{N},
 \end{array}
 \hspace{0.3cm} \nu  >0, \: p+q=1,
\end{equation} 
where for $\nu=1$ (\ref{WIIfracB}) yields the {\it geometric waiting-time distribution $\mathbb{P}(Z=k) = pq^{k-1}$} of the {\it Bernoulli process} \cite{PachonPolitoRicciuti2018} (See Definition 3.1 with Eqs. (53), (54) in that paper)
with $w_B^{(\nu )}(0)=0$. We further employed here the Pochhammer symbol $ (\nu )_m$ defined in (\ref{Pochhammer}). The waiting-time distribution (\ref{WIIfracB}) has generating function
\begin{equation}
\label{dfgen}
{\bar w}_B^{(\nu )}(u) = \sum_{k=1}^{\infty} p^{\nu }q^{k-1} (-1)^{k-1}\left(\begin{array}{l} -\nu  \\ k-1 \end{array}\right) u^k = \frac{up^{\nu }}{(1-qu)^{\nu }} = \frac{u\xi^{\nu}}{(\xi +1-u)^{\nu}}
\end{equation}
where we have put $\xi=\frac{p}{q}$ $\left(p=\frac{\xi}{\xi+1}, q= \frac{1}{\xi+1}\right)$ 
and ${\bar w}_B^{(\nu )}(u)\big|_{u=1}=1$ reflects normalization. For $\nu=1$, (\ref{dfgen}) recovers the generating function of the standard Bernoulli counting process.
Now we generate a new process in the above described fashion. The new discrete-time process
hence with (\ref{genfucombi}) has waiting-time generating function
\begin{equation}
\label{generalizedfrac}
\begin{array}{l} 
\ds \mathbb{E} u^Z= {\bar {\cal \psi}}_{\alpha}^{(\nu )}(u) =
{\bar w}_B^{(\nu )}({\bar w}_{\alpha}(u))=
{\bar w}_B^{(\nu )}\bigg(\left(1-(1-u)^{\alpha}\right)
\bigg)
\\ \\ \ds \hspace{1.3cm}
= (1-(1-u)^{\alpha}) {\bar \varphi}_{\alpha}^{(\nu )}(u),\hspace{0.5cm} {\bar \varphi}_{\alpha}^{(\nu )}(u) = 
 \frac{(\frac{p}{q})^{\nu}}{\left(\frac{p}{q}+(1-u)^{\alpha}\right)^{\nu }} ,\hspace{0.5cm} \nu  >0, \hspace{0.5cm} 0<\alpha \leq 1 .
 \end{array}
\end{equation}
For $\nu =1$ and $0<\alpha<1$, formula (\ref{generalizedfrac}) recovers the generating function of 
a `discrete-time Mittag-Leffler distribution' (of so-called `type A') DML$_A$ where this process
has been named
`{\it fractional Bernoulli process (type A)}', see \cite{PachonPolitoRicciuti2018} for details of classification scheme `type A and B discrete-time processes'. 
The waiting time distribution $\varphi_{\alpha}^{(\nu=1)}(t,\xi)$ defines a discrete-time approximation of the Mittag-Leffler waiting-time distribution \cite{PilaiJayakumar1995}.
Indeed (\ref{generalizedfrac}) is generating function of a discrete-time waiting time distribution which is for $\nu>0$ and $0<\alpha<1$ a {\it generalization of discrete-time fractional Bernoulli process (of `type A')} and recovers for $\nu=1$, $\alpha=1$ the generating function of the standard Bernoulli counting process. 
\\[1ex]
Our goal now is to show that the renewal process defined by generating function (\ref{generalizedfrac})
is a discrete-time version of the Prabhakar process. 
To this end we expand (\ref{generalizedfrac}) as follows
\begin{equation}
\label{expandgenfracpoisson}
\begin{array}{l} \ds
{\bar \varphi}_{\alpha}^{(\nu)}(u) =  \xi^{\nu}(1-u)^{-\alpha\nu} (1+\xi(1-u)^{-\alpha})^{-\nu } = (1+\xi^{-1}(1-u)^{\alpha})^{-\nu} , \hspace{0.5cm} \xi=\frac{p}{q}\\ \\ \ds =  \left\{\begin{array}{clc} \ds 
\sum_{m=0}^{\infty} \left(\begin{array}{l} -\nu  \\ m \end{array}\right) 
\xi^{m+\nu}(1-u)^{-\alpha(m+\nu)} , & \xi^{\frac{1}{\alpha}}|(1-u)|^{-1} <1   \\ \\ 
\ds \sum_{m=0}^{\infty} \left(\begin{array}{l} -\nu  \\ m \end{array}\right) 
\xi^{-m}(1-u)^{\alpha m} , &  \xi^{\frac{1}{\alpha}}|(1-u)|^{-1} > 1 \end{array}\right.
\\ \\ \ds
{\bar {\cal \psi}}_{\alpha}^{(\nu )}(u) =  \left(1-(1-u)^{\alpha}\right)\, {\bar \varphi}_{\alpha}^{(\nu )}(u)  \\ \\ \ds \hspace{1.3cm} 
=
\left\{\begin{array}{clc} \ds 
\sum_{m=0}^{\infty} \frac{(-1)^m (\nu )_{m} }{m!} \xi^{m+\nu}
\left\{(1-u)^{-(m+\nu)\alpha }-(1-u)^{-(m-1+\nu)\alpha  }\right\} , &\xi^{\frac{1}{\alpha}}|(1-u)|^{-1} < 1 \\ \\ \ds 
\sum_{m=0}^{\infty} \frac{(-1)^m (\nu )_{m} }{m!}\xi^{-m}
\left\{(1-u)^{\alpha m}-(1-u)^{\alpha(m+1)}\right\} , & \xi^{\frac{1}{\alpha}}|(1-u)|^{-1} > 1 \end{array}\right.
\end{array}
\end{equation}
where these series converge absolutely
and $ (\nu )_m$ indicates the above introduced Pochhammer symbol (\ref{Pochhammer}). 
By putting $u=e^{-hs}$ and $\xi(h)=\xi_0h^{\alpha}$ ($h>0$) we see that
$ \lim_{h\rightarrow 0} (1-e^{-sh})^{-1}h \xi_0^{\frac{1}{\alpha}} = s^{-1} \xi_0^{\frac{1}{\alpha}}$ thus the Laplace variable in the limit $h\to 0$ 
has to fulfill $|s|> \xi_0^{\frac{1}{\alpha}}$, and $|s|< \xi_0^{\frac{1}{\alpha}}$
in above cases, respectively.
We notice that in (\ref{expandgenfracpoisson}) 
the powers $(1-u)^{-\mu}$ can be seen as the generating functions of expected numbers of Sibuya hits (See Appendix \ref{Sib}).
\begin{figure*}[!t]
\begin{center}
\includegraphics*[width=1.0\textwidth]{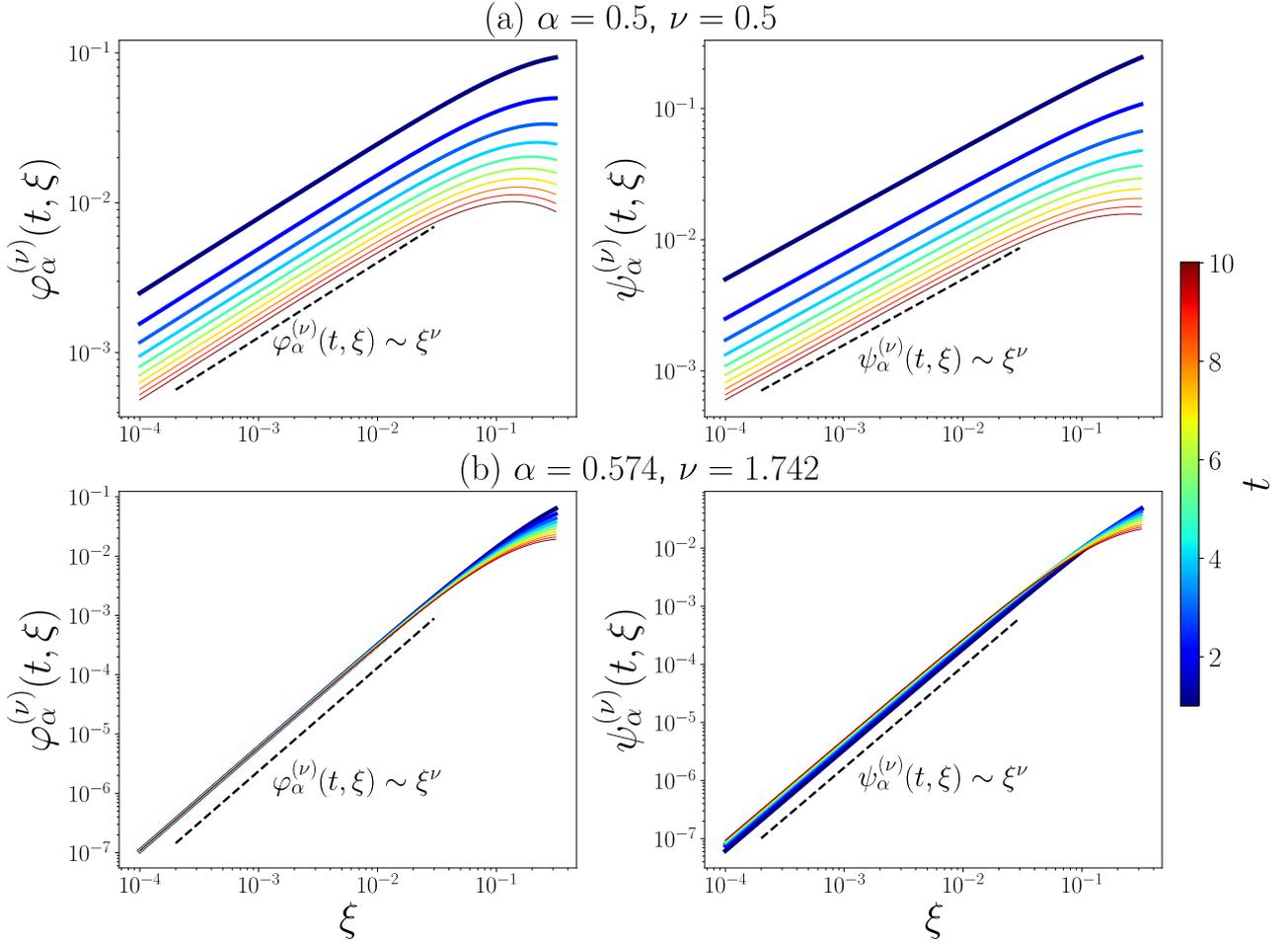}
\end{center}
\vspace{-5mm}
\caption{\label{Fig_1}\it Discrete-time densities $\varphi_{\alpha}^{(\nu )}(t,\xi)$ and $\psi_{\alpha}^{(\nu)}(t,\xi)$ versus $\xi$ for different values of $t$ indicated in the colorbar $t=1,2,\ldots,10$ calculated using expression (\ref{probabilitiesDFPP}). 
For $\xi \rightarrow 0$ they approach a power-law $\sim \xi^{\nu}$ indicated by dashed lines (See also in (\ref{probabilitiesDFPP})).
}
\end{figure*} 
Generating function (\ref{expandgenfracpoisson}) yields
the probabilities
\begin{equation}
\label{probabilitiesDFPP}
\begin{array}{l}
\ds
\varphi_{\alpha}^{(\nu )}(t,\xi)= \frac{1}{t!} \frac{d^t}{du^t}{\bar {\cal \varphi}}_{\alpha}^{(\nu )}(u)|_{u=0} = 
\left\{\begin{array}{clc} \ds 
\frac{\xi^{\nu}}{t!}\sum_{m=0}^{\infty} \frac{(-1)^m (\nu )_{m} \xi^m}{m!}  \frac{\Gamma[\alpha(m+\nu)+t]}{\Gamma[\alpha(m+\nu)]} , & 0 < \xi < 1 \\ \\
\ds \frac{(-1)^t}{t!}\sum_{m=0}^{\infty} \frac{(-1)^m (\nu )_{m} \xi^{-m}}{m!}
\frac{\Gamma[\alpha m+1]}{\Gamma[\alpha m-t+1]} , & \xi > 1
\end{array}\right.
\\ \\
\ds 
\mathbb{P}(Z=t) =
{\cal \psi}_{\alpha}^{(\nu )}(t,\xi) = \frac{1}{t!} \frac{d^t}{du^t}{\bar {\cal \psi}}_{\alpha}^{(\nu )}(u)|_{u=0}
\\ \\ \hspace{2cm} \ds = \left\{\begin{array}{clc} \ds \frac{\xi^{\nu}}{t!}
\sum_{m=0}^{\infty} \frac{(-1)^m(\nu )_{m} \xi^{m}}{m!} 
\left(\frac{\Gamma[\alpha(m+\nu)+t]}{\Gamma[\alpha(m+\nu)]}-\frac{\Gamma[\alpha(m-1+\nu)+t]}{\Gamma[\alpha(m-1+\nu)]}    \right)  & 0 < \xi < 1 \\ \\ \ds 
\frac{(-1)^t}{t!}\sum_{m=0}^{\infty} \frac{(-1)^m (\nu )_{m} \xi^{-m}}{m!}
\left( \frac{\Gamma[\alpha m+1]}{\Gamma[\alpha m-t+1]} -
\frac{\Gamma[\alpha (m+1)+1]}{\Gamma[\alpha (m+1)-t+1]}\right)
 &  \xi > 1 \end{array} \right.
\end{array} \hspace{-0.5cm} t\in \mathbb{N}_0 .
\end{equation}
where these series converge absolutely. This can be seen in view of the asymptotic behavior of the
coefficients which scale for $m\rightarrow \infty $ as $\xi^m m^{t+\nu-1}$ for $0<\xi<1$ and as
$\xi^{-m} m^{t+\nu-1}$ for $\xi>1$.
We utilize throughout this paper for all distributions synonymous notations $f(t)=f(t,\xi)$ where notation $f(t,\xi)$ is used only when it is necessary to consider the $\xi$-dependence (for instance when analyzing the continuous-space limit).
\\[1ex]
The discrete-time densities of (\ref{probabilitiesDFPP}) are plotted in Figure (\ref{Fig_1}) for different values of $t$ 
as functions of the parameter $\xi$ which defines a time scale in the process. The densities 
behave like a power law similar to $\varphi_{\alpha}^{(\nu )}(0,\xi) = \frac{\xi^{\nu}}{(1+\xi)^{\nu}} \sim \xi^{\nu}$
for $\xi$ and $t$ small; the power-law is depicted with dashed lines in Figure \ref{Fig_1}.
\\[1ex]
Now we show that {\it both} of 
the distributions in (\ref{probabilitiesDFPP}) are approximations
of the Prabhakar-Mittag-Leffler density (\ref{Prabdensity}), 
but only ${\cal \psi}_{\alpha}^{(\nu )}(t,\xi)$ per construction fulfills the desired initial condition \\ ${\cal \psi}_{\alpha}^{(\nu )}(t,\xi)\big|_{t=0}=0$. 

\subsection{\small CONTINUOUS-TIME LIMIT}
\label{Continuous-time-GFPP}

We recommend to consult Appendix \ref{crucialdefinitions} where we outline properties of the shift operator ${\hat T}_{(..)}$ which we are extensively using to define `well-scaled' continuous-time limit procedures.
Further we mention that throughout the analysis to follow we utilize as synonymous 
notations $\lim x \to a\pm 0$ and $\lim x \to a\pm $
for left- and right sided limits, respectively.
\\[1ex]
Let us introduce the (scaled) `{\it discrete-time waiting time density}' (See (\ref{discrete-time-delta})-(\ref{thepropertyogdelta-contilimit}))
\begin{equation}
\label{densityformulation}
\chi_{\alpha,\nu }(t)_{h} =  {\bar {\cal \psi}}_{\alpha}^{(\nu)}({\hat T}_{-h})\delta_h(t)
= \sum_{k=1}^{\infty} {\cal \psi}_{\alpha}^{(\nu)}(k,\xi_0h^{\alpha}) \delta_h(t-kh) = \frac{1}{h} {\cal \psi}_{\alpha}^{(\nu)}\left(\frac{t}{h},\xi_0h^{\alpha}\right) ,\hspace{0.5cm} h>0 ,\hspace{0.5cm}  t \in h \mathbb{N}_0 
\end{equation}
generalizing the notion of (continuous-time) waiting-time density to the discrete-time cases. 
We employ in this paper the notation $(..)(t)_h$ for scaled quantities 
defined on $t\in h\mathbb{Z}$ and skip subscript $h$ when $h=1$ 
(Appendix \ref{crucialdefinitions}).
Note that ${\bar {\cal \psi}}_{\alpha}^{(\nu)}({\hat T}_{-h}) = 
\sum_{k=0}^{\infty} {\cal \psi}_{\alpha}^{(\nu)}(k,\xi) {\hat T}_{-hk}$ is the operator function obtained
by replacing $u \rightarrow {\hat T}_{-h}$ in the generating function of (\ref{expandgenfracpoisson}). Bear in mind that the shift operator is such that $({\hat T}_{-h})^a f(t) ={\hat T}_{-ah} f(t)=f(t-ah)$ , $a\in \mathbb{R}$ and ${\hat T}_{-h}=e^{-hD_t}$ ($D_t=\frac{d}{dt}$). 
In (\ref{densityformulation}) occur the probabilities $\mathbb{P}(Z=k)={\cal \psi}_{\alpha}^{(\nu)}(k,\xi=\xi_0h^{\alpha})$ ($ k \in \mathbb{N}_0$) of (\ref{probabilitiesDFPP}) and the discrete-time $\delta$-distribution $\delta_h(\tau)$ ($\tau \in h\mathbb{Z}$) is
defined in (\ref{discrete-time-delta}) in Appendix \ref{crucialdefinitions}. Note that the multiplier $h^{-1}$ on the right-hand side of (\ref{densityformulation}) comes into play due to the definition (\ref{discrete-time-delta}) of the discrete-time $\delta$-distribution $\delta_{h}(t)$ guaranteeing the discrete-time densities indeed have physical dimension $sec^{-1}$.
We can then write for (\ref{densityformulation}) the distributional relations 
\begin{equation}
\label{limit-consists}
\begin{array}{l}
\ds 
{\bar {\cal \psi}}_{\alpha}^{(\nu)}({\hat T}_{-h})\delta_h(t) = 
\bigg(1-(1-e^{-hD_t})^{\alpha}\bigg) \frac{\xi(h)^{\nu}}{(\xi(h)+(1-e^{-hD_t})^{\alpha})^{\nu}}\delta_h(t) ,\hspace{1cm} \xi(h)= \xi_0h^{\alpha} \\ \\
\ds \lim_{h\rightarrow 0}{\bar {\cal \psi}}_{\alpha}^{(\nu)}({\hat T}_{-h})\delta_h(t)  =  
  \frac{\xi_0^{\nu}}{(D_t^{\alpha}+\xi_0)^{\nu}}\delta(t)   =
\frac{\xi_0^{\nu}D_t^{-\alpha\nu}}{(1+\xi_0 D_t^{-\alpha})^{\nu}}\delta(t)   ,\hspace{1cm} t \in \mathbb{R}
\end{array}
\end{equation}
where in the limiting 
process only fractional integral operators $D_t^{-\beta}$ ($\beta > 0$) occur. Indeed the discrete-time density
(\ref{limit-consists}) can be conceived as a `{\it generalized fractional integral}' 
(See Appendix \ref{appendix1}, with relations
(\ref{discrete-time-delta}) - (\ref{Lapladeltah}) and consult also \cite{PachonPolitoRicciuti2018,michelCFM2011}).
The scaling of the constant $\xi(h)$ is chosen such that the 
limit $h\rightarrow 0$ of (\ref{limit-consists})
exists. Clearly the continuous-time limit in (\ref{limit-consists}) exists if and only if
$\lim_{h\rightarrow 0} \xi(h)(1-e^{-hD_t})^{-\alpha} $ exists and hence $\xi(h)=\xi_0 h^{\alpha}$ is the required scaling where $\xi_0 >0$
is an arbitrary positive dimensional constant of physical dimension $\sec^{-\alpha}$ and independent of $h$.
We notice that (\ref{limit-consists}) indeed is a distributional representation of 
Prabhakar-Mittag-Leffler density (\ref{Prabdensity}) which follows in view of Laplace transform of (\ref{limit-consists}), namely 
${\tilde  \chi}_{\alpha,\nu}(s) = \frac{\xi_0^{\nu}}{(\xi_0+s^{\alpha})^{\nu}}$ coinciding with Laplace transform (\ref{GFPPcontinuous}) of the Prabhakar density.
To obtain the limiting density explicitly we
introduce the rescaled variable $\tau_k =hk$ kept finite for $h\rightarrow 0$. 
Hence in (\ref{probabilitiesDFPP}) $k=\frac{\tau_k}{h} \in \mathbb{N}$ becomes very large for $h\rightarrow 0$ thus we can use the asymptotic expression 
$\frac{(\beta)_k}{k!} =\frac{\Gamma(\beta+k)}{\Gamma(k+1)\Gamma(\beta)} \sim \frac{k^{\beta-1}}{\Gamma(\beta)}$ holding for $k$ large.
 Then consider the scaling behavior of the coefficients in (\ref{probabilitiesDFPP}) in the expansions converging for $0<\xi<1$
(as $\xi(h)=\xi_0h^{\alpha} \rightarrow 0$), namely
\begin{equation}
\label{scalinrelations}
\begin{array}{l} \ds 
\frac{((m+\nu)\alpha)_{k}}{k!}  (\xi(h))^{m+\nu} \sim  (\xi_0h^{\alpha})^{m+\nu} \frac{k^{\alpha (m+\nu)-1}}{\Gamma(\alpha (m+\nu))} = 
h \xi_0^{m+\nu} \frac{\tau_k^{\alpha (m+\nu)-1}}{\Gamma(\alpha (m+\nu))} ,\\ \\
\ds \frac{((m+\nu-1)\alpha)_{k}}{k!} (\xi(h))^{m+\nu} \sim  h^{1+\alpha} \xi_0^{m+\nu} \frac{\tau_k^{\alpha(m+\nu-1)-1}}{\Gamma(\alpha(m+\nu-1))} .
\end{array}
\end{equation}
It follows from (\ref{densityformulation}) 
(and see also (\ref{thepropertyogdelta})-(\ref{limitexist}))
that multiplying (\ref{scalinrelations}) by $h^{-1}$ yields densities which remain finite in the continuous-time limit $h\rightarrow 0$.
Another important thing here is that the second coefficient tends to zero by a factor 
$h^{\alpha}$ faster than the first one. 
Generally we observe that
terms of the form $\xi^{\mu}(1-u)^{-\alpha\mu+\lambda} \rightarrow 0$ ($\lambda >0$)  
giving rise to terms scaling as $\sim h^{\lambda} \rightarrow 0$ (where $\tau_k=hk$ and $\xi(h)=\xi_0h^{\alpha}$), namely (See Appendix \ref{appendix1})

\begin{equation}
\label{terms}
\frac{1}{h}(\xi(h))^{\mu} \frac{(\alpha\mu-\lambda)_{k}}{k!} \sim 
\xi_0^{\mu} h^{\alpha\mu-1}  \frac{k^{\alpha\mu-\lambda-1}}{\Gamma(\alpha\mu-\lambda)} =
 h^{\lambda} \frac{\tau_{k}^{\alpha\mu-\lambda-1}}{\Gamma(\alpha\mu-\lambda)} \sim h^{\lambda} 
 \frac{\tau^{\alpha\mu-\lambda-1}}{\Gamma(\alpha\mu-\lambda)} \sim h^{\lambda} \rightarrow 0
\end{equation}
vanishing in the continuous-time limit $h\rightarrow 0$. Hence, for the continuous-time limit, only the part 
$\varphi_{\alpha}^{(\nu )}(k)$ is relevant as $1-(1-u)^{\alpha} \rightarrow 1$.
Then 
consider (\ref{densityformulation}) in the limit $h\to 0$ by using 
(\ref{scalinrelations}) to arrive at
\begin{equation}
\label{densityformulation2}
\chi_{\alpha,\nu }(t)_{ct} = \lim_{h\rightarrow 0} \chi_{\alpha,\nu}(t)_h = \lim_{h\rightarrow 0}
\sum_{k=1}^{\infty} {\cal \psi}_{\alpha}^{(\nu)}(k) \delta_h(t-kh) = 
 \lim_{h\rightarrow 0} \frac{1}{h} {\cal \psi}_{\alpha}^{(\nu)}\left(\frac{t}{h},\xi_0h^{\alpha}\right)
\end{equation}
with
\begin{equation}
\label{asymptsca}
\begin{array}{l} 
\ds 
\chi_{\alpha,\nu }(t)_{ct} = 
 \lim_{h\rightarrow 0}\frac{1}{h} {\cal \psi}_{\alpha}^{(\nu)}\left(\frac{t}{h},\xi_0h^{\alpha}\right) \\ \\ \hspace{1.5cm}\ds = \lim_{h\rightarrow 0}
\sum_{m=0}^{\infty} \frac{(-1)^{m} (\nu )_{m} }{m!}\left( 
 \xi_0^{m+\nu} \frac{t^{\alpha (m+\nu)-1}}{\Gamma(\alpha (m+\nu))} - h^{\alpha} \xi_0^{m+\nu} \frac{t^{\alpha(m+\nu-1)-1}}{\Gamma(\alpha(m+\nu-1))}\right)   \\ \\
 \hspace{1.5cm}\ds  = \sum_{m=0}^{\infty} \frac{(-1)^{m} (\nu )_{m} }{m!}
 \xi_0^{m+\nu} \frac{t^{\alpha (m+\nu)-1}}{\Gamma(\alpha (m+\nu))} ,\hspace{1cm} t \in \mathbb{R}^{+} .
 \end{array}
\end{equation}
Throughout this paper we commonly use subscript notation $(..)_{ct}$ for continuous-time limit distributions.
We can also obtain this result from (\ref{densityformulation2}) with $h=\tau_{k+1}-\tau_k \rightarrow {\rm d}\tau$ and $\tau_k=hk \rightarrow \tau$ and by using the limiting property $\delta_h(t-kh) \rightarrow \delta(t-\tau)$ (see again (\ref{thepropertyogdelta})-(\ref{contilimim-cmul})). 
Hence we can also write (\ref{densityformulation2}) in the form 
\begin{equation}
\label{limit}
\chi_{\alpha,\nu }(t)_{ct} = \lim_{h \to 0} \int_0^{\infty} {\rm d}\tau \delta(t-\tau) \sum_{m=0}^{\infty} \frac{(-1)^{m} (\nu )_{m} }{m!}\left(\frac{\xi_0^{m+\nu}\tau^{\alpha (m+\nu)-1}}{\Gamma((m+\nu)\alpha)} + h^{\alpha} \frac{\tau^{\alpha(m+\nu-1)-1}}{\Gamma(\alpha(m+\nu-1))}\right) .
\end{equation}
The second term tends to zero as $h^{\alpha}$ thus we recover for $h\rightarrow 0$ the density (\ref{asymptsca}), namely
\begin{equation}
\label{densityofGFPP}
\chi_{\alpha,\nu}(t)_{ct} = \xi_0^{\nu}t^{\alpha\nu-1} \sum_{m=0}^{\infty} \frac{(\nu )_{m} }{m!}\frac{(-\xi_0t^{\alpha})^m}{\Gamma((m+\nu)\alpha)} = 
\xi_0^{\nu}t^{\alpha\nu-1}E_{\alpha,\alpha\nu}^{\nu}(-\xi_0t^{\alpha}) 
,\hspace{0.5cm} 0<\alpha\leq 1, \hspace{0.5cm} \nu >0
\end{equation}
which indeed is the Prabhakar-Mittag-Leffler density (\ref{Prabdensity}). In this way we have shown that 
the waiting-time probabilities (\ref{probabilitiesDFPP}) are discrete-time approximations 
of the Prabhakar-Mittag-Leffler density (\ref{Prabdensity}), and the underlying discrete-time counting process indeed is a 
discrete-time version of the GFPP.
\\[1ex]
\begin{figure*}[!t]
\begin{center}
\includegraphics*[width=1.0\textwidth]{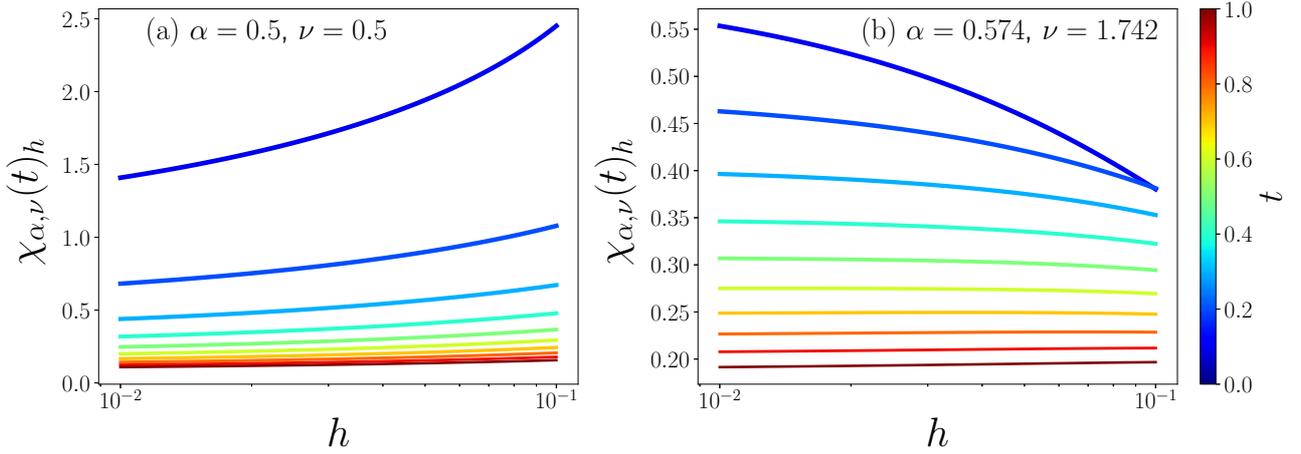}
\end{center}
\vspace{-5mm}
\caption{\label{Fig_2}\it The plots show discrete-time waiting-time density $\chi_{\alpha,\nu}(t)_{h}$ versus $h$ for different times $t$. (a) $\alpha=0.5$ and $\nu=0.5$, (b) $\alpha=0.574$ and $\nu=1.742$. 
The colorbar represents the values of time $t$. The smaller $h$ the more the GFPP Prabhakar waiting time density (\ref{Prabdensity}) is approached (See Eq. (\ref{asymptsca})).
}
\end{figure*}
In Figure \ref{Fig_2} is depicted the behavior of the well-scaled discrete-time density (\ref{densityformulation}) versus 
$h$ for different values of time $t$. For $h\rightarrow 0$ the values converge to the continuous-time Prabhakar density (\ref{Prabdensity}) 
(See (\ref{densityformulation2}), (\ref{asymptsca})). These plots show for fixed $h$ 
monotonically decreasing values of $\chi_{\alpha,\nu}(t)_h$ with time $t$.
This monotonic behavior for the values $\nu$ and $\alpha$ used in these plots reflects the complete monotonicity of the continuous-time limit: The Prabhakar density (\ref{Prabdensity}) is completely monotonic for $\alpha\nu \leq 1$ with $0< \alpha \leq 1$ (See \cite{GiustiPolitoMainardi-etal2020,MainardiGarrappa2015} and references therein for a discussion of complete monotonicity of Prabhakar functions).
\\[1ex]
Now we can define the class of {\it Prabhakar discrete-time processes}: We call a discrete-time
renewal process with waiting times following the distribution $\mathbb{P}(Z=t)= \chi(t)_1$ ($t \in \mathbb{N}_0$) 
{\it Prabhakar} if a well-scaled continuous-time limit $\lim_{h\rightarrow 0} \chi(t)_h= \chi_{\alpha,\nu}(t)_{ct}$ exists (in the sense of (\ref{thepropertyogdelta})-(\ref{contilimitexistsscaling})) to the Prabhakar-Mittag-Leffler density (\ref{Prabdensity}). The remaining part of the paper
is devoted to the analysis of {\it Prabhakar-discrete time processes} and related stochastic motions.

\subsection{\small GENERALIZATION}
\label{generalization}
From the above introduced limiting procedures we can infer that further discrete-time generalizations of the Prabhakar process can be obtained
by the following class of generating functions 
\begin{equation}
\label{genfugen}
{\bar 
\theta}_{\alpha}^{(\nu)} (u) = {\bar f}(u) {\bar \varphi}_{\alpha}^{(\nu)}(u) ,
\hspace{0.5cm}{\bar \varphi}_{\alpha}^{(\nu)}(u)= \frac{\xi^{\nu}}{(\xi+(1-u)^{\alpha})^{\nu}} ,\hspace{0.3cm} \nu>0, \hspace{0.3cm} \alpha \in (0,1] .
\end{equation}
In this expression 
\begin{equation}
\label{hgen}
{\bar f}(u) = \mathbb{E} u^t = u \sum_{k=1}^{\infty} f(k)u^{k-1}
\end{equation} 
which can be seen as generating function of any waiting time-distribution, (i.e.\ with ${\bar f}(u)|_{u=1}=1$) which fulfills
the desired initial condition
${\bar f}(u)|_{u=0} = f(t)\big|_{t=0} =0$ (with $f(t) \geq 0$ on $t\in \mathbb{N}$), we call such a distribution with zero initial condition here simply `$f$-distribution'. It is important to notice that ${\bar f}(u)$ does not contain any scaling parameter $\xi$, thus in the continuous-time limit ${\bar f} \rightarrow 1$. We will see a little later that the continuous-time limit is 
uniquely governed by the `relevant part' ${\bar \varphi}_{\alpha}^{(\nu)}(u)$ in (\ref{genfugen}).
We can generate (\ref{genfugen}) by the following composition procedure 
leading to Eq. (\ref{genfucombi}), see also \cite{PachonPolitoRicciuti2018}. Consider the strictly increasing integer time random variable
\begin{equation}
\label{processI}
\begin{array}{l}
\ds 
 X_n=X_1 ,\hspace{1cm} n=1 \\ \\ \ds
X_n = X_1 + \sum_{j=2}^n Z_j ,\hspace{0.5cm} X_1, Z_j \in \mathbb{N} ,\hspace{0.5cm} n>1
\end{array}
\end{equation}
and $X_0=0$.
The first step $X_1$ is a strictly positive random integer following an $f$-distribution
thus $\mathbb{E} u^{X_1} ={\bar f}(u)$ whereas the increments
$Z_j$ for $j=2,\ldots n$ are IID {\it Sibuya}$(\alpha)$ with $\mathbb{E} u^Z = 1-(1-u)^{\alpha}$ (See Appendix \ref{Sib} for details). The integer random variable (\ref{processI})
has generating function
\begin{equation}
\label{generpro1}
\mathbb{E} u^{X_n} =({\bar w}_I(u))^{(n)} = \mathbb{E}u^{X_1} (\mathbb{E} u^Z)^{n-1} =
{\bar f}(u) \left(1-(1-u)^{\alpha} \right)^{n-1}.
\end{equation}
For the distribution of the events $n$ in (\ref{processI}) 
we choose again (\ref{WIIfracB}).
The so defined process has then generating function
\begin{equation}
\label{hasthenthegen}
\begin{array}{l}
\ds 
{\bar {\cal W}}_{\alpha}^{\nu}(u)= \sum_{n=1}^{\infty} w_B^{(\delta)}(n) ({\bar w}_I(u))^{(n)}
= \sum_{n=1}^{\infty} w_B^{(\delta)}(n){\bar f}(u)  \left(1-(1-u)^{\alpha} \right)^{n-1} \\ \\ \ds \hspace{1.2cm} = \frac{{\bar f}(u)}{1-(1-u)^{\alpha}}\sum_{n=1}^{\infty} w_B^{(\delta)}(n)\left(1-(1-u)^{\alpha} \right)^n \\ \\ \ds
\hspace{1.2cm} = \frac{{\bar f}(u)}{1-(1-u)^{\alpha}}{\tilde w}_B^{(\nu)}\left(1-(1-u)^{\alpha}\right) =
{\bar f}(u) \frac{(\frac{p}{q})^{\nu}}{\left(\frac{p}{q}+(1-u)^{\alpha}\right)^{\nu}} .  
\end{array}
\end{equation}
In the last line we utilized generating functions (\ref{dfgen}) and  (\ref{generalizedfrac}). Clearly (\ref{hasthenthegen}) is a waiting-time generating function of type (\ref{genfugen}).
Now it is only a small step to prove that ${\bar {\cal W}}_{\alpha}^{\nu}(u)$ 
converges to the Prabhakar density for $h\rightarrow 0$. 
Generating function ${\bar f}(u)$ ($|u| \in [0,1]$) can be written as follows
\begin{equation}
\label{as-follows}
{\bar f}(u)= \sum_{t=1}^{\infty}f(t)(1-(1-u))^{t} = 1+\sum_{k=1}^{\infty}g_k(1-u)^k ,\hspace{0.5cm} g_k = (-1)^k
\sum_{t=1}^{\infty} \left(\begin{array}{l} t \\ k \end{array}\right) f(t). 
\end{equation} 
We confirm by ${\bar f}(u)|_{u=1} =\sum_{t=1}^{\infty}f(t)=1$ the normalization of the $f$-distribution. Further, we have ${\bar f}(u)|_{u=0}=\sum_{k=0}^{\infty}g_k = \sum_{t=1}^{\infty}f(t)(1-1)^t=0$.
Clearly, in view of the property (\ref{terms}) it follows that terms $(1-u)^{k}$ produce contributions that tend for $h\rightarrow 0$ to zero 
as $h^k$, namely with $u=e^{-hs}$ 
we have with (\ref{as-follows}) ${\bar f}(e^{-hs}) \sim 1 +\sum_{k=1}^{\infty}g_kh^ks^k
\sim 1+O(h) \to 1$.
\section{\small DISCRETE-TIME VERSIONS OF PRABHAKAR-GENERALIZED POISSON DISTRIBUTION}
\label{DT}
In this section our goal is to analyze a particular important case of Prabhakar discrete-time counting process. We define this process by the strictly increasing random walk
\begin{equation}
\label{renewal-chain}
J_n = \sum_{n=1}^n Z_j ,\hspace{1cm} Z_j \in \mathbb{N} ,\hspace{1cm} J_0=0 ,\hspace{1cm} Z_j \in \mathbb{N}
\end{equation}
where the $Z_j$ are the IID copies of $Z$ 
(interpreted as waiting time in the related counting process) following a Prabhakar type discrete-time distribution 
$\mathbb{P}(Z=k)=\theta_{\alpha}^{(\nu)}(t)$ with generating function of type (\ref{genfugen}). 
The analysis to follow can be extended to any $f$-distribution of the previous section.
As a proto-example we analyze here the most simple 
generating function of this type namely with ${\bar f}(u)=u$, thus
\begin{equation}
\label{mostsimplegen}
\mathbb{E} u^Z =
{\bar \theta}_{\alpha}^{(\nu)}(u) = \frac{\xi^{\nu} u}{(\xi+(1-u)^{\alpha})^{\nu}}.
\end{equation}
For $\nu=1$ ($0<\alpha<1$) (\ref{mostsimplegen}) recovers generating function of the {\it fractional Bernoulli counting process (of `type B')} introduced in \cite{PachonPolitoRicciuti2018} (Eq. (78) therein). 
We call the discrete-time counting process with waiting-time generating function (\ref{mostsimplegen}) the `{\it Prabhakar discrete time process}' ({\it PDTP}).
The PDTP stands out by generalizing fractional Bernoulli (type B), and
for $\nu=1$, $\alpha=1$ (\ref{mostsimplegen}) recovers the generating function
${\bar \theta}_{1}^{(1)}(u) = \frac{p u}{1-qu}$ ($\xi=\frac{p}{q}$ and $p+q=1$) 
of the standard {\it Bernoulli-process}. For $\alpha=1$ and
$\nu>0$, (\ref{mostsimplegen}) coincides with generating function (\ref{dfgen}). 
\\[1ex]
The goal is now to derive explicitly the state probabilities of the PDTP and to show
that the PDTP converges to the continuous-time Prabhakar renewal process (GFPP) under suitable scaling assumptions. 
Note also that the PDTP waiting time distribution has the convenient property that it is the distribution 
$\{\varphi_{\alpha}^{(\nu)}(t) \}$ just shifted 
by one time unit into positive time-direction. This shows the (shift)-operator representation
\begin{equation}
\label{simply}
\begin{array}{l}
\theta_{\alpha}^{(\nu)}(t,\xi) = {\bar \theta}_{\alpha}^{(\nu)}({\hat T}_{-1})\delta_{t0} = 
{\hat T}_{-1} \frac{\xi^{\nu} }{(\xi+(1-{\hat T}_{-1})^{\alpha})^{\nu}}\delta_{t0} = 
{\hat T}_{-1} \varphi_{\alpha}^{(\nu)}(t,\xi) 
\\ \\ \ds 
\hspace{1cm} = \varphi_{\alpha}^{(\nu)}(t-1,\xi) = \frac{1}{(t-1)!}\frac{d^{t-1}}{du^{t-1}}{\bar \varphi}_{\alpha}^{(\nu)}(u)|_{u=0}
,\hspace{0.5cm}{\bar \varphi}_{\alpha}^{(\nu)}(u)= \frac{\xi^{\nu}}{(\xi+(1-u)^{\alpha})^{\nu}} 
\end{array}
\end{equation}
where we always utilize causality, i.e.\ all distributions vanish for negative times. 
The discrete-time density $\varphi_{\alpha}^{(\nu)}(t,\xi)$ ($t\in \mathbb{N}_0$)
is evaluated explicitly in Eq.~(\ref{probabilitiesDFPP}).
On the other hand relation (\ref{simply}) is simply 
reconfirmed by the Leibniz-rule
\begin{equation}
\label{leibniz}
\begin{array}{l} \ds 
\theta_{\alpha}^{(\nu)}(t,\xi) = \frac{1}{t!}\frac{d^{t}}{du^t}\left(u {\bar \varphi}_{\alpha}^{(\nu)}(u)\right)\Big|_{u=0}   \\ \\ \ds 
\hspace{0.5cm} = \frac{1}{t!}\left(u\frac{d^{t}}{du^t}
{\bar \varphi}_{\alpha}^{(\nu)}(u) +t \frac{d^{t-1}}{du^{t-1}} {\bar \varphi}_{\alpha}^{(\nu)}(u)+0\right)\Big|_{u=0}  =
\frac{1}{(t-1)!}\frac{d^{t-1}}{du^{t-1}}{\bar \varphi}_{\alpha}^{(\nu)}(u)\Big|_{u=0} 
\end{array} t\in \mathbb{N}_0 .
\end{equation}
Let us first derive further related distributions such as survival and state probabilities.
To this end consider the probability for {\it at least} one arrival within $[0,t]$, namely
\begin{equation}
\label{at-leastone}
\begin{array}{l} \ds 
\Psi_{\alpha}^{(\nu)}(t,\xi)= \sum_{k=1}^t \theta_{\alpha}^{(\nu)}(k) = \sum_{k=1}^t \varphi_{\alpha}^{(\nu)}(k-1) = \frac{1}{t!} \frac{d^{t}}{du^{t}}
\left(\frac{u {\bar \varphi}_{\alpha}^{(\nu)}(u)}{(1-u)}\right)\Big|_{u=0}  \\ \\ \ds \hspace{1.3cm} =
\frac{1}{(t-1)!}\frac{d^{t-1}}{du^{t-1}} \left(\frac{{\bar \varphi}_{\alpha}^{(\nu)}(u)}{(1-u)}\right)\Big|_{u=0} \, ,\hspace{1cm} t \in \mathbb{N}
\end{array}
\end{equation}
with $\Psi_{\alpha}^{(\nu)}(t,\xi)\big|_{t=0}=0$ since $\theta_{\alpha}^{(\nu)}(t,\xi)\big|_{t=0}=0$ where the generating function of $\Psi_{\alpha}^{(\nu)}(t,\xi)$ is
\begin{equation}
\label{genfuofpsi}
{\bar \Psi}_{\alpha}^{(\nu)}(u) =\sum_{t=0}^{\infty} u^t  \Psi_{\alpha}^{(\nu)}(t) = \frac{u {\bar \varphi}_{\alpha}^{(\nu)}(u)}{(1-u)} .
\end{equation}
Then the survival probability $\Phi^{(0)}_{\alpha,\nu}(t,\xi)$ is
\begin{equation}
\label{survivalprob}
\Phi^{(0)}_{\alpha,\nu}(t,\xi)= \mathbb{P}(J_1>t) = 1-\Psi_{\alpha}^{\nu}(t) = \sum_{k=t+1}^{\infty} \theta_{\alpha}^{(\nu)}(k) =  \sum_{k=t+1}^{\infty} \varphi_{\alpha}^{(\nu)}(k-1)
\end{equation}
with
the generating function
\begin{equation}
\label{wtimelarger}
{\bar \Phi}^{(0)}_{\alpha,\nu}(u) = \sum_{t=0}^{\infty}u^t (1-\Psi_{\alpha}^{\nu}(t))  = \frac{1-u {\bar \varphi}_{\alpha}^{(\nu)}(u)}{(1-u)} , \hspace{1cm} |u| < 1 
\end{equation} 
fulfilling the desired initial condition ${\bar \Phi}^{(0)}_{\alpha,\nu}(u)\big|_{u=0} = 
\Phi^{(0)}_{\alpha,\nu}(t,\xi)\big|_{t=0}=1$ saying that the waiting time $J_1$ for the first arrival is strictly positive.
Then by simple conditioning arguments we obtain the generating function ${\bar \Phi}^{(n)}_{\alpha,\nu}(u)$ of the 
state probabilities $\Phi^{(n)}_{\alpha,\nu}(t,\xi)$ ($n, t \in \mathbb{N}_0$), i.e. the probabilities 
for $n$ arrivals within $[0,t]$ as
\begin{equation}
\label{genfustate}
\begin{array}{l} 
\ds 
{\bar \Phi}^{(n)}_{\alpha,\nu}(u) = {\bar \Phi}^{(0)}_{\alpha,\nu}(u) 
\left( u {\bar \varphi}_{\alpha}^{(\nu)}(u)\right)^n = \frac{(1-u{\bar \varphi}_{\alpha}^{(\nu)}(u))}{(1-u)} u^n 
{\bar \varphi}_{\alpha}^{(n \nu)}(u)    ,\hspace{0.5cm} n \in \{0,1,2,\ldots \} \\ \\ \ds
{\bar \Phi}^{(n)}_{\alpha,\nu}(u) = \frac{u^n{\bar \varphi}_{\alpha}^{(n\nu)}(u)}{1-u} - \frac{u^{n+1}{\bar \varphi}_{\alpha}^{((n+1)\nu)}(u)}{1-u} 
\end{array} \hspace{-0.2cm} |u| < 1,
\end{equation} 
where $ ({\bar \varphi}_{\alpha}^{(\nu)}(u))^n 
 = \frac{\xi^{n\nu}}{(\xi+(1-u)^{\alpha})^{n\nu}} = {\bar \varphi}_{\alpha}^{(n \nu)}(u)$. We also mention the normalization of the state probabilities which can be seen by means of
the general relation
\begin{equation}\label{normastate}
 \frac{1}{t!}\frac{d^t}{du^t}\left\{\sum_{n=0}^{\infty} {\bar \Phi}^{(n)}_{\alpha,\nu}(u)\right\}\Big|_{u=0} = \sum_{n=0}^{\infty} \Phi^{(n)}_{\alpha,\nu}(t,\xi) =\frac{1}{t!}\frac{d^t}{du^t}\frac{1}{1-u}\Big|_{u=0} = 1 ,\hspace{0.5cm} t \in \mathbb{N}_0 .
\end{equation}
Note that relation (\ref{genfustate}) includes $n=0$ where ${\bar \varphi}_{\alpha}^{0}(u) =1$ which has a distribution of the form of a discrete-time $\delta$-distribution (See (\ref{discrete-time-delta}) with $h=1$)
\begin{equation}
\label{0dis}
\varphi_{\alpha}^{(0)}(t,\xi)  = \delta_1(t) = \delta_{t0} 
\end{equation}
thus $\theta_{\alpha}^{(0)}(t,\xi)={\hat T}_{-1}\delta_1(t)=\delta_{t1}$ (See also Eq. (\ref{simply})).
The `state-probabilities' are then obtained from generating function (\ref{genfustate}) as
\begin{equation}
\label{withtheprob} 
\Phi_{\alpha,\nu}^{(n)}(t,\xi) =
\mathbb{P}(J_n \leq t) = \frac{1}{t!}\frac{d^t}{du^t} {\bar \Phi}^{(n)}_{\alpha,\nu}(u)\Big|_{u=0} 
,\hspace{1cm} t, n \in \mathbb{N}_0 .
\end{equation} 
The representation (\ref{genfustate}) is especially convenient for an explicit evaluation of 
$\Phi^{(n)}_{\alpha,\nu}(t,\xi)$.
Be reminded that the state probabilities are shifted distributions where we 
account for (\ref{0dis}) to arrive at
\begin{equation}
\label{shiftet-distr}
\begin{array}{l}
\ds 
 \frac{{\bar \varphi}_{\alpha}^{(n\nu)}({\hat T}_{-1})}{1-{\hat T}_{-1}}{\hat T}_{-n} 
 \varphi_{\alpha}^{(0)}(t) = \frac{{\bar \varphi}_{\alpha}^{(n\nu)}({\hat T}_{-1})}{1-{\hat T}_{-1}} \delta_{1}(t-n) \\ \\ \hspace{3.7cm}= {\hat T}_{-n} {\cal P}_{\alpha,\nu}^{(n)}(t,\xi) =
{\cal P}_{\alpha,\nu}^{(n)}(t-n,\xi) ,\hspace{0.5cm} {\cal P}_{\alpha,\nu}^{(n)}(k,\xi) = {\cal P}_{\alpha,\nu n}(k,\xi) 
\\ \\ \ds  
 {\cal P}_{\alpha,\mu}(k,\xi) = 
\frac{1}{k!} \frac{d^k}{du^k}\frac{{\bar \varphi}_{\alpha}^{(\mu)}(u)}{1-u}\big|_{u=0} ,\hspace{0.5cm} 
{\bar \varphi}_{\alpha}^{(\mu)}(u)= \frac{\xi^{\mu} }{(\xi+(1-u)^{\alpha})^{\mu}} .
\end{array}
\end{equation}
This result is also obtained from the Leibniz-rule which yields 
\begin{equation}
\label{Pprobab}
\begin{array}{l} \ds
{\cal P}_{\alpha,n \nu}(t-n,\xi) = \frac{1}{t!}\frac{d^t}{du^t}\left(\frac{u^n{\bar \varphi}_{\alpha}^{(n \nu)}(u)}{1-u}\right)\Big|_{u=0} \\ \\ \ds
\hspace{2.1cm} = \frac{1}{t!}\sum_{k=0}^t 
\left(\begin{array}{l} t \\ k \end{array}\right) 
\frac{d^k}{du^k}u^n\Big|_{u=0}\frac{d^{t-k}}{du^{t-k}}\left(\frac{{\bar \varphi}_{\alpha}^{(n \nu)}(u)}{1-u}\right)\Big|_{u=0}
\\ \\ \ds \hspace{2.1cm}=  \frac{1}{t!} \frac{t!}{(t-n)!n!} \frac{d^n}{du^n}(u^n) \left( \frac{d^{t-n}}{du^{t-n}}   
\frac{{\bar \varphi}_{\alpha}^{(n \nu)}(u)}{1-u} \right) \Big|_{u=0}
\\ \\ \ds \hspace{2.1cm} =\frac{1}{(t-n)!} \frac{d^{t-n}}{du^{t-n}} \left(  
\frac{{\bar \varphi}_{\alpha}^{(n \nu)}(u)}{1-u} \right) \Big|_{u=0} .
\end{array}
\end{equation}
We hence can write for the state-probability distribution
\begin{equation}
\label{states}
\Phi^{(n)}_{\alpha,\nu}(t,\xi) = {\cal P}_{\alpha,\nu n}(t-n,\xi) - {\cal P}_{\alpha,\nu (n+1)}(t-n-1,\xi) , \hspace{0.5cm} n, t \in \mathbb{N}_0 .
\end{equation} 
To evaluate this expression 
we account for the expansions with respect to $(1-u)^{-1}\xi^{\frac{1}{\alpha}}$, namely
\begin{equation}
\label{takeinto}
\frac{{\bar \varphi}_{\alpha}^{(\mu)}(u)}{1-u} = 
\left\{\begin{array}{clc} \ds 
\sum_{m=0}^{\infty} \left(\begin{array}{l} -\mu  \\ m \end{array}\right) 
\xi^{m+\mu}(1-u)^{-\alpha(m+\mu)-1} , & \xi^{\frac{1}{\alpha}}|(1-u)|^{-1} <1  
\\ \\ 
\ds \sum_{m=0}^{\infty} \left(\begin{array}{l} -\mu  \\ m \end{array}\right) 
\xi^{-m}(1-u)^{\alpha m -1} , &  \xi^{\frac{1}{\alpha}}|(1-u)|^{-1} > 1 .\end{array}\right.
\end{equation}
Then we get
\begin{equation}
\label{wegetthen}
{\cal P}_{\alpha,\mu}(k,\xi) = 
\frac{1}{k!} \frac{d^k}{du^k}\left(\frac{{\bar \varphi}_{\alpha}^{(\mu)}(u)}{1-u}\right)\Big|_{u=0} = \left\{\begin{array}{clc} \ds 
\frac{\xi^{\mu}}{k!}\sum_{m=0}^{\infty} \frac{(-1)^m (\mu )_{m} \xi^m}{m!}  \frac{\Gamma[\alpha(m+\mu)+1+k]}{\Gamma[\alpha(m+\mu)+1]} , & 0 < \xi < 1 \\ \\
\ds \frac{(-1)^k}{k!}\sum_{m=0}^{\infty} \frac{(-1)^m (\mu )_{m} \xi^{-m}}{m!}
\frac{\Gamma(\alpha m)}{\Gamma(\alpha m-k)} , & \xi > 1
\end{array}\right.
\end{equation}
with $k \in \mathbb{N}_0$ and $(\rho)_m$ denotes the Pochhammer symbol (\ref{Pochhammer}). 
With relations (\ref{states}) and (\ref{wegetthen}) we can
write the state-probabilities as
\begin{equation}
\label{state-prob}
\begin{array}{lll} \ds 
\Phi^{(n)}_{\alpha,\nu}(t,\xi) &
\\ \\ 
 \ds  =
\frac{\Theta(t-n)\xi^{n\nu}}{(t-n)!}\sum_{m=0}^{\infty} \frac{(-1)^m (n\nu)_{m} \xi^m}{m!}  \frac{\Gamma[\alpha(m+n\nu)+1+t-n]}{\Gamma[\alpha(m+n\nu)+1]} & \\ \\ \ds \hspace{0.5cm} - 
\frac{\Theta(t-n-1)\xi^{(n+1)\nu}}{(t-n-1)!}\sum_{m=0}^{\infty} \frac{(-1)^m ([n+1]\nu)_{m} \xi^m}{m!}  \frac{\Gamma[\alpha[m+(n+1)\nu]+t-n]}{\Gamma(\alpha[m+(n+1)\nu]+1)} ,
 \hspace{1cm}  0<\xi<1 &
\\ \\
\ds = 
\frac{\Theta(t-n)(-1)^{t-n}}{(t-n)!}
\sum_{m=0}^{\infty} \frac{(-1)^m (n\nu )_{m} \xi^{-m}}{m!}
\frac{\Gamma[\alpha m]}{\Gamma[\alpha m-(t-n)]} \\ \\ \hspace{0.5cm} \ds -\frac{\Theta(t-n-1)(-1)^{t-n-1}}{(t-n-1)!}
\sum_{m=0}^{\infty} \frac{(-1)^m ([n+1]\nu )_{m} \xi^{-m}}{m!}
\frac{\Gamma[\alpha m]}{\Gamma[\alpha m-(t-n-1)]}
, \hspace{1cm} \xi > 1
 \\ \\ \\ \ds 
 \Phi^{(n)}_{\alpha,\nu}(t,\xi) = 0 ,\hspace{1cm} t < n  
\end{array} 
\end{equation} 
where $n,t =\{0,1,2,\dots\} \in \mathbb{N}_0$. In this expression the (discrete-time) Heaviside functions (with $\Theta(0)=1$, see (\ref{heavisideintegers})) reflect causality of ${\cal P}_{\alpha,\mu}(k,\xi)$ in (\ref{states}) such that $\Phi^{(n)}_{\alpha,\nu}(t,\xi) = 0$ for $t<n$ and hence fulfills initial condition $\Phi^{(n)}_{\alpha,\nu}(t,\xi)\big|_{t=0} = 0 $ for $n>0$. 
We notice that (\ref{state-prob}) 
is non-zero for $t\geq n$ (starting with $\Phi^{(n)}_{\alpha,\nu}(t,\xi)\big|_{t=n}={\cal P}_{\alpha,n\nu}(k,\xi)\big|_{k=0}=\frac{\xi^{n\nu}}{(1+\xi)
^{n\nu}}$ which gives for $n=0$ the initial condition $\Phi^{(0)}_{\alpha,\nu}(t,\xi)\big|_{t=0}=1$).
Keep in mind that we utilize the synonymous notations $\Phi^{(n)}_{\alpha,\nu}(t) = \Phi^{(n)}_{\alpha,\nu}(t,\xi)$, the latter when it is necessary to consider the dependence of parameter $\xi$ (for instance in the continuous-time limit).
\\[1ex]
It is especially instructive to consider contained special cases in (\ref{state-prob}), namely fractional Bernoulli $\nu=1$ with $0<\alpha<1$ (subsequent Eq. 
(\ref{state-prob-disceet-Laskin})) 
as well as standard Bernoulli $\nu=1$ with $\alpha=1$ (subsequent Eqs. (\ref{B1}), (\ref{B2})).
Consider now the survival probability, i.e.\ $n=0$ in (\ref{state-prob}), namely
\begin{equation}
\label{survivalgfppdis}
\begin{array}{l} \ds 
\Phi^{(0)}_{\alpha,\nu}(t,\xi) = 1 - \Psi_{\alpha}^{(\nu)}(t,\xi) 
,\hspace{0.5cm} t\in  \mathbb{N}_0 \\ \\ \ds 
 \Phi^{(0)}_{\alpha,\nu}(t,\xi)\big|_{t=0} = 1 
\end{array}
\end{equation}
where $\Psi_{\alpha}^{\nu}(t,\xi)$ is the probability of 
{\it at least one event} within $[0,t]$ 
\begin{equation}
\label{psialphanu}
\begin{array}{l} \ds
\Psi_{\alpha}^{(\nu)}(t,\xi) =  \Theta(t-1){\cal P}_{\alpha,\nu}(t-1,\xi) ,\hspace{0.5cm} t \geq 1 \\ \\
\ds \Psi_{\alpha}^{(\nu)}(t,\xi)|_{t=0} = 0 
\end{array}
\end{equation}
and has generating function (\ref{genfuofpsi}). Since
$\Psi_{\alpha}^{\nu}(t,\xi)|_{t=0}=0$ we have for initial condition of the survival probability $\Phi^{(0)}_{\alpha,\nu}(t,\xi)\big|_{t=0}=1$.
Thus we identify for the state-probabilities (\ref{state-prob}) the important initial condition
\begin{equation}
\label{iniconddesriable}
\Phi^{(n)}_{\alpha,\nu}(t,\xi)\big|_{t=0} = 
{\bar \Phi}^{(n)}_{\alpha,\nu}(u)\big|_{u=0}  = \delta_{n0} ,\hspace{1cm} 
n = \{ 0,1,2, \ldots \} \in  \mathbb{N}_0 .
\end{equation}
The initial condition of this form indeed is crucial for many applications of discrete-time 
renewal processes which come along as Cauchy initial-value problems. By this reason we have constructed generating function
(\ref{mostsimplegen}) such that it fulfills initial condition  ${\bar \theta}_{\alpha}^{(\nu)}(u)\big|_{u=0}=0$. We will come back to this important issue later on in the context of `discrete-time random walks' (Section \ref{CTRM-Prabhakar}).
\\[1ex]
In order to verify that the state probabilities (\ref{state-prob}) approximate the continuous-time state 
probabilities of the
Prabhakar process, let us consider the continuous-time limiting process more closely
(Appendices \ref{crucialdefinitions}, especially (\ref{thepropertyogdelta})-(\ref{contilimim-cmul}) for shift-operator properties and general limiting procedures). The continuous-time limit state probabilities are determined by the limiting behavior of the well-scaled
state probabilities in the sense of relation (\ref{cumuldiscrete-time}) 
\begin{equation}
\label{shift-rep-state}
\begin{array}{l} \ds 
\Phi_{\alpha,\nu}^{(n)}(t)_{ct} = \lim_{h\rightarrow 0}  \frac{h}{1-{\hat T}_{-h}}\left\{
{\hat T}_{-nh} {\bar \varphi}_{\alpha}^{(n\nu)}({\hat T}_{-h})  -{\hat T}_{-(n+1)h} 
{\bar \varphi}_{\alpha}^{((n+1)\nu)}({\hat T}_{-h}) \right\} \delta_h(t) \\ \\ \ds =
 D_{t}^{-1} \left( \frac{\xi_0^{n\nu}}{(\xi_0+D_{t}^{\alpha})^{n\nu}} -  
\frac{\xi_0^{(n+1)\nu}}{(\xi_0+D_{t}^{\alpha})^{(n+1)\nu}}\right) \delta(t) .
\end{array}
\end{equation}
Laplace transforming this relation indeed recovers the Laplace transform
of the Prabhakar continuous-time state probabilities (\cite{TMM-APR-PhysicaA2020}, Eq. (36)).
This continuous-time limit is obtained explicitly by performing the well-scaled limit 
(\ref{contilimim-cmul}) in (\ref{state-prob}) by accounting for the fact that the state 
probabilities are dimensionless cumulative distributions, namely
\begin{equation}
\label{theform}
\Phi_{\alpha,\nu}^{(n)}(t)_{ct} = \lim_{h\rightarrow 0} \Phi_{\alpha,\nu}^{(n)}\left(\frac{t}{h},\xi_0h^{\alpha}\right).
\end{equation}
\begin{figure*}[!t]
\begin{center}
\includegraphics*[width=1.0\textwidth]{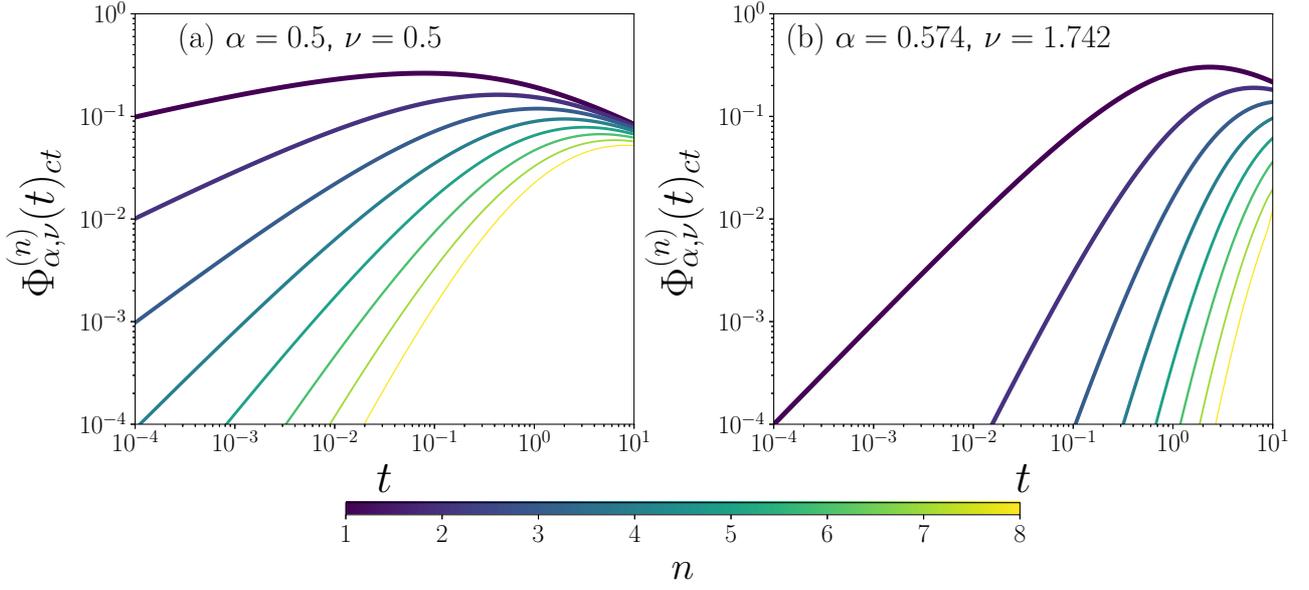}
\end{center}
\vspace{-5mm}
\caption{\label{Fig_3} State probabilities $\Phi_{\alpha,\nu}^{(n)}(t)_{ct}$ of the GFPP representing the continuous-time limit of the PDTP state probabilities versus $t$ for different values of $n$. (a) $\alpha=0.5$ and $\nu=0.5$, (b) $\alpha=0.574$ and $\nu=1.742$. In the colorbar we represent the values $n=1,2,\ldots,8$. The results were obtained numerically using $\xi_0=1$ with Eq.\ (\ref{state-prob-scalintg-limit}).
}
\end{figure*} 
By using then
the asymptotic relation of the
Pochhammer symbol for $k=\frac{t}{h}$ large $\frac{(\mu)_k}{k!} =\frac{\Gamma(\mu+k)}{\Gamma(\mu)\Gamma(k+1)}\sim \frac{k^{\mu-1}}{\Gamma(\mu)}$ in the expansion (\ref{state-prob}) for 
$0<\xi<1$ (as $\xi(h)=\xi_0h^{\alpha} \to 0$) we arrive at
\begin{equation}
\label{state-prob-scalintg-limit}
\begin{array}{l} \ds
\Phi^{(n)}_{\alpha,\nu}(t)_{ct} = (\xi_0t^{\alpha})^{n\nu} \sum_{m=0}^{\infty}\frac{(-\xi_0t^{\alpha})^m}{m!}
\left\{\frac{(n\nu)_m}{\Gamma(\alpha m+\alpha\nu n+1)} - \frac{(\xi_0t^{\alpha})^{\nu}((n+1)\nu)_m}{\Gamma(\alpha m+\alpha\nu(n+1)+1)} 
\right\} \\ \\  \ds
 \hspace{0.5cm} = (\xi_0t^{\alpha})^{n\nu} \left\{ E_{\alpha, \alpha n \nu +1}^{n\nu}(-\xi_0t^{\alpha}) -  
 (\xi_0t^{\alpha})^{\nu} E_{\alpha, \alpha (n+1) \nu +1}^{(n+1)\nu}(-\xi_0t^{\alpha})\right\} ,\hspace{0.1cm} n\in \mathbb{N}_0, 
 \hspace{0.5cm} t \in \mathbb{R}^{+}
 \end{array}
\end{equation}
where in this expression appears the {\it Prabhakar function} $E_{a,b}^c(z)$ (\ref{genmittag-Leff}). Expression (\ref{state-prob-scalintg-limit}) indeed coincides with
the state probabilities of the {\it continuous-time Prabhakar counting process (Generalized Fractional Poisson process - GFPP)} 
\cite{PolitoCahoy2013,TMM-APR-PhysicaA2020} (see Eq. (2.8) in \cite{PolitoCahoy2013} and Eq. (38) in  \cite{TMM-APR-PhysicaA2020}). 
\\[1ex]
In Figure \ref{Fig_3} we draw the GFPP state probabilities (\ref{state-prob-scalintg-limit}) for different states $n$. The state probabilities exhibit for large $t$ an universal power-law limit which is independent of $n$ (See Eq. (\ref{state-asymp})).
The larger $n$ (indicated with brighter colors) the smaller the state probability is for the same time $t$. This general behavior can be understood with the intuitive picture that states with higher $n$ are less `occupied' at the same time $t$.
\\[1ex]
We show in \cite{TMM-APR-PhysicaA2020} that for $\nu=1$ relation (\ref{state-prob-scalintg-limit}) 
recovers Laskin's {\it fractional Poisson distribution} introduced in \cite{Laskin2003} which 
is also the scaling limit of discrete-time state probabilities (\ref{state-prob}) for $\nu=1$.
In order to see explicitly the connection with the fractional Poisson process 
we write the state probabilities (\ref{state-prob}) for $\nu=1$ and $0<\xi < 1$ in the form
\begin{equation}
\label{state-prob-disceet-Laskin}
\begin{array}{l} \ds
\ds 	\Phi^{(n)}_{\alpha,1}(t,\xi) =  
	\xi^n\frac{(\alpha n+1)_{t-n}}{(t-n)!}\Theta(t-n)
	+\sum_{m=1}^{\infty}(-1)^m\xi^{m+n}  \\ \\ \ds 
\times	\left(\Theta(t-n)\frac{(n)_m}{m!}\frac{(\alpha(n+m)+1)_{t-n}}{(t-n)!}+ 
	\Theta(t-n-1)\frac{(n+1)_{m-1}}{(m-1)!} \frac{(\alpha(n+m)+1)_{t-n-1}}{(t-n-1)!} \right) 
	,\hspace{0.1cm}   \\ \\ \ds
		\Phi^{(n)}_{\alpha,1}(t,\xi) = 0 ,\hspace{1cm} t < n 
	\end{array}
\end{equation}
where $n,t =\{0,1,2,\dots\} \in \mathbb{N}_0$.
For $0<\alpha<1$ these are the state-probabilities of the {\it fractional Bernoulli process} (type B), and
for $\alpha=1$ this expression recovers the state probabilities of standard Bernoulli shown subsequently in Eqs. (\ref{B1}), (\ref{B2}).
The state probabilities (\ref{state-prob-disceet-Laskin}) have
with the same limiting procedure 
the continuous-time limit
\begin{equation}
\label{fractionalPoissondistribution}
\begin{array}{l} \ds
\Phi^{(n)}_{\alpha,1}(t)_{ct}= \lim_{h \rightarrow 0} \Phi^{(n)}_{\alpha,1}\left(\frac{t}{h},\xi_0h^{\alpha}\right)    =\frac{(\xi_0t^{\alpha})^n}{\Gamma(\alpha n+1)} 
+ \sum_{m=1}^{\infty}\frac{(-1)^m(\xi_0 t^{\alpha})^{m+n}}{\Gamma(\alpha(m+n)+1)}
\left(\frac{(n)_m}{m!}+\frac{(n+1)_{m-1}}{(m-1)!}\right) \\ \\ \ds 
\hspace{1.4cm} = \frac{(\xi_0 t^{\alpha})^n}{n!} \sum_{m=0}^{\infty}\frac{(n+m)!}{m!}
\frac{(-\xi_0 t^{\alpha})^m}{\Gamma(\alpha(m+n)+1)} ,\hspace{1cm} t \in \mathbb{R}^{+},
\end{array}
\end{equation}
with $\frac{(n)_m}{m!} +\frac{(n+1)_{m-1}}{(m-1)!}= \frac{(n+m)!}{n!m!}$.
We identify (\ref{fractionalPoissondistribution}) indeed with Laskin's {\it fractional 
Poisson distribution} \cite{Laskin2003} which also is recovered for $\nu=1$ 
in the expression (\ref{state-prob-scalintg-limit}) for GFPP state probabilities. It follows that the state probabilities (\ref{state-prob-disceet-Laskin}) of fractional Bernoulli are discrete-time approximations of the fractional Poisson distribution (\ref{fractionalPoissondistribution}).
\\[1ex]
Now let us consider $\alpha=1$, $\nu=1$, i.e. the case of (standard) Bernoulli more closely. 
The generating functions of the state probabilities (\ref{genfustate}) then take the form
\begin{equation}
\label{B1}
{\bar \Phi}_{1,1}^{(n)}(u) = \frac{(\xi+1)(\xi u)^n}{(\xi+1-u)^{n+1}}
,\hspace{1cm} \xi =\frac{p}{q} ,\hspace{0.5cm} p+q=1 ,\hspace{0.5cm} n \in \mathbb{N}_0
\end{equation}
where the state probabilities yield (with $p= \frac{\xi}{\xi+1}$ and  $q =\frac{1}{\xi+1}$)
\begin{equation}
\label{B2}
\begin{array}{l} \ds
\Phi_{1,1}^{(n)}(t,\xi) = \frac{(\xi+1)\xi^n}{(t-n)!} \frac{d^{t-n}}{du^{t-n}}
  \frac{1}{(\xi+1-u)^{n+1}}\big|_{u=0} \\ \\ \ds \hspace{0.5cm}  =\frac{(\xi+1)\xi^n}{(\xi+1)^{n+1+t-n}}
\left(\begin{array}{l} n+t-n \\ t-n \end{array}\right) = \left(\begin{array}{l} t \\ n \end{array}\right) \frac{\xi^n}{(\xi+1)^{t}} = \left(\begin{array}{l} t \\ n \end{array}\right) p^n q^{t-n} ,\hspace{0.1cm} t \geq n \\ \\ \text{and } \Phi_{1,1}^{(n)}(t,\xi)  \ds
= 0  ,\hspace{1cm} t < n
\end{array}  (t, n \in \mathbb{N}_0) .
\end{equation}
We identify (\ref{B2}) with the Binomial distribution (i.e.\ the state distribution of the Bernoulli process).
The continuous-time limit of (\ref{B2}) is obtained from $p=\frac{h\xi_0}{1+h\xi_0}$ and 
$q=\frac{1}{1+h\xi_0}$ thus performing the well-scaled limit (\ref{theform}) yields
\begin{equation}\label{contlim-bern}
\begin{array}{l} \ds 
\Phi^{(n)}_{1,1}(t)_{ct} = \lim_{h\rightarrow 0} \, \frac{1}{n!}
\frac{t}{h}\,\left(\frac{t}{h}-1\right)\, \ldots \,\left(\frac{t}{h}-n+1\right)\, h^n \,\xi_0^n\, (1+h\xi_0)^{-\frac{t}{h}}, \hspace{1cm}  t \in h\mathbb{N}_0 \\ \\ \ds
 \hspace{0.5cm} = \frac{(\xi_0t)^n}{n!}e^{-\xi_0 t} , \hspace{0.5cm} t \in \mathbb{R}^{+} ,\hspace{0.5cm} n \in \mathbb{N}_0
 \end{array}
\end{equation}
which is the {\it Poisson distribution} also recovered in (\ref{fractionalPoissondistribution}) for $\alpha=1$ and 
in
(\ref{state-prob-scalintg-limit}) for $\nu=1$, $\alpha=1$ (Consult also \cite{TMM-APR-PhysicaA2020,MichelitschRiascosGFPP2019}). 
This reflects the well-known fact that the standard Bernoulli process converges in a well-scaled continuous-time limit to standard Poisson (see e.g.\ \cite{PachonPolitoRicciuti2018} and many others).

\subsection{\small ASYMPTOTIC FEATURES}
\label{asymfeat}
In many applications the asymptotic features for large and small observation times are of interest. Clearly the asymptotic behavior for 
large $t$ in a discrete-time distribution is determined by the 
leading power in $(1-u)$ in the limit $u\rightarrow 1-0$ in the generating function.
The generating function (\ref{mostsimplegen}) behaves then as
\begin{equation}
\label{varphiexpan}
\begin{array}{l} \ds
{\bar \theta}_{\alpha}^{(\nu)}(u) = u 
\left(1+\frac{1}{\xi}(1-u)^{\alpha}\right)^{-\nu}  =  u \sum_{m=0}^{\infty}
(-1)^m\frac{(\nu)_m}{m!}\xi^{-m}(1-u)^{\alpha m}  \\ \\ \ds {\bar \theta}_{\alpha}^{(\nu)}(u) \sim {\bar \varphi}_{\alpha}^{(\nu)}(u) \sim \left(1-\frac{\nu}{\xi}(1-u)^{\alpha}+ O(1-u)^{\alpha}\right) , \hspace{0.5cm} (u \to 1-0), \hspace{0.3cm} \alpha \in (0,1),\hspace{0.1cm} \nu>0 .
\end{array}
	\end{equation}
The asymptotic behavior of the waiting time density (\ref{simply}) for $t$ large is hence governed by the term $-\frac{\nu}{\xi}(1-u)^{\alpha}$. We notice that this term is (up to the scaling multiplier $\frac{\nu}{\xi}$) the same as in the generating function of Sibuya $(\alpha)$ (See Appendix \ref{Sib}). 
Thus we get for large $t$ the fat-tailed behavior
\begin{equation}
\label{appro}
\begin{array}{l}
\ds 
\varphi_{\alpha}^{(\nu)}(t,\xi) \sim  \theta_{\alpha}^{(\nu)}(t,\xi) \sim \sum_{m=1}^{\infty}
(-1)^m\frac{(\nu)_m}{m!}\xi^{-m}\frac{(-\alpha m)_t}{t!}   \sim \frac{\nu}{\xi}(-1)^{t-1}\left(\begin{array}{l} \alpha \\ t \end{array}\right)   \sim  
 \frac{\alpha\nu}{\xi} \frac{t^{-\alpha -1}}{\Gamma(1-\alpha)},   \\ \\ \ds
 \hspace{0.5cm}\varphi_{\alpha}^{(\nu)}(t)_{ct} = \theta_{\alpha}^{(\nu)}(t)_{ct}  = \lim_{h\rightarrow 0} \frac{1}{h}\varphi_{\alpha}^{(\nu)}\left(\frac{t}{h},\xi_0h^{\alpha}\right) \sim \frac{\alpha\nu}{\xi_0} \frac{t^{-\alpha -1}}{\Gamma(1-\alpha)} ,\hspace{2cm}  \alpha \in (0,1)
 \end{array} 
\end{equation}
where the last line refers to the large-time asymptotic behavior of the continuous-time limit which is of the same type.
We notice that this asymptotic distribution is (due to the multiplier $\frac{\nu}{\xi}$) a scaled version of {\it Sibuya($\alpha$)}. It follows that the PDTP converges for large observation times to a scaled version of the Sibuya counting process (Appendix \ref{Sibuyarenewal}). This is true for all $\nu>1$ including $\nu=1$ for the {\it fractional Bernoulli} counting process.
The asymptotic relation (\ref{appro}) indeed is in accordance with the well-known fat-tailed asymptotic behavior of the Prabhakar density (\ref{densityofGFPP}) of the GFPP in the continuous-time case 
(See e.g. \cite{TMM-APR-PhysicaA2020}, Eq. (35)). It is clear that the same tail
(\ref{appro}) occurs in a fractional Bernoulli process (with rescaled parameter 
$\xi' = \frac{\xi}{\nu}$) with discrete-time waiting time density being an approximation of the Mittag-Leffler density. Indeed (\ref{appro}) also is the fat-tail of the Mittag-Leffler density reflecting
the asymptotic `Mittag-Leffler universality' \cite{GorenfloMainardi2006}.
\\[1ex]
Let us also consider the asymptotic behavior of the state probabilities. With the same argument we
obtain from the generating function (\ref{genfustate}) the leading power in $(1-u)$ for $u\rightarrow 1-0$,
namely
\begin{equation}
\label{the-leading-power-state-prob}
\begin{array}{l} \ds 
{\bar \Phi}_{\alpha,\nu}^{(n)}(u) = \frac{u^n}{(1-u)}
\left\{\left(1+\frac{(1-u)^{\alpha}}{\xi}\right)^{-n\nu} - u \left(1+\frac{(1-u)^{\alpha}}{\xi}\right)^{-(n+1)\nu}\right\} \\ \\ \ds 
{\bar \Phi}_{\alpha,\nu}^{(n)}(u) \sim {\bar M}^{\alpha,\nu}(u) \sim  \frac{\nu}{\xi}(1-u)^{\alpha-1} ,\hspace{1cm} (u \to 1-0)
\end{array} \alpha \in (0,1).
\end{equation}
For the particular case $\alpha=1$, we obtain the finite value ${\bar \Phi}_{1,\nu}^{(n)}(u)|_{u=1} = {\bar M}^{1,\nu}(u)|_{u=1}=1+\frac{\nu}{\xi}$. In Eq. (\ref{the-leading-power-state-prob}) we have accounted for the behavior of the memory generating function ${\bar M}^{\alpha,\nu}(u)$ introduced in (\ref{casezeo}) for $u \rightarrow 1-0$, with the same
leading term as the state probabilities, namely
${\bar M}^{\alpha,\nu}(u) \sim  \frac{\nu}{\xi}(1-u)^{\alpha-1}$. 
We observe in relation (\ref{the-leading-power-state-prob}) that state probability (and memory-) generating functions $\forall n \in \mathbb{N}_0$
at $u=1$ have a positive {\it weak singularity} $\sim (1-u)^{\alpha-1}$ ($\alpha\in (0,1)$).
This weak singularity indeed is the origin
of the asymptotic large-time power-law decay of the PDTP state probabilities $\forall n \in \mathbb{N}_0$ and of the `memory function' (introduced subsequently in (\ref{memoryop})), namely
\begin{equation}
\label{state-asymp}
\begin{array}{l} \ds
 \Phi_{\alpha,\nu}^{(n)}(t,\xi)
\sim  {\cal M}_{\alpha,\nu}(t)  \sim \frac{\nu}{\xi}(-1)^t
\left(\begin{array}{l} \alpha -1 \\\hspace{0.3cm} t \end{array}\right) \sim \frac{\nu}{\xi} \frac{t^{-\alpha}}{\Gamma(1-\alpha)} , \\ \\ \ds
\Phi_{\alpha,\nu}^{(n)}(t)_{ct} = \lim_{h\rightarrow 0} \Phi_{\alpha,\nu}^{(n)}\left(\frac{t}{h},\xi_0h^{\alpha}\right) \rightarrow
{\cal M}_{\alpha,\nu}(t)_{ct} \sim
\frac{\nu}{\xi_0} \frac{t^{-\alpha}}{\Gamma(1-\alpha)} 
\end{array} \hspace{0.5cm} (t \,\, {\rm large}) \hspace{0.5cm} \alpha \in (0,1) 
\end{equation}
where the last line refers to the continuous-time limit.
The power-law (\ref{state-asymp}) is universal and reflects the long-time memory and non-Markovian feature in the PDTP for $\alpha \in (0,1)$. This asymptotic relation is in accordance with the large-time limiting behavior of the continuous-time GFPP state probabilities \cite{TMM-APR-PhysicaA2020}. Physically two `extreme' regimes are noteworthy: (i) $\alpha\rightarrow 0+$: There we have
$$\ds \lim_{\alpha \rightarrow 0+} \frac{t^{-\alpha}}{\Gamma(1-\alpha)}  = (-1)^t
\left(\begin{array}{l} \alpha -1 \\\hspace{0.3cm} t \end{array}\right)\Big|_{\alpha=0} = \Theta(t)= 1 $$
thus $\Phi_{0+,\nu}^{(n)}(t,\xi) \sim  \frac{\nu}{\xi}$ where extremely long waiting times occur thus the states $n$ `live' extremely long and the memory of the PDTP becomes infinite.
\\[1ex]
In the second `extreme' regime (ii) we have $\alpha\rightarrow 1-0$ where $
\frac{t^{-\alpha}}{\Gamma(1-\alpha)} \sim \delta(t) = 0 $ ($t$ large) indicating lack of memory (for $\nu=1$) or short-time memory (for $\nu \neq 1$).
This feature also is reflected by the non-singular behavior of the state probability generating functions of standard Bernoulli (\ref{B1}) which yield at $u=1$ the finite value
${\bar \Phi}_{1,1}^{(n)}(u)\big|_{u=1}=\frac{\xi+1}{\xi}$ ($\forall n \in \mathbb{N}_0$).
\\[1ex]
We observe in these results the following general property of discrete-time counting processes: If the state probability generating 
functions (and memory generating function) are {\it weakly singular} 
$\sim (1-u)^{\lambda}$ ($\lambda \in (-1,0)$)
at $u=1$, then the discrete-time counting process is non-Markovian and has long-time power-law memory as in (\ref{state-asymp}) with fat-tailed waiting-time density (\ref{appro}).
If in contrast the state probability generating functions (and memory generating function) at $u=1$ are finite, then the 
discrete-time counting process either is memoryless with the Markovian property or is non-Markovian and has only a short-time memory with a light-tailed waiting time density.

\subsection{\small GENERALIZED FRACTIONAL DIFFERENCE EQUATIONS GOVERNING THE PDTP STATE PROBABILITIES}
\label{stae-diffeqs}
The goal of this part is to derive recursive evolution equations 
that are solved by the PDTP state probabilities. To this end we utilize the correspondence of generating functions and their operator representations. Whenever we deal with operator-functions of the shift operator ${\hat T}_{-h}$ we refer to a renewal chain (\ref{renewal-chain}) with rescaled waiting times 
$Z_j \in h\mathbb{N}$ ($h>0$) of the PDTP with generating function (\ref{mostsimplegen}).
The state probability generating functions (\ref{genfustate}) fulfill
\begin{equation}
\label{recursiveequations}
\frac{ {\bar \Phi}^{(n)}_{\alpha,\nu}(u) }{ {\bar \varphi}_{\alpha}^{(\nu)}(u)} = u {\bar \Phi}^{(n-1)}_{\alpha,\nu}(u)  ,\hspace{0.5cm} n \in \mathbb{N}
\end{equation}
and for $n=0$ we have
\begin{equation}
\label{casezeo}
{\bar M}^{\alpha,\nu}(u)=
\frac{{\bar \Phi}^{(0)}_{\alpha,\nu}(u) }{{\bar \varphi}_{\alpha}^{(\nu)}(u)}  
=
\frac{1}{1-u}\left\{\left(1+\frac{(1-u)^{\alpha}}{\xi}\right)^{\nu}-u\right\} .
\end{equation}
where we refer (\ref{casezeo}) to as `memory generating function'.
These simple relations allow us to obtain recursive equations for the state probabilities.
We introduce the operator
\begin{equation}
\label{shiftD}
{\hat {\cal D}}_{\alpha}^{\nu}(h)  = \frac{1}{{\bar \varphi}_{\alpha}^{(\nu)}({\hat T}_{-h})} = \left(1+\frac{(1-{\hat T}_{-h})^{\alpha}}{\xi}\right)^{\nu} ,\hspace{1cm} {\hat {\cal D}}_{\alpha}^{\nu}:={\hat {\cal D}}_{\alpha}^{\nu}(1) =
 \frac{1}{{\bar \varphi}_{\alpha}^{(\nu)}({\hat T}_{-1})}
\end{equation}
where we skip the argument in ${\hat {\cal D}}_{\alpha}^{\nu}={\hat {\cal D}}_{\alpha}^{\nu}(h=1)$ when we refer to integer discrete-time processes with $h=1$.
We utilize as synonymous notations ${\hat A}^{-1} =\frac{1}{{\hat A}}$ and write for the unity operator (zero-shift) simply ${\hat 1}=1$, namely ${\hat T}_a\big|_{a=0}= 1$.
We can then rewrite (\ref{recursiveequations}) and (\ref{casezeo}) compactly in operator form (See also Appendix \ref{crucialdefinitions})
\begin{equation}
\label{recursive-diff}
  {\hat {\cal D}}_{\alpha}^{\nu} \cdot \Phi^{(n)}_{\alpha,\nu}(t) = 
 \Phi^{(n-1)}_{\alpha,\nu}(t-1) +\delta_{n0} {\cal M}_{\alpha,\nu}(t)  ,\hspace{1cm} t, n \in \mathbb{N}_0 .
\end{equation}
Bear in mind the correspondence
$u \leftrightarrow {\hat T}_{-h}$ thus $\frac{1}{u} \leftrightarrow {\hat T}_{+h}$, 
see Appendix \ref{crucialdefinitions} and we use in this formulation
that $\Phi^{(n)}_{\alpha,\nu}(t)=0$ if at least one of the variables $t,n <0$.
(\ref{recursive-diff}) can be conceived as {\it generalized fractional 
Kolmogorov-Feller difference equations} governing the time evolution of the PDTP state probabilities. We
will determine them subsequently in explicit form including their continuous-time limit representations.
For $n=0$ (by accounting for $\Phi^{(-1)}(t)=0$) we have with (\ref{casezeo}) the `{\it memory function}'
\begin{equation}
\label{memoryop}
{\cal M}_{\alpha,\nu}(t)  =  {\hat {\cal D}}_{\alpha}^{\nu} \cdot 
\Phi^{(0)}_{\alpha,\nu}(t) = {\hat M}^{\alpha,\nu}_1\delta_1(t) =
 \frac{1}{1-{\hat T}_{-1}}{\hat {\cal D}}_{\alpha}^{\nu} \delta_{1}(t) - \Theta(t-1) ,\hspace{0.5cm} t \in \mathbb{Z}
\end{equation}
fulfilling initial condition 
${\cal M}_{\alpha,\nu}(t)\big|_{t=0}= M^{\nu,\alpha}(u)\big|_{u=0}= \frac{(\xi+1)^{\nu}}{\xi^{\nu}}$.
We use here that $\frac{h}{1-{\hat T}_{-h}} \delta_h(t) = \sum_{k=0}^{\infty} h\delta_h(t-kh) =\Theta(t)$ thus $\frac{h {\hat T}_{-h} }{1-{\hat T}_{-h}} \delta_h(t) =
\Theta(t-h)$, see Appendix \ref{crucialdefinitions}.
For later use we introduced in (\ref{memoryop}) the `memory operator' which has with $\xi(h)=\xi_0h^{\alpha}$ the well-scaled representation
\begin{equation}
\label{memoryoperator}
{\hat M}^{\alpha,\nu}_h = \frac{h}{1-{\hat T}_{-h}}\left\{{\hat {\cal D}}_{\alpha}^{\nu}(h) 
-{\hat T}_{-h}\right\} 
\end{equation}
which defines by ${\cal M}_{\alpha,\nu}(t)_h = {\hat M}^{\alpha,\nu}_h\delta_h(t)$ ($t\in h\mathbb{Z}$, $\xi=\xi_0h^{\alpha}$) the well-scaled memory function maintaining the initial condition of (\ref{memoryop}). 
We observe that in the case of Bernoulli $\alpha=1$, $\nu=1$ 
the memory function becomes `local', namely ${\hat M}^{1,1}_h\delta_h(t)=\frac{(1+\xi)}{\xi}h\delta_{h}(t)$ 
vanishing for $t>0$ which reflects the loss of memory
in the standard Bernoulli process
(See 
\cite{PachonPolitoRicciuti2018} for general aspects).
Now we rewrite above introduced shift-operator functions (\ref{shiftD}) and (\ref{memoryoperator})
in such a way that only generalized fractional integrals 
and derivatives emerge in the continuous-time limit (involving kernels which are at most {\it weakly singular} and hence integrable), namely
\begin{equation}
\label{Drerite}
\begin{array}{l} \ds 
{\hat {\cal D}}_{\alpha}^{\nu}(h) = \xi_0^{-\nu} h^{-\nu\alpha} (1-{\hat T}_{-h})^{\alpha\nu} 
\left(1+\xi_0h^{\alpha}(1-{\hat T}_{-h})^{-\alpha}\right)^{\nu}
\\ \\ \ds \hspace{1.1cm} = \xi_0^{-\nu} h^{-\ceil{\nu\alpha}} (1-{\hat T}_{-h})^{\ceil{\alpha\nu}}
 h^{-(\alpha\nu-\ceil{\alpha\nu})} (1-{\hat T}_{-h})^{\alpha\nu-\ceil{\alpha\nu}}
\left(1+\xi_0h^{\alpha}(1-{\hat T}_{-h})^{-\alpha}\right)^{\nu} \\ \\ \ds
\hspace{1.1cm} = \xi_0^{-\nu} h^{-\ceil{\nu\alpha}} (1-{\hat T}_{-h})^{\ceil{\alpha\nu}} {\hat B}_{\alpha,\nu}(h) 
\end{array}
\end{equation} 
where we introduced the {\it ceiling function} $\ceil{\mu}$ which indicates the smallest 
integer larger or equal to $\mu$ (See Appendix \ref{appendix1}). 
We emphasize that for the continuous-time limit representations to be deduced 
in this part
we have the scaling $\xi(h)=\xi_0h^{\alpha} \rightarrow 0$ and therefore write expansions which converge for $0<\xi<1$ (See e.g. expansions (\ref{probabilitiesDFPP}) and (\ref{state-prob})) without always explicitly mentioning this issue.
In Appendix \ref{appendix1} (See especially Eqs.\ (\ref{further-important})-(\ref{RLint})) it is shown that the operator $\lim_{h\rightarrow 0} h^{-\nu\alpha} (1-{\hat T}_{-h})^{\alpha\nu}= D_t^{\alpha\nu}$ that occurs in (\ref{Drerite}) converges 
to the Riemann-Liouville fractional derivative of order $\alpha\nu$.
Relation (\ref{Drerite}) contains the operator function
\begin{equation}
\label{Boperator}
{\hat B}_{\alpha,\nu}(h) = h^{-(\alpha\nu-\ceil{\alpha\nu})} (1-{\hat T}_{-h})^{\alpha\nu-\ceil{\alpha\nu}}
\left(1+\xi_0h^{\alpha}(1-{\hat T}_{-h})^{-\alpha}\right)^{\nu} 
\end{equation}
which can be conceived as a discrete-time {\it generalized fractional integral operator}. We call an operator a `{\it generalized fractional integral}' if it contains only non-negative powers of the discrete integral 
operator $h(1-{\hat T}_{-h})^{-1}$, see especially Appendix \ref{appendix1} and consult also \cite{michelCFM2011}. 
In view of (\ref{memoryoperator}) we further introduce
\begin{equation}
\label{theinter}
{\hat K}_{\alpha\nu}^{(0)}(h)=\frac{h}{1-{\hat T}_{-h}}{\hat {\cal D}}_{\alpha}^{\nu}(h) =  \xi_0^{-\nu} h^{1-\ceil{\nu\alpha}} 
(1-{\hat T}_{-h})^{\ceil{\alpha\nu} -1} {\hat B}_{\alpha,\nu}(h) .
\end{equation}
All these operators are operator functions of the shift operator ${\hat T}_{-h}=e^{-hD_t}$ 
and can be considered as discrete-time versions of generalized fractional derivatives or integrals, respectively. 
The discrete-time memory function (\ref{memoryop}) then can be evaluated as (See also Appendix \ref{crucialdefinitions}) 
\begin{equation}
\label{discrete-timememoryf}
\begin{array}{l}
\ds 
{\cal M}_{\alpha,\nu}(t)  = 
\frac{1}{1-{\hat T}_{-1}}{\hat {\cal D}}_{\alpha}^{\nu} \delta_{1}(t) -{\hat T}_{-1}\Theta(t) =
{\cal K}^{(0)}_{\alpha,\nu}(t) -\Theta(t-1) \\ \\
\ds  {\cal K}^{(0)}_{\alpha,\nu}(t) = 
\frac{1}{1-{\hat T}_{-1}}{\hat {\cal D}}_{\alpha}^{\nu} \delta_{1}(t) =  
\xi^{-\nu}\left(1-{\hat T}_{-1}\right)^{\ceil{\alpha\nu}-1} {\cal B}_{\alpha,\nu}(t)  ,\hspace{2cm} t\in \mathbb{Z} \\ \\ 
 \ds  {\cal B}_{\alpha,\nu}(t) = {\hat B}_{\alpha,\nu}(1) \delta_1(t) =
\sum_{m=0}^{\infty} \frac{(-\nu)_m}{m!} (-\xi)^m 
\frac{(\alpha m + \ceil{\alpha\nu} -\alpha \nu)_t}{t!}  ,\hspace{1cm} t\in \mathbb{N}_0 .
 \end{array} 
\end{equation}
Now let us consider 
\begin{equation}
\label{exapandas}
\begin{array}{l} \ds
 {\hat B}_{\alpha,\nu}(h)\Phi(t)= \sum_{m=0}^{\infty} \frac{(-\nu)_m}{m!} (-\xi_0)^m 
 \left(h
 (1-{\hat T}_{-h})^{-1}\right)^{\alpha m + \ceil{\alpha\nu} -\alpha \nu}\Phi(t)  \\ \\ 
 \ds \hspace{2.1cm} = \sum_{m=0}^{\infty} \frac{(-\nu)_m}{m!} (-\xi_0)^m 
 \sum_{k=0}^{\infty}\frac{(\alpha m + \ceil{\alpha\nu} -\alpha \nu)_k}{k!}h^{(\alpha m + \ceil{\alpha\nu} -\alpha \nu)}\Phi(t-kh) 
 \end{array} \hspace{0.5cm} t \in h\mathbb{Z}
\end{equation}
where we put always $\xi = \xi_0 h^{\alpha} < 1$ in order to have well-scaled operators with existing continuous-time limits.
This relation is a generalized fractional difference equation with memory where 
the whole history $\{\Phi(t-kh)\}$ ($t\in h\mathbb{N}_0$) of the causal function $\Phi(t)$ contributes.
Now since $\alpha m + \ceil{\alpha\nu} -\alpha \nu >0$ if $\alpha\nu \notin \mathbb{N}$ and  
$\alpha m + \ceil{\alpha\nu} -\alpha \nu \geq 0$ if $\alpha\nu \in \mathbb{N}$ we can use for
the continuous-time limit that $
\frac{(\alpha m + \ceil{\alpha\nu} -\alpha \nu)_k}{k!}h^{(\alpha m + \ceil{\alpha\nu} -\alpha \nu)} \rightarrow  \frac{\tau^{\alpha m + \ceil{\alpha\nu} -\alpha \nu-1}}{\Gamma(\alpha m + \ceil{\alpha\nu} -\alpha \nu)} {\rm d\tau} $ ($h\rightarrow {\rm d}\tau$ and $\tau=kh$, Appendix \ref{appendix1}) thus
\begin{equation}
\label{wehave}
\begin{array}{l} \ds 
D_t^{-(\alpha m + \ceil{\alpha\nu} -\alpha \nu)}  \Phi(t)_{ct} =
\lim_{h\rightarrow 0} \left(h
 (1-{\hat T}_{-h})^{-1}\right)^{\alpha m + \ceil{\alpha\nu} -\alpha \nu} \Phi(t) \\ \\ \ds \hspace{3.8cm}  =  \int_0^t\Phi(t-\tau)_{ct} 
 \frac{\tau^{( \alpha m + \ceil{\alpha\nu} -\alpha \nu) -1}}{\Gamma(\alpha m + \ceil{\alpha\nu} -\alpha \nu)}{\rm d}\tau
 \end{array}
 \end{equation}
 are Riemann-Liouville fractional integrals of orders $(\alpha m + \ceil{\alpha\nu} -\alpha \nu)$. 
 Note that for $\alpha\nu \in \mathbb{N}$ where  $\ceil{\alpha\nu} -\alpha \nu =0$ the kernel of the order $m=0$ becomes $\lim_{\alpha\nu \rightarrow \ceil{\alpha\nu}}
 \frac{\tau^{(\ceil{\alpha\nu} -\alpha \nu) -1}}{\Gamma(\ceil{\alpha\nu} -\alpha \nu)} =\delta(\tau)$ a Dirac's $\delta$-distribution thus (\ref{wehave}) recovers $\Phi(t)_{ct}$. The continuous-time limit of (\ref{exapandas}) then yields
\begin{equation}
\label{conntilimBh}
\begin{array}{l} 
\ds 
\lim_{h\rightarrow 0} {\hat B}_{\alpha,\nu}(h)\Phi(t) =  \int_0^t{\cal B}_{\alpha,\nu}(\tau)_{ct}\Phi(t-\tau)_{ct}{\rm d}\tau  \\ \\ \ds \hspace{2.8cm}  =
\sum_{m=0}^{\infty} \frac{(-\nu)_m}{m!} (-\xi_0)^m 
 \int_{0}^{\infty} \frac{\tau^{\alpha m + \ceil{\alpha\nu} -\alpha \nu-1}}{\Gamma(\alpha m + \ceil{\alpha\nu} -\alpha \nu)} \Phi(t-\tau)_{ct}{\rm d}\tau
 \end{array} 
\end{equation}
with continuous-time limit kernel
\begin{align}
\label{kernelwehave}
\ds
{\cal B}_{\alpha,\nu}(\tau)_{ct}
&= \begin{cases} \ds \tau^{\ceil{\alpha\nu} -\alpha \nu-1}\sum_{m=0}^{\infty} 
\frac{(-\nu)_m}{m!}
\frac{(-\xi_0\tau^{\alpha})^m}{\Gamma(\alpha m + \ceil{\alpha\nu} -\alpha \nu)}, & \alpha\nu \notin \mathbb{N} \\ \\ \ds
\delta(\tau)+ \sum_{m=1}^{\infty} 
\frac{(-\nu)_m}{m!}
\frac{(-\xi_0)^m\tau^{\alpha m-1}}{\Gamma(\alpha m)}, & \alpha\nu \in \mathbb{N}
\end{cases} \qquad \tau \in \mathbb{R}^{+} \notag \\
&=\begin{cases} \ds  
\tau^{\ceil{\alpha\nu} -\alpha \nu-1} E_{\alpha,\ceil{\alpha\nu} -\alpha \nu}^{-\nu}(-\xi_0\tau^{\alpha}) , & \alpha\nu \notin \mathbb{N} \\ \\ \ds
  D_{\tau}\left(\Theta(\tau)E_{\alpha,1}^{-\nu}(-\xi_0\tau^{\alpha})\right)   , & \alpha\nu \in \mathbb{N}
\end{cases} \qquad \tau \in \mathbb{R}^{+}
\end{align}
containing the Prabhakar function $E_{a,b}^c(z)$ (\ref{genmittag-Leff}).
Then it follows that
\begin{align}
\label{Doperatorappliedonfunction}
\ds 
{\hat {\cal D}}_{\alpha}^{\nu}(h)\Phi(t) 
&= \ds   \xi_0^{-\nu} \left(h^{-1}(1-{\hat T}_{-h})\right)^{\ceil{\alpha\nu}}
\sum_{m=0}^{\infty} \frac{(-\nu)_m}{m!} (-\xi_0)^m \notag \\
 &\times \sum_{k=0}^{\infty}\frac{(\alpha m + \ceil{\alpha\nu} -\alpha \nu)_k}{k!}
 h^{(\alpha m + \ceil{\alpha\nu} -\alpha \nu)}\Phi(t-kh)  , \hspace{1cm} t \in h\mathbb{Z}
\end{align}
is well-scaled in the sense of the limiting equation (\ref{contilimim-cmul}).
Thus with (\ref{theinter}) we first get for the continuous-time limit 
\begin{equation}
\label{thelimitfunction}
\begin{array}{l} \ds
 {\cal K}^{(0)}_{\alpha,\nu}(t)_{ct}= \lim_{h\rightarrow 0} {\hat K}_{\alpha,\nu}^{(0)}(h)\delta_h(t) \\ \\ \ds \hspace{1.6cm}  =  \xi_0^{-\nu}
D_t^{\ceil{\alpha\nu}-1} \int_0^t{\cal B}_{\alpha,\nu}(\tau)_{ct}\delta(t-\tau){\rm d}\tau
= \xi_0^{-\nu} D_t^{\ceil{\alpha\nu}-1}{\cal B}_{\alpha,\nu}(t)_{ct} 
\end{array} \hspace{0.5cm} t \in \mathbb{R}^{+} .
\end{equation}
Note that, in this limiting procedure derivatives of integer orders come into play with $\ceil{\alpha\nu}-1 =0$ for $0<\alpha\nu \leq 1$ and $\ceil{\alpha\nu}-1 \in \mathbb{N}$ 
for $\alpha\nu >1$. 
We notice that kernel (\ref{kernelwehave})
has Laplace transform ${\cal L}({\cal B}_{\alpha,\nu}(t)_{ct})(s) = s^{-\ceil{\alpha\nu}}\left(s^{\alpha}+\xi_0\right)^{\nu}$. The continuous-time limit kernel (\ref{kernelwehave}) is for $\alpha\nu \notin \mathbb{N}$ by the term $m=0$ weakly singular. 
We then get for the continuous-time limit of (\ref{Doperatorappliedonfunction})
\begin{equation}
\label{ct}
\lim_{h\rightarrow 0} {\hat {\cal D}}_{\alpha}^{\nu}(h) \cdot \Phi(t) 
= \xi_0^{-\nu} D_t^{\ceil{\alpha\nu}} 
\int_0^t {\cal B}_{\alpha,\nu}(t-\tau)_{ct}\Phi(\tau)_{ct}{\rm d}\tau .
\end{equation}
The results (\ref{thelimitfunction}), (\ref{ct}) are in accordance with the continuous-time 
expressions derived for the GFPP (Eqs. (66)-(68) in \cite{TMM-APR-PhysicaA2020}), where the definition of
${\cal K}^{(0)}_{\alpha,\nu}(t)_{ct}$ there differs by a multiplier $\xi_0^{\nu}$.
\\[1ex]
Having derived these continuous-time limits, our goal is now to deduce the recursive generalized fractional
difference equations governing the state probabilities $\Phi^{(n)}_{\alpha,\nu}(t)$ of the PDTP.
To this end let us rewrite Eq. (\ref{recursive-diff})
more explicitly.
Taking into account the evaluations (\ref{discrete-timememoryf}), (\ref{exapandas}) and (\ref{Doperatorappliedonfunction}) this equation takes the form (with $h=1$)
\begin{equation}
\label{explicitbackward}
\begin{array}{l} \ds 
 \xi^{-\nu}\left(1-{\hat T}_{-1}\right)^{\ceil{\alpha\nu}}
\sum_{m=0}^{\infty} \frac{(-\nu)_m}{m!} (-\xi)^m 
 \sum_{k=0}^{\infty}\frac{(\alpha m + \ceil{\alpha\nu} -\alpha \nu)_k}{k!}\Phi_{\alpha,\nu}^{(n)}(t-k) 
 \\ \\ \ds
 = \Phi_{\alpha,\nu}^{(n-1)}(t-1)
 +\delta_{n0} \left\{\xi^{-\nu} \left(1-{\hat T}_{-1}\right)^{\ceil{\alpha\nu}-1}
\sum_{m=0}^{\infty} \frac{(-\nu)_m}{m!} (-\xi)^m 
\frac{(\alpha m + \ceil{\alpha\nu} -\alpha \nu)_t}{t!} - \Theta(t-1)\right\} 
\end{array}
 \end{equation}
where $n, t \in \mathbb{N}_0$ and we wrote this equation for $0<\xi<1$.
For later use and further analysis of the structure it appears instructive to write this equation by accounting 
for 
(\ref{discrete-timememoryf}) as
\begin{equation}
\label{alsointheform}
 (1-{\hat T}_{-1}) \sum_{k=0}^{\infty} {\cal K}^{(0)}_{\alpha,\nu}(k)\cdot \Phi^{(n)}_{\alpha,\nu}(t-k) = 
\Phi^{(n-1)}_{\alpha,\nu}(t-1) +\delta_{n0}\left({\cal K}^{(0)}_{\alpha,\nu}(t)-\Theta(t-1)\right) ,\hspace{0.2cm} n,t \in \mathbb{N}_0 .
\end{equation}
We observe that the $k$-series (due to causality and property $\Phi_{\alpha,\mu}^{(n)}(t-k)=0$ for $t-k<n$) has upper limit $k=t-n$.
The Equations (\ref{recursive-diff}), (\ref{explicitbackward}) and (\ref{alsointheform}) indeed are equivalent discrete-time versions
of the {\it generalized fractional Kolmogorov-Feller equations} governing the 
state-probabilities in a PDTP.
These equations are non-local discrete-time convolutions with non-Markovian long-time memory effects where the whole history $\{\Phi_{\alpha,\nu}^{(n)}(t-k)\}$ comes into play.
We will come back to this issue again later on.
With (\ref{thelimitfunction}) and (\ref{ct}) we can directly write the well-scaled continuous-time limit of 
(\ref{alsointheform}) in the form
\begin{equation}
\label{continuous-time-rep}
\xi_0^{-\nu} D_t^{\ceil{\alpha\nu}} \int_0^t {\cal B}_{\alpha,\nu}(t-\tau)_{ct}\Phi^{(n)}_{\alpha,\nu}(\tau)_{ct}{\rm d}\tau
= \Phi^{(n-1)}_{\alpha,\nu}(t)_{ct} +\delta_{n0}{\cal M}_{\alpha,\nu}(t)_{ct}  ,\hspace{0.5cm} t\in \mathbb{R}^{+}
\end{equation}
which is solved by the generalized fractional Poisson distribution functions (GFPP state probabilities) $\Phi^{(n)}_{\alpha,\nu}(t)_{ct}$ of (\ref{state-prob-scalintg-limit})
with the continuous-time limit kernel
${\cal B}_{\alpha,\nu}(t)_{ct}$ (\ref{kernelwehave}). This equation contains the continuous-time limit of the memory function
\begin{equation}
\label{momrycontilim}
{\cal M}_{\alpha,\nu}(t)_{ct} = \lim_{h\rightarrow 0} {\hat M}^{\alpha,\nu}_h\delta_h(t)=  \xi_0^{-\nu} D_t^{\ceil{\alpha\nu}-1}{\cal B}_{\alpha,\nu}(t)_{ct} -\Theta(t)
\end{equation}
where ${\hat M}^{\alpha,\nu}_h$ is the well-scaled memory operator of (\ref{memoryoperator}).
\\[1ex]
Continuous-time limit memory function (\ref{momrycontilim}) and `generalized fractional Kolmogorov-Feller equations' (\ref{continuous-time-rep}) indeed are in accordance with the results for the GFPP \cite{TMM-APR-PhysicaA2020,MichelitschRiascosGFPP2019}.
\\[2ex]
Let us analyze now more closely the important regime $\nu=1$ with $\alpha \in (0,1]$. 
\\[1ex]
\noindent {\it (i) \underline{Fractional Bernoulli counting process} (of type B), $\nu=1$ and $\alpha \in (0,1)$:}
\newline\newline
We obtain then for the above introduced quantities
\begin{equation}
\label{Kkernel}
\begin{array}{l} \ds 
{\cal K}_{\alpha,1}^{(0)}(t)= \left(\frac{1}{\xi}(1-{\hat T}_{-1})^{\alpha-1}+ (1-{\hat T}_{-1})^{-1}\right)  \delta_{1}(t) \\ \\
\ds \hspace{1.3cm} = \frac{1}{\xi}(1-{\hat T}_{-1})^{\alpha-1}\delta_{1}(t) +\Theta(t)
\end{array} t \in \mathbb{Z}
\end{equation} 
thus 
\begin{equation}
\label{T1D}
 (1-{\hat T}_{-1}) {\cal K}_{\alpha,1}^{(0)}(t) = D_{\alpha}^1\delta_{1}(t) = \left( \frac{1}{\xi} (1-{\hat T}_{-1})^{\alpha} +1\right)\delta_{t0}.
\end{equation}
The memory function (\ref{memoryop}) then yields 
\begin{equation}
\label{memoryfu}
\begin{array}{l} \ds
{\cal M}_{\alpha,1}(t) = {\cal K}_{\alpha,1}^{(0)}(t)-\Theta(t-1) = \left(\frac{1}{\xi}(1-{\hat T}_{-1})^{\alpha-1}+ 1\right) \delta_{1}(t) ,\hspace{0.5cm} 0<\alpha \leq 1 \\ \\ \ds {\cal M}_{\alpha,1}(t) = 
\frac{1}{\xi} \frac{1}{t!} (1-\alpha)_{t} + \delta_1(t)  = \frac{(-1)^{t}}{\xi} 
\left(\begin{array}{l} \alpha-1 \\ \hspace{0.5cm} t \end{array}\right) +\delta_{t0} ,\hspace{0.5cm} t\in \mathbb{N}_0  \\ \\ \ds
{\cal M}_{1,1}(t) = \frac{\xi+1}{\xi} \delta_{t0} ,\hspace{0.5cm} \alpha=1 ,\hspace{0.5cm} t\in \mathbb{N}_0 
\end{array}
\end{equation}
and the well-scaled form
\begin{equation}
\label{well-scaled-memo-alpha}
{\cal M}_{\alpha,1}(t)_h =
 \frac{h}{1-{\hat T}_{-h}}({\cal D}_{\alpha}^1(h)-{\hat T}_{-h})\delta_h(t) = \left(1+\frac{1}{\xi}(1-{\hat T}_{-h})^{\alpha-1} \right) h\delta_h(t), \hspace{0.5cm} t \in h\mathbb{Z}.
\end{equation}
For later use we notice the continuous-time limit
${\cal M}_{\alpha,1}(t)_{ct}= \lim_{h\rightarrow 0}{\cal M}_{\alpha,1}(\frac{t}{h},\xi_0h^{\alpha})= \frac{1}{\xi_0}
\frac{t^{-\alpha}}{\Gamma(1-\alpha)}$ reflecting long-time power-law memory of fractional Bernoulli process.
\\
We used here
 $(1-{\hat T}_{-h})^{\alpha-1}\delta_{h}(t) = h^{-1}\frac{(1-\alpha)_k}{k!}\big|_{k=\frac{t}{h}}$ (See Appendix \ref{crucialdefinitions}).
 \\[1ex]
The discrete-time fractional Kolmogorov-Feller equation then reads ($\delta_1(t)=\delta_{t0}$)
\begin{align}
\label{fracKolmFeller}
&(1-{\hat T}_{-1})^{\alpha} \Phi_{\alpha,1}^{(n)}(t) \notag \\
&=  \xi \Phi_{\alpha,1}^{(n-1)}(t-1) - \xi \Phi_{\alpha,1}^{(n)}(t)
+\delta_{n0} \bigg\{ (-1)^{t}
\left(\begin{array}{l} \alpha-1 \\ \hspace{0.5cm} t \end{array}\right) +\xi \delta_{t0}\bigg\}  ,\hspace{0.15cm} \alpha \in (0,1) ,\hspace{0.15cm} t \in \mathbb{N}_0 .
\end{align}
This equation coincides with the fractional difference equation for the fractional Bernoulli counting process (type B) given in \cite{PachonPolitoRicciuti2018} (See {\it Proposition 7}, with Eqs. (81)-(82) therein).
Eq. (\ref{fracKolmFeller}) is solved by the state probabilities (\ref{state-prob-disceet-Laskin}).
Indeed for $\nu=1$ generating function of the state probabilities (\ref{genfustate}) with ${\bar \varphi}_{\alpha}^{(1)}(u) = \frac{\xi}{\xi+(1-u)^{\alpha}}$ yields
\begin{equation}
\label{indeedyieldsste}
{\bar \Phi}_{\alpha,1}^{(n)}(u)= \frac{(u\xi)^n[\xi+(1-u)^{\alpha-1}]}{[\xi+(1-u)^{\alpha}]^{n+1}} ,\hspace{1cm} n \in \mathbb{N}_0
\end{equation}
which was also given in \cite{PachonPolitoRicciuti2018}. 
On the other hand (\ref{indeedyieldsste}) recovers
for $\alpha=1$ the generating function (\ref{B1}) of standard Binomial distribution.
We can also recover (\ref{fracKolmFeller}) from Eq. (\ref{explicitbackward}) for $\nu=1$ by accounting for $(-1)_m$ is nonzero only for $m=0$ and $m=1$
and where $\ceil{\alpha}=1$. The Eq. (\ref{fracKolmFeller}) is the discrete-time fractional Kolmogorov-Feller equation of this process where the memory function for $t$ large behaves as ${\cal M}_{\alpha,1}^{(0)}(t) \sim \frac{1}{\xi}\frac{t^{-\alpha}}{\Gamma(1-\alpha)}$ thus fractional Bernoulli has a long-time power-law memory 
and non-Markovian features.
\\[1ex]
Now let us consider more closely the continuous-time limit of Eq. (\ref{fracKolmFeller}). First let us write Eq. (\ref{recursive-diff}) by accounting for the memory operator (\ref{memoryoperator}) and with (\ref{shiftD}) in well-scaled form ($\xi=\xi_0h^{\alpha}$) to arrive at
\begin{equation}
\label{we-arrive-at-frac}
\begin{array}{l} \ds 
(1-{\hat T}_{-h})^{\alpha}\Phi_{\alpha,1}^{(n)}(t)_h \\ \\\hspace{0.5cm}  = \xi_0h^{\alpha} 
\Phi_{\alpha,1}^{(n-1)}(t-h)_h - \xi_0h^{\alpha} \Phi_{\alpha,1}^{(n)}(t)_h
+\delta_{n0}\left(\xi_0h^{\alpha} +(1-T_{-h})^{\alpha-1}\right) h\delta_h(t) 
\end{array} t \in h\mathbb{Z}.
\end{equation}
The solution of this equation can be written in well-scaled operator representation (See Eq. (\ref{shift-rep-state}) for $\nu=1$)
\begin{equation}
\label{oprepML}
\Phi_{\alpha,1}^{(n)}(t)_h =
\frac{(\xi_0h^{\alpha} +(1-{\hat T}_{-h})^{\alpha -1})}
{(\xi_0h^{\alpha}+(1-{\hat T}_{-h})^{\alpha})^{n+1}}(\xi_0h^{\alpha})^n\, {\hat T}_{-hn} \, h \delta_h(t)  ,\hspace{1cm} t \in h\mathbb{Z}.
\end{equation}
The limit $h\rightarrow 0$ of (\ref{oprepML}) yields for continuous-time limit
$\Phi_{\alpha,1}^{(n)}(t)_{ct}= \frac{\xi_0^nD_t^{\alpha-1}}{(\xi_0+D_t^{\alpha})^{n+1}}\delta(t)$ 
(with 
Laplace transform ${\cal L}\{\Phi_{\alpha,1}^{(n)}(t)_{ct}\}(s) = \frac{\xi_0^n s^{\alpha-1}}{(\xi_0+s^{\alpha})^{n+1}}$) which can be identified with Laskin's fractional Poisson distribution (\ref{fractionalPoissondistribution})
where $n=0$ recovers the standard Mittag-Leffler survival probability  \cite{Laskin2003,TMM-APR-PhysicaA2020} (among many others).
The continuous-time limit of Eq. (\ref{we-arrive-at-frac}) gives
\begin{equation}
\label{fracKolmFellerwithh}
\lim_{h\rightarrow 0} 
h^{-\alpha}(1-{\hat T}_{-h})^{\alpha} \Phi_{\alpha,1}^{(n)}(t)_h 
= \lim_{h\rightarrow 0} \xi_0\left(\Phi_{\alpha,1}^{(n-1)}(t-h)_h -  \Phi_{\alpha,1}^{(n)}(t)_h\right)
+\delta_{n0} \left(\frac{t^{-\alpha}}{\Gamma(1-\alpha)} +\xi_0 h \, \delta_h(t) \right) .
\end{equation}
By using $\lim_{h\rightarrow 0} h \delta_h(t)  = 0$, however $\delta_h(t) \rightarrow \delta(t)$ together with
$$ 
h^{-\alpha} (1-T_{-h})^{\alpha-1} (h\delta_h(t)) = 
h^{-(\alpha-1)} (1-T_{-h})^{\alpha-1} \delta_h(t) \rightarrow D_t^{\alpha-1}\delta(t) =
\frac{t^{-\alpha}}{\Gamma(1-\alpha)} ,
$$
we obtain as continuous-time limit of (\ref{we-arrive-at-frac}) the fractional differential equation
\begin{equation}
\label{frac-KolFeller}
D_{t}^{\alpha} \cdot \Phi_{\alpha,1}^{(n)}(t)_{ct} = 
+ \xi_0 \Phi_{\alpha,1}^{(n-1)}(t)_{ct} -\xi_0 \Phi_{\alpha,1}^{(n)}(t)_{ct} 
+\delta_{n0}\frac{t^{-\alpha}}{\Gamma(1-\alpha)}  , \hspace{1cm} t \in \mathbb{R}^{+} ,\hspace{1cm} n \in \mathbb{N}_0 .
\end{equation}
In this equation occurs the Riemann-Liouville fractional derivative $D_{t}^{\alpha}$ 
of order $\alpha$ ($0<\alpha <1$)
(See Appendix \ref{appendix1} for details). Eq. (\ref{frac-KolFeller}) coincides with the {\it fractional Kolmogorov-Feller equation}
given by Laskin \cite{Laskin2003} and is solved by the continuous-time state probabilities of Laskin's
fractional Poisson distribution (\ref{fractionalPoissondistribution}). 
\\ \\
{\it (ii) \underline{Bernoulli counting process}, $\nu=1$, $\alpha=1$:}
\\ \\ 
Eq. (\ref{fracKolmFeller}) reduces then to
\begin{equation}
\label{fracKolmFelleralpha1}
(1-{\hat T}_{-1})\Phi_{1,1}^{(n)}(t) 
= \xi \Phi_{1,1}^{(n-1)}(t-1) - \xi \Phi_{1,1}^{(n)}(t) + (\xi+1)\delta_{n0}\delta_{t0}
,\hspace{0.5cm} \alpha=1 , \hspace{0.5cm} t \in \mathbb{Z} 
\end{equation}
where the memory function ${\cal M}_{1,1}(t)= \frac{(\xi+1)}{\xi}\delta_{t0} =0$ is null for $t > 0$. 
This observation explains our choice of the name `memory function' for 
${\cal M}_{\alpha,\nu}(t)$.
The Bernoulli process is memoryless and Markovian and Eq. (\ref{fracKolmFelleralpha1}) is solved by the Binomial distribution (\ref{B2}). 
Eq. (\ref{fracKolmFelleralpha1}) indeed is the Kolmogorov-Feller difference equation of the 
corresponding Bernoulli process with $\nu=1$, $\alpha=1$. To see this let us consider (\ref{fracKolmFelleralpha1}) for $t=0$ and $n=0$:
\begin{equation}
\label{consntnull}
\Phi_{1,1}^{(0)}(0)-\Phi_{1,1}^{(0)}(-1) =  \xi \Phi_{1,1}^{(-1)}(-1) 
- \xi \Phi_{1,1}^{(0)}(0) + (\xi+1)
\end{equation}
by using $\Phi_{1,1}^{(-1)}(0) = \Phi_{1,1}^{(-1)}(-1) =0$ we recover in this equation 
the initial condition $\Phi_{1,1}^{(0)}(0)=1$. 
For $t\geq 1$ from Eq. (\ref{fracKolmFelleralpha1}) we have then for $n=0$
\begin{equation}
\label{Nzeroeq}
\Phi_{1,1}^{(0)}(t)
= \frac{1}{\xi+1}\Phi_{1,1}^{(0)}(t-1) ,\hspace{1cm} t \geq 1.
\end{equation}
Now it follows from this recursion
and the initial condition $\Phi_{1,1}^{(0)}(0)=1$ that the survival 
probability in the Bernoulli process is
\begin{equation}
\label{iniconsol}
\Phi_{1,1}^{(0)}(t)=\frac{1}{(\xi+1)^t} = q^t  ,\hspace{1cm} t \geq 0
\end{equation}
also in accordance with (\ref{B2}) for $n=0$. 
Then let us finally consider (\ref{fracKolmFelleralpha1}) for $n \geq 1$ and $t \geq n$ which writes
\begin{equation}
\label{writesas}
\Phi_{1,1}^{(n)}(t) = 
\frac{1}{\xi+1}\left(\Phi_{1,1}^{(n)}(t-1)+\xi \Phi_{1,1}^{(n-1)}(t-1)\right) ,\hspace{0.5cm} t\geq n\geq 1
\end{equation} 
which is indeed solved by (\ref{B2}), namely
\begin{equation}
\label{issolved}
\begin{array}{l} \ds 
\Phi_{1,1}^{(n)}(t) = \frac{1}{\xi+1}
\left[\left(\begin{array}{l} t-1 \\ n \end{array}\right) \frac{\xi^n}{(\xi+1)^{t-1}}
+\xi \left(\begin{array}{l} t-1 \\ n-1 \end{array}\right)\frac{\xi^{n-1}}{(\xi+1)^{t-1}}\right] \\ \\ \ds
\hspace{0.5cm} = \frac{\xi^n}{(\xi+1)^t}\left[\left(\begin{array}{l} t-1 \\ n \end{array}\right)+\left(\begin{array}{l} t-1 \\ n-1 \end{array}\right)  \right] \\ \\
\ds \hspace{0.5cm} = \left(\begin{array}{l} t \\ n \end{array}\right) \frac{\xi^n}{(\xi+1)^t}
= \left(\begin{array}{l} t \\ n \end{array}\right) p^nq^{t-n} ,\hspace{1cm} t\geq n \\ \\
\ds \hspace{0.5cm}  = \left(\begin{array}{l} t \\ n \end{array}\right) p^nq^{t-n} \Theta(t-n)
\end{array}
\end{equation}
and $\Phi_{1,1}^{(n)}(t)= 0$ for $t<n$.
Hence (\ref{issolved}) together with (\ref{iniconsol}) shows that (\ref{fracKolmFelleralpha1})
indeed is solved by the
state probabilities of the (standard) Bernoulli process (\ref{B2}).
We also verify easily the normalization $$\sum_{n=0}^{\infty}\Phi_{1,1}^{(n)}(t) =
\sum_{n=0}^{t}\left(\begin{array}{l} t \\ n \end{array}\right) p^nq^{t-n}=(p+q)^t=1.$$
\\
The continuous-time limit of 
(\ref{fracKolmFelleralpha1}) then yields (See also Eq. (\ref{we-arrive-at-frac}) for $\alpha\rightarrow 1-0$)
\begin{align}
\label{fracKolmFelleralpha1-continuum-lim}
&\lim_{h\rightarrow 0} (1-{\hat T}_{-h})\Phi_{1,1}^{(n)}(t)_h \notag \\
&= \xi_0 h \,\Phi_{1,1}^{(n-1)}(t-h)_h - \xi_0 h \,\Phi_{1,1}^{(n)}(t)_h + 
\delta_{n0}(\xi_0h +1) h \,\delta_h(t)
,\hspace{0.5cm} \alpha=1 , \hspace{0.5cm} t \in h\mathbb{N}_0 .
\end{align}
Thus we get with 
$(1-{\hat T}_{-h}) \rightarrow hD_t$ and $\delta_h(t) \rightarrow \delta(t)$ the equation
\begin{equation}
\label{standardPoissonKF}
D_{\tau}\Phi_{1,1}^{(n)}(t)_{ct} = \xi_0 \Phi_{\alpha,1}^{(n-1)}(\tau)_{ct}
-\xi_0 \Phi_{\alpha,1}^{(n)}(t)_{ct} +\delta_{n0}\delta(t) 
,\hspace{1cm} t \in \mathbb{R}
\end{equation}
which is indeed solved by the Poisson distribution $\Phi_{1,1}^{(n)}(t)_{ct}=\Theta(t)\frac{(\xi_0t)^n}{n!}e^{-\xi_0t}$ (See also limiting equation (\ref{contlim-bern}))
and with $\Phi_{1,1}^{(0)}(t)_{ct}=e^{-\xi_0t} \Theta(t)$ thus $D_t \Phi_{1,1}^{(0)}(t)_{ct} = 
-\xi_0 \Phi_{1,1}^{(0)}(t)_{ct}+\delta(t)$ occurs in (\ref{standardPoissonKF}) for $n=0$.
The limiting equation (\ref{standardPoissonKF}) reflects the fact that the Bernoulli process is a discrete-time version of the Poisson process, see also \cite{PachonPolitoRicciuti2018} among many others and consult also 
\cite{TMM-APR-PhysicaA2020} (Eqs. (79), (80) in that paper.).
This limiting equation can be recovered from (\ref{frac-KolFeller}) in the limit $\alpha\rightarrow 1-0$ when we account for that $\lim_{\alpha\rightarrow 1-0} D_{\tau}^{\alpha} = D_{\tau}$ reproduces the standard first-order 
derivative (See also Appendix \ref{appendix1}) and memory function $\lim_{\alpha\rightarrow 1-0}{\cal M}_{\alpha,1}(t)_{ct} = \lim_{\alpha\rightarrow 1-0}\frac{1}{\xi_0}\frac{t^{-\alpha}}{\Gamma(1-\alpha)} =\frac{1}{\xi_0}\delta(t)$ becomes a Dirac distribution (See also Eq. (\ref{state-asymp})) reflecting the memoryless Markovian feature of standard Poisson.
\subsection{\small EXPECTED NUMBER OF ARRIVALS IN A PDTP}
\label{avaerga-number}
For many applications especially in problems of diffusive motion, the expected number
of arrivals within $[0,t]$ is of utmost importance.
For the PDTP this quantity is defined as
\begin{equation}
\label{exparrive}
\langle n \rangle_{\alpha,\nu}(t)  = \sum_{n=0}^{\infty} n \Phi_{\alpha,\nu}^{(n)}(t) ,\hspace{1cm} t \in \mathbb{N}_0
\end{equation}
with the state probabilities $\Phi_{\alpha,\nu}^{(n)}(t)$ of Eq. (\ref{state-prob}) where
the upper limit in this summation is $n=t$.
We introduce the important generating function 
\begin{equation}
\label{intro-new-gen-fu}
\begin{array}{l} \ds 
{\bar {\cal G}}_{\alpha,\nu}(u,v) = {\bar {\cal G}}_{\alpha,\nu}(u,\xi,v) =
{\bar \Phi}_{\alpha,\nu}^{(0)}(u) \sum_{n=0}^{\infty} (vu{\bar \varphi}_{\alpha}^{\nu}(u))^n 
= \frac{1-u{\bar \varphi}_{\alpha}^{\nu}}{1-u}\frac{1}{1-uv{\bar \varphi}_{\alpha}^{\nu}(u)}  \\ \\
\ds \hspace{1.7cm} = \frac{1}{(1-u)}
\frac{([1+\frac{1}{\xi}(1-u)^{\alpha}]^{\nu} -u)}{([1+\frac{1}{\xi}(1-u)^{\alpha}]^{\nu} -uv)} = \frac{M^{\alpha,\nu}(u)}{([1+\frac{1}{\xi}(1-u)^{\alpha}]^{\nu} -uv)} ,
\end{array} |v|<1,\hspace{0.5cm} |u| \leq 1 .
\end{equation}
We will see subsequently (Section \ref{CTRM-Prabhakar}) that
this generating function is a key-quantity for stochastic motions with PDTP waiting times.
We refer this generating function to as `{\it PDTP generating function}' as it contains the complete stochastic information such as state probabilities and the expected number of arrivals and also memory generating function (\ref{casezeo}).
\\[1ex]
It is now only a small step to directly extract from (\ref{intro-new-gen-fu}) the asymptotic behavior in the time domain for large $t$. Consider (\ref{intro-new-gen-fu}) for $u\rightarrow 1-0$ which yields for $\alpha\in (0,1)$ the same weakly singular behavior as in the PDTP state probability generating functions and memory generating function (See Eqs. (\ref{the-leading-power-state-prob}), (\ref{state-asymp})), namely

\begin{equation}
\label{asympgenerfunctionPDTP}
{\bar {\cal G}}_{\alpha,\nu}(u,v) \sim \frac{{\bar M}^{\alpha,\nu}(u)}{1-v} 
\sim   \frac{\nu}{\xi(1-v)}(1-u)^{\alpha-1} ,\hspace{0.5cm} u\rightarrow 1-0 ,\hspace{0.5cm} |v|<1 , 
\hspace{0.5cm} \alpha \in (0,1) .
\end{equation}
In the regime $\alpha=1$ we have the finite value ${\bar {\cal G}}_{1,\nu}(u,v)|_{u=1}=\frac{1}{1-v}(1+\frac{\nu}{\xi})$.
Accounting for (\ref{state-asymp}) we then get for the PDTP generating function 
the same type of asymptotic power-law behavior for $t$ large, namely
\begin{equation}
\label{asympt}
{\cal G}_{\alpha,\nu}(t,\xi,v) \sim \frac{\nu}{\xi(1-v)}\frac{t^{-\alpha}}{\Gamma(1-\alpha)} 
\hspace{0.5cm} \alpha \in (0,1) .
\end{equation}
For subsequent use let us also write the scaled time domain representation
\begin{equation}
\label{operatortimescalegenPDTP}
{\cal G}_{\alpha,\nu}(t,\xi,v)_h =
{\cal G}_{\alpha,\nu}\left(\frac{t}{h},\xi,v\right) = \frac{h}{(1-{\hat T}_{-h})}
\frac{([1+\frac{1}{\xi}(1-{\hat T}_{-h})^{\alpha}]^{\nu} -{\hat T}_{-h})}{([1+\frac{1}{\xi}(1-{\hat T}_{-h})^{\alpha}]^{\nu} -v{\hat T}_{-h})} \delta_h(t) ,\hspace{1cm} t \in h \mathbb{Z} 
\end{equation}
which retains the initial condition of the process with $h=1$, namely
${\cal G}_{\alpha,\nu}(t,\xi,v)_h\big|_{t=0}=
{\bar {\cal G}}_{\alpha,\nu}(u,v)\big|_{u=0}=1$.
With the scaling $\xi=\xi_0h^{\alpha}$ (\ref{operatortimescalegenPDTP}) and ${\hat T}_{-h}=e^{-hD_t}$ the well-scaled continuous-time limit writes
\begin{equation}
\label{thewell-scaled}
{\cal G}(t,\xi_0,v)_{ct} = \lim_{h\rightarrow 0}
{\cal G}_{\alpha,\nu}\left(\frac{t}{h},\xi_0h^{\alpha},v\right) = D_t^{-1} \frac{[1+\frac{1}{\xi_0}D_t^{\alpha}]^{\nu} -1}{[1+\frac{1}{\xi_0}D_t^{\alpha}]^{\nu} -v} \, \delta(t) ,\hspace{1cm} t \in \mathbb{R}
\end{equation}
with Laplace transform 
\begin{equation}
\label{Laplace-thewell-scaled}
{\cal L}\{{\cal G}(t,\xi_0,v)_{ct}\}(s) = 
\frac{1}{s} \frac{[1+\frac{1}{\xi_0}s^{\alpha}]^{\nu} -1}{[1+\frac{1}{\xi_0}s^{\alpha}]^{\nu} -v} .
\end{equation}
The expected number of arrivals (\ref{exparrive}) are then given by
\begin{equation}
\label{gen-expnum}
\langle n \rangle_{\alpha,\nu}(t)  = 
\frac{1}{t!}\frac{\partial^{t+1}}{\partial v\partial u^t}{\bar {\cal G}}_{\alpha,\nu}(u,v)\big|_{u=0,v=1} = 
\frac{1}{t!}\frac{\partial^{t}}{\partial u^t}{\cal N}_{\alpha,\nu}(u)
 \big|_{u=0}
\end{equation}
where 
\begin{equation}
\label{Negen}
{\cal N}_{\alpha,\nu}(u) =  \frac{u{\bar \varphi}_{\alpha}^{\nu}(u)}{(1-u)(1- u{\bar \varphi}_{\alpha}^{\nu}(u))}
\end{equation}
is the generating function for the expected number of arrivals in (\ref{gen-expnum}).
Note that all series are convergent for $|v|<1$, $ |u|\leq 1$ (and for $|v|\leq 1$ when $|u|<1$) where
$|{\bar \varphi}_{\alpha}^{\nu}(u)|\leq 1$. For the further evaluation it is useful to expand (\ref{Negen}) as 
\begin{equation}
\label{expandinthatser}
{\cal N}_{\alpha,\nu}(u) = \frac{1}{(1-u)}\sum_{n=1}^{\infty}u^n {\bar \varphi}_{\alpha}^{n \nu}(u).
\end{equation}
From (\ref{gen-expnum}) together with relations (\ref{Pprobab}), (\ref{wegetthen}) follows 
that the terms $n>t$ are zero thus
\begin{equation}
\label{wegetthenuse}
\langle n \rangle_{\alpha,\nu}(t) = 
\sum_{n=1}^t {\cal P}_{\alpha,n\nu}(t-n) =
\sum_{m=0}^{\infty} Q_{\alpha,\nu}^{(m)}(t)  ,\hspace{1cm} t,n \in \mathbb{N}_0
\end{equation}
with
\begin{equation}
\label{Qm}
Q_{\alpha,\nu}^{(m)}(t) =
\sum_{n=1}^t \frac{(-1)^m}{m!}
\xi^{\nu n+m} \frac{(n\nu)_m}{(t-n)!} \left(\alpha(m+ n\nu)+1\right)_{t-n}   ,\hspace{1cm} t,n \in \mathbb{N}_0
\end{equation}
where ${\cal P}_{\alpha,n\nu}(t-n)$ is the probability of (\ref{wegetthen}) for {\it at least} $n$ arrivals within $[0,t]$ and we assumed here that $0<\xi<1$.
\\[1ex]
Now we derive the continuous-time limit by taking into account that
the expected number of arrivals (\ref{exparrive}) is a dimensionless (in the sense of relation (\ref{contilimim-cmul}) a cumulative)
distribution. To this end we account for
\begin{equation}
\lim_{h\rightarrow 0}
 \frac{\xi_0^{\nu n+m} h^{\alpha(\nu n+m)}}{(\frac{t}{h}-n)!} \left(\alpha(m+ n\nu)+1\right)_{\frac{t}{h}-n} =
\frac{(\xi_0t^{\alpha})^{(m+\nu n)}}{\Gamma(\alpha m +\alpha \nu n +1)} ,\hspace{1cm} t \in \mathbb{R} .
\end{equation}
The continuous-time limit of (\ref{Qm}) then is obtained by virtue of relation 
(\ref{contilimim-cmul}) leading to
\begin{equation}
\label{Qmct}
Q_{\alpha,\nu}^{(m)}(t)_{ct} = \lim_{h\rightarrow 0}Q\left(\frac{t}{h},\xi_0h^{\alpha}\right)  = \sum_{n=1}^{\infty} \frac{(-1)^m}{m!}(n\nu)_m
\frac{(\xi_0t^{\alpha})^{(m+\nu n)}}{\Gamma(\alpha m +\alpha \nu n +1)} .
\end{equation}
Hence the continuous-time limit of (\ref{wegetthenuse}) yields
\begin{equation}
\label{wegetthenuseCT}
\begin{array}{l} 
\ds
\langle n \rangle_{\alpha,\nu} (t)_{ct}  = \sum_{n=1}^{\infty} (\xi_0t^{\alpha})^{\nu n}  \sum_{m=0}^{\infty} \frac{(n\nu)_m}{m!}
 \frac{(-\xi_0 t^{\alpha})^{m}}{\Gamma(\alpha m +\alpha \nu n +1)} \\ \\
 \ds = \sum_{n=1}^{\infty} \xi_0^{\nu n} t^{\alpha \nu n} 
 E_{\alpha, \alpha \nu n +1}^{\nu n}(-\xi_0 t^{\alpha}) = 
 \sum_{n=0}^{\infty} \xi_0^{\nu (n+1)} t^{\alpha \nu (n+1)} 
 E_{\alpha, \alpha \nu (n+1) +1}^{\nu (n+1)}(-\xi_0 t^{\alpha}) 
 \end{array} \hspace{0.5cm} t \in \mathbb{R}^{+}
\end{equation}
where $E_{a,b}^c(z)$ denotes the Prabhakar function (\ref{genmittag-Leff}). The result (\ref{wegetthenuseCT}) coincides with the expression 
given in \cite{PolitoCahoy2013} (Eq. (2.36) in that paper).
\\[1ex]
In diffusion problems especially of interest is the asymptotic behavior for large $t$. With the same argument as in Section (\ref{asymfeat}) we infer that this 
asymptotic behavior is contained
in the dominating term for $u\rightarrow 1$, i.e. in the leading power of $(1-u)$. To cover this part
we rewrite generating function (\ref{Negen})
\begin{equation}
\label{genfuexpnumarr}
\begin{array}{l}
 \ds 
{\cal N}_{\alpha,\nu}(u)  = \frac{1}{(1-u)} \frac{u{\bar \varphi}_{\alpha}^{\nu}}{(1-u{\bar \varphi}_{\alpha}^{\nu})} \sim \frac{1}{(1-u)} \frac{{\bar \varphi}_{\alpha}^{\nu}}{(1-{\bar \varphi}_{\alpha}^{\nu})}  \hspace{1cm} (u \to 1-0)\\ \\ \ds
 {\cal N}_{\alpha,\nu}(u) \sim  \frac{\xi}{\nu}(1-u)^{-\alpha-1} \hspace{1cm} (u \to 1-0) 
\end{array} \alpha \in (0,1] .
\end{equation}
where this asymptotic formula also includes $\alpha=1$.
Hence the large-time behavior of (\ref{exparrive}) yields
\begin{equation}
\label{tlargelimitexpnum}
\langle n \rangle_{\alpha,\nu}(t) \sim \frac{\xi}{\nu} \frac{(1+\alpha)_t}{t!} 
\sim \frac{\xi}{\nu} \frac{t^{\alpha}}{\Gamma(1+\alpha)}  ,\hspace{0.5cm} (t\,\, {\rm large}) 
,\hspace{0.5cm} \alpha \in (0,1] .
\end{equation}
This expression contains, by accounting for (\ref{contilimim-cmul}), also the 
continuous-time limit for large observation times $\langle n \rangle_{\alpha,\nu}(t)_{ct}=\lim_{h\rightarrow 0} \langle n \rangle_{\alpha,\nu}(\frac{t}{h},\xi_0h^{\alpha})=\frac{\xi_0}{\nu} \frac{t^{\alpha}}{\Gamma(1+\alpha)}$. We note that the average number of arrivals within [0,t] of (\ref{tlargelimitexpnum}) is a dimensionless function of $t$ ($\xi_0$ has physical dimension $\sec^{-\alpha}$.). Also we mention that $\nu$ for large $t$
enters the power-law (\ref{tlargelimitexpnum}) only as a scaling parameter where the smaller $\nu$ the shorter the waiting times leading to an increase of 
(\ref{tlargelimitexpnum}) for a fixed $t$. 
We also notice that the same $t^{\alpha}$-power law for the mean square displacement is found in ``Prabhakar diffusion'' in the multi-dimensional infinite space. For details we invite the reader to consult the recent articles \cite{TMM-APR-PhysicaA2020,MichelitschRiascosGFPP2019}.
\\[1ex]
It is further noteworthy that in the fractional interval $\alpha \in (0,1)$ the power-law (\ref{tlargelimitexpnum}) reflects asymptotic self-similar scaling. In other words: For long observation times the event stream converges to a stochastic fractal on the time line with scaling dimension $0<\alpha<1$ indicating a disjoint `dust like' fractal. For $\alpha=1$ this behavior turns into linear law indicating compact coverage of the time line by the events. Indeed for $\alpha=1$ relation (\ref{tlargelimitexpnum}) recovers the linear behavior of standard Bernoulli, and of standard Poisson in the continuous-time limit.
The asymptotic power-law (\ref{tlargelimitexpnum}) together with  
(\ref{thewell-scaled}), (\ref{Laplace-thewell-scaled})
are in accordance
with the results obtained for the GFPP \cite{TMM-APR-PhysicaA2020}.
\\[1ex]
Particularly instructive is again the standard Bernoulli case $\alpha=1$ and $\nu=1$. 
(\ref{Negen}) takes then the particularly simple form
\begin{equation}
\label{simpleformbernoulliarrivals}
{\cal N}_{1,1}(u) = \frac{\xi}{(\xi+1)} \frac{u}{(1-u)^2} =\frac{\xi}{\xi+1} \sum_{k=1}^{\infty}k u^k .
\end{equation}
Thus $\left(\frac{p}{q}=\xi\right)$
\begin{equation}
\label{nbernoulli}
\langle n \rangle_{1,1}(t) =\frac{\xi}{\xi+1} t = p t ,\hspace{1cm} t \in \mathbb{N}_0 .
\end{equation}
We can easily reconfirm this well-known result by employing (\ref{exparrive}) accounting for the 
Bernoulli process state probabilities (\ref{B2}). Then we get
$$
\langle n \rangle_{1,1}(t)= \sum_{n=0}^{t} 
\left( \begin{array}{l} t \\ n \end{array}\right) np^nq^{t-n} = p\frac{d}{dp}(p+q)^t = pt .
$$
By virtue of Eq. (\ref{contilimim-cmul}), we can directly derive the continuous-time limit as
\begin{equation}
\label{conlimBernoulli}
\langle n \rangle_{1,1}(t)_{ct} = \lim_{h\rightarrow 0} \frac{\xi_0h }{(\xi_0h +1)} \frac{t}{h} = 
\xi_0 t ,\hspace{1cm} t \in \mathbb{R}^+ 
\end{equation}
in accordance with the above considerations.
This result indeed is the well-known linear law for the expected number of arrivals of the
standard Poisson process (see  e.g.\ \cite{TMM-APR-PhysicaA2020} among many others)
and can be reconfirmed by
$\langle n \rangle_{1,1}(t)_{ct}= e^{-\xi_0 t} \sum_{n=0}^{\infty} n\frac{(\xi_0t)^n}{n!} = \frac{d}{dv}e^{-\xi_0 t} \sum_{n=0}^{\infty} v^n\frac{(\xi_0t)^n}{n!}\big|_{v=1}
=\frac{d}{dv} e^{\xi_0 t(v-1)}\big|_{v=1}
=\xi_0 t$.
\\[1ex]
The results of this subsection are particularly useful in problems of diffusive particle motions with PDTP waiting times between the particle jumps. We devote subsequent section to this subject.

\section{\small PRABHAKAR DISCRETE TIME RANDOM WALK ON UNDIRECTED GRAPHS}
\label{CTRM-Prabhakar} 
The goal of the present section is to analyze the stochastic motion on undirected graphs 
governed by the PDTP.
We consider a walker who performs random steps between connected nodes 
where the IID waiting times between the steps are drawn from a discrete-time renewal process.
If in contrast the sojourn times on the nodes follow a continuous-time renewal process, then the resulting motion is a
Montroll-Weiss continuous-time random walk (CTRW) \cite{MontrollWeiss1965}. 
Instead we consider here walks where the waiting times between the steps follow a PDTP. The so defined walk is a discrete-time generalization 
of the classical Montroll-Weiss CTRW.
We call such a Montroll-Weiss type of walk {\it `Discrete-time random walk' (DTRW)}.
\\[1ex]
Consider now a random walker on an undirected connected graph of $N$ states (nodes).
To characterize the walk on the network we introduce the $N\times N$ 
one-step transition matrix 
(also referred to as stochastic matrix)
${\mathbf H} =(H_{ij})$ where the matrix elements $0 \leq H_{ij} \leq 1$ indicate the 
conditional probabilities
that the walker who is sitting on node $i$ in one step moves to node $j$. The transition matrix is 
normalized as $\sum_{j=1}^N H_{ij}=1$ (i.e. row-stochastic). From row-stochasticity of $\mathbf{H}$ follows row-stochasticity of its matrix powers $\mathbf{H}^n$ ($n\in \mathbb{N}_0$) with the probability $(\mathbf{H}^n)_{ij}$ that the walker in $n$ steps moves from node $i$ to $j$ (For a proof we refer to \cite{TMM-APR-ISTE2019}). The matrix $\mathbf{H}^n$ is the $n$-step transition matrix.
The one-step transition matrix relates the topological information
of the network with the random walk and is defined by \cite{NohRieger2004,RiascosMateos2014,Newman2010,TMM-APR-ISTE2019}
\begin{equation}
\label{One-stp}
H_{ij} = \frac{1}{K_i}A_{ij} = \delta_{ij} - \frac{L_{ij}}{K_i} .
\end{equation}
In this relation the adjacency matrix $A_{ij}=A_{ji}$ (symmetric in an undirected network) is introduced with $A_{ij}=1$ if the 
pair of nodes $ij$ is connected by an edge, and $A_{ij}=0$ otherwise. Further to force the walker 
to change the node in each step we have $A_{jj}=0 ,j=1,\ldots N$, i.e. 
there are no self-connections. The matrix $\mathbf{L} =(L_{ij})$ is referred to as Laplacian-matrix with
$L_{ij}=K_i\delta_{ij}-A_{ij}$. Due to this structure we have $L_{ii}=K_i$ where $K_i=\sum_ {j=1}^N A_{ij}$ 
is the degree of node $i$ indicating the number of neighbor (connected) nodes with node $i$. We see in 
(\ref{One-stp}) that the inverse degrees play the role of row-normalization factors.
In an undirected graph
the adjacency matrix and Laplacian matrices are symmetric whereas the transition-matrix (\ref{One-stp})
in general is not if there is a pair of nodes $i,j$ such that $K_i \neq K_j$.
The spectral structure of Laplacian and
transition matrix is analyzed in details in 
\cite{TMM-APR-ISTE2019} (and see the references therein). We assume here a connected (ergodic) graph. In such a graph the transition matrix has spectral representation 
\begin{equation}
\label{Montroll-W}
{\mathbf H} = |v_1\rangle\langle {\bar v}_1| + \sum_{m=2}^N \lambda_m |v_m\rangle\langle {\bar v}_m|
\end{equation}
with unique eigenvalue $\lambda_1=1$ and $|\lambda_m|<1$ ($m = 2,\ldots,N$)
where $| v_1\rangle\langle {\bar v}_1|$ indicates the stationary (invariant)
distribution with $\lim_{n\rightarrow \infty} {\mathbf H}^n =|v_1 \rangle \langle {\bar v}_1|$. We ignore here for simplicity cases of so called bipartite graphs where an eigenvalue $-1$ occurs. For an outline consult \cite{TMM-APR-ISTE2019} and the references therein.
For our convenience we employ in this section 
Dirac's $\langle bra| ket\rangle$-notation. In this notation ${\mathbf S}=|a\rangle\langle b|$ stands for the $N\times N$ matrix
which has the elements $S_{ij}= \langle i|a\rangle \langle b|j\rangle$
and $\langle i|c\rangle = (\langle c|i\rangle)^{*}$ where 
$(..)^{*}$ indicates complex conjugation.
Further $|v_m\rangle$ denote the right- and $\langle {\bar v}_m|$ 
the left- eigenvectors of the generally non-symmetric transition matrix ${\mathbf H}$. 
We have the properties
$ \langle {\bar v}_n|v_m\rangle =\delta_{mn}$ with the $N\times N$ unity matrix ${\mathbf 1} = \sum_{n=1}^N |v_n\rangle \langle {\bar v}_n|$.
\\[1ex]
We analyze a Montroll-Weiss
DTRW where the IID waiting times between the steps follow a PDTP discrete-time density with generating function (\ref{mostsimplegen}).
Let ${\mathbf P}(t)= (P_{ij}(t))$ be the transition matrix of this walk. 
The element $P_{ij}(t)$ indicates the probability that the walker who is sitting at $t=0$ on node $i$ is present at time $t$ on node $j$.
Since the walker is always somewhere on the network
this transition matrix as well is row-stochastic, 
i.e. fulfills $0\leq P_{ij}(t)\leq 1$ and is normalized as $\sum_{j=1}^N P_{ij}(t)=1$.
We assume the initial condition that at $t=0$ the walker is sitting at node $i$ which is expressed by $P_{ij}(t)\big|_{t=0}= \delta_{ij}$.
Then by simple conditioning arguments we take into account that the walker can move from $i$ to $j$ within the time interval $[0,t]$ in
$n=0,1,,\ldots \in \mathbb{N}_0$ steps where the occurrence of $n$ steps is governed by the PDTP state probabilities
${\mathbb P}(J_n \leq t)=\Phi^{(n)}_{\alpha,\nu}(t)$ (i.e. the probabilities that the walker makes $n$ steps within $[0,t]$).
The transition matrix of a DTRW can then be written as a Cox-series \cite{Cox1967}, namely
\begin{equation}
\label{Coxseries}
{\mathbf P}^{\alpha,\nu}(t) = \sum_{n=0}^{\infty} {\mathbf H}^n \Phi^{(n)}_{\alpha,\nu}(t) ,\hspace{1cm} 
P_{ij}(t)\big|_{t=0}= \delta_{ij} ,\hspace{1cm} t \in \mathbb{N}_0
\end{equation}
containing here the PDTP state probabilities $\Phi^{(n)}_{\alpha,\nu}(t)$ of (\ref{state-prob}). We refer this walk to as {\it Prabhakar DTRW}.
In a general DTRW the state-probabilities of the respective discrete-time counting process are replacing the PDTP state-probabilities.
Transition matrix (\ref{Coxseries}) is a matrix function of ${\mathbf H}$. It follows from (\ref{Coxseries}) that in all equations subsequently derived the matrices have a common base of eigenvectors with ${\mathbf H}$ and hence are commuting among each other and
with ${\mathbf H}$. 
So 
${\mathbf P}^{\alpha,\nu}(t)\cdot {\mathbf H} = {\mathbf H}\cdot {\mathbf P}^{\alpha,\nu}(t)$ is commuting since our
initial condition is ${\mathbf P}(0)={\mathbf 1}$ (commuting with ${\mathbf H}$; however, ${\mathbf P}(t)$ and ${\mathbf H}$ do not 
commute if the initial condition is such that ${\mathbf H} \cdot {\mathbf P}(0) \neq {\mathbf P}(0)\cdot {\mathbf H}$, i.e. if initial condition ${\mathbf P}(0)$ and ${\mathbf H}$ do not have the same base of eigenvectors).
Keep in mind that $ P_{ij}(t)\big|_{t=0} = \delta_{ij}$ indeed requires the state probabilities
to fulfill the initial conditions $\Phi_{\alpha,\nu}^{(n)}(t)\big|_{t=0} = \delta_{n0}$ which per construction is fulfilled (See Eq. (\ref{iniconddesriable})). 
\\[1ex]
For the proofs to follow we make use of the generating function defined by $\ds {\bar {\mathbf P}}^{\alpha,\nu}(u)=\sum_{t=0}^{\infty}{\mathbf P}^{\alpha,\nu}(t) u^t$ of the transition matrix (\ref{Coxseries})
\begin{equation}
\label{genermontrollB}
\begin{array}{l} \ds 
{\bar {\mathbf P}}^{\alpha,\nu}(u)  =
\frac{1- u{\bar \varphi}_{\alpha}^{(\nu)} }{1-u} \sum_{n=0}^{\infty} {\mathbf {H}}^n u^n{\bar \varphi}_{\alpha}^{\nu n}(u)   = \frac{1- u{\bar \varphi}_{\alpha}^{(\nu)} }{1-u}\left({\mathbf 1}- u{\bar \varphi}_{\alpha}^{(\nu)} {\mathbf H}\right)^{-1} \\   \hspace{1.4cm} = {\bar M}^{\alpha,\nu}(u)\left([1+\frac{1}{\xi}(1-u)^{\alpha}]^{\nu}{\mathbf 1} -u{\mathbf H}\right)^{-1}  = {\bar {\cal G}}_{\alpha,\nu}(u,\xi,{\mathbf H})
\end{array}
\end{equation}
which we identify with the PDTP generating function (\ref{intro-new-gen-fu}) with matrix argument $v \rightarrow {\mathbf H}$ containing memory generating function ${\bar M}^{\alpha,\nu}(u)$ of (\ref{casezeo}).
Series (\ref{genermontrollB}) is convergent for $|u|<1$ since $ u{\bar \varphi}_{\alpha}^{(\nu)}  = |\mathbb{E} \, u^t| <  \mathbb{E}\,\, 1=1$ and ${\mathbf H}$ has eigenvalues $|\lambda_j|\leq 1$. We can easily confirm that the transition matrix (\ref{genermontrollB})
is row-stochastic by $\sum_{j=1}^N P^{\alpha,\nu}_{ij}(u) = 
{\bar \Phi}_{\alpha,\nu}^{(n)}(u) = \frac{1}{1-u}$.
Using (\ref{Montroll-W}) with (\ref{genermontrollB}) yields for the transition matrix (\ref{Coxseries}) the canonical representation
\begin{equation}
\label{simple-Eq}
{\mathbf P}^{\alpha,\nu}(t) = {\cal G}_{\alpha,\nu}(t,{\mathbf H}) = |v_1\rangle\langle {\bar v}_1|+\sum_{m=2}^N 
{\cal G}_{\alpha,\nu}(t,\lambda_m) |v_m\rangle\langle {\bar v}_m| 
,\hspace{1cm} t \in \mathbb{N}_0
\end{equation}
with time-domain representation ${\cal G}_{\alpha,\nu}(t,\lambda)={\cal G}(t,\xi,\lambda)$ given in (\ref{operatortimescalegenPDTP}). We considered
that for $\lambda_1=1$ we have $${\cal G}_{\alpha,\nu}(t,\lambda)\big|_{\lambda=1} = 1 ,\hspace{0.5cm}\forall t \in \mathbb{N}_0 .$$
It follows from the asymptotic behavior of the state probabilities, memory function (\ref{state-asymp}) and PDTP generating function (asymptotic relations (\ref{asympgenerfunctionPDTP}) and (\ref{asympt})) that the stationary distribution for large $t$ is approached by the power-law
\begin{equation}
\label{asform}
 {\mathbf P}^{\alpha,\nu}(t) \sim |v_1\rangle\langle {\bar v}_1| 
 + \frac{\nu}{\xi}\frac{t^{-\alpha}}{\Gamma(1-\alpha)}\sum_{m=2}^N \frac{|v_m\rangle\langle {\bar v}_m|}{1-\lambda_m} 
  ,\hspace{0.5cm} \alpha \in (0,1) 
 \end{equation}
and for $\alpha\rightarrow 1-0$ we have in this asymptotic relation $\frac{t^{-\alpha}}{\Gamma(1-\alpha)}\rightarrow \delta(t)= 0$ ($t$ large).
(\ref{asform}) shows the non-Markovian long-time memory feature of the Prabhakar DTRW.
Noteworthy here the already mentioned limit $\alpha\rightarrow 0+$ where extremely long waiting times between the steps occur with $\frac{t^{-\alpha}}{\Gamma(1-\alpha)} \sim \Theta(t) =1 $. The walk then becomes `infinitely slow' thus there is a range for $t$ large where $ {\mathbf P}^{\alpha,\nu}(t) \sim |v_1\rangle\langle {\bar v}_1| 
 + \frac{\nu}{\xi}\sum_{m=2}^N \frac{|v_m\rangle\langle {\bar v}_m|}{1-\lambda_m} $ and thus
 the walk `struggles' to take the stationary distribution $|v_1\rangle\langle {\bar v}_1|$ (which eventually is taken since $\frac{t^{-\alpha}}{\Gamma(1-\alpha)} \rightarrow 0+$ for infinitesimally small positive $\alpha$). 
\\[1ex]
Now our goal is to derive the evolution equation that governs the PDTP transition matrix.
To this end we account for
\begin{equation}
\label{theformofmemoryeq}
\frac{{\bar {\cal G}}_{\alpha,\nu}(u,\lambda)}{{\bar \varphi}_{\alpha}^{(\nu)}(u)} =  u \lambda  
{\bar {\cal G}}_{\alpha,\nu}(u,\lambda) + {\bar M}^{\alpha,\nu}(u) 
\end{equation}
with memory generating function ${\bar M}^{\alpha,\nu}(u)$ (\ref{casezeo}). 
Rewriting (\ref{theformofmemoryeq}) in operator form yields
\begin{equation}
\label{shiftreprentationKoFeEq}
 {\hat {\cal D}}_{\alpha}^{\nu} {\cal G}_{\alpha,\nu}(t,\lambda)  = \lambda  {\cal G}_{\alpha,\nu}(t-1,\lambda) + {\cal M}_{\alpha,\nu}(t) ,\hspace{1cm} t \in \mathbb{N}_0  
\end{equation}
containing the memory function ${\cal M}_{\alpha,\nu}(t)$ defined in (\ref{discrete-timememoryf}). 
We then can write (\ref{shiftreprentationKoFeEq}) in matricial representation as
\begin{equation}
\label{straight}
{\hat {\cal D}}_{\alpha}^{\nu} {\mathbf P}^{\alpha,\nu}(t) - {\cal M}_{\alpha,\nu}(t) 
{\mathbf 1}
= {\mathbf H} \cdot {\mathbf P}^{\alpha,\nu}(t-1)   ,\hspace{1cm} P^{\alpha,\nu}_{ij}(0) =\delta_{ij} .
\end{equation}
This equation is the {\it Kolmogorov-Feller generalized fractional difference equation}
that governs the Prabhakar DTRW on the network. Be reminded that with our initial condition ${\mathbf P}(0)={\mathbf 1}$ the matrices on the right-hand side commute.
Since
${\cal M}_{\alpha,\nu}(t)\big|_{t=0} = {\cal D}_{\alpha}^{\nu}\delta_{1}(t)\big|_{t=0} = \frac{(\xi+1)^{\nu}}{\xi^{\nu}}$ fulfill the same initial condition, Eq. (\ref{straight}) recovers for $t=0$ the initial condition of the transition matrix (where due to causality the right-hand side for $t=0$ is null).
\\[1ex]
Eq.\ (\ref{straight}) refers to the general class 
of equations governing discrete-time semi-Markov chains given in \cite{PachonPolitoRicciuti2018} (See Theorem 3.4.).
We can conceive (\ref{straight}) as the discrete-time Cauchy problem 
which is solved by transition matrix (\ref{simple-Eq}) where the complete history $\{ P^{\alpha,\nu}_{ij}(t-k) \}$ ($0\leq k\leq t$) of the walk comes into play. This shows the following representation of (\ref{straight}), namely
\begin{equation}
\label{equithis}
(1-{\hat T}_{-1})\sum_{k=0}^{\infty}{\cal K}_{\alpha,\nu}^{(0)}(k)
P^{\alpha,\nu}_{ij}(t-k) 
- \delta_{ij}\left[{\cal K}_{\alpha,\nu}^{(0)}(t)-\Theta(t-1)\right]  =  \sum_{r=1}^NH_{ir} P^{\alpha,\nu}_{rj}(t-1) ,\hspace{0.5cm} t\in \mathbb{N}_0 
\end{equation}
also reflecting Eq. (\ref{alsointheform}).
This equation contains on the right-hand side the topological information of the graph (See (\ref{One-stp})).
We observe as a consequence of causality ($P^{\alpha,\nu}_{ij}(t-k)= 0$ for $k>t$) that the upper limit of the $k$-summation on the left-hand side is $k=t$.
Eqs. (\ref{straight}), (\ref{equithis}) are equivalent representations of the discrete-time Kolmogorov-Feller equations governing the stochastic motion of a Prabhakar DTRW.
These equations are explicit accounting for relations (\ref{discrete-timememoryf}) and 
(\ref{Doperatorappliedonfunction}).
\\ \\ 
\noindent {\it Continuous-time limit}
\\ \\
Then by the same scaling arguments as in previous sections and outlined in Appendix \ref{crucialdefinitions}, it is not a big deal to establish the continuous time
limit of these equations. In view of (\ref{alsointheform}) having continuous-time limit (\ref{continuous-time-rep})
we arrive at
\begin{equation}
\label{contitimeKFeqs}
\begin{array}{l} \ds 
\xi_0^{-\nu} D_t^{\ceil{\alpha\nu}} \int_0^t{\cal B}_{\alpha,\nu}(\tau)_{ct}
P^{\alpha,\nu}_{ij}(t-\tau)_{ct}{\rm d}\tau -\delta_{ij}\left(\xi_0^{-\nu}D_t^{\ceil{\alpha\nu}-1} 
{\cal B}_{\alpha,\nu}(t)_{ct}-\Theta(t)\right) =  
\sum_{r=1}^NH_{ir} P^{\alpha,\nu}_{rj}(t)_{ct}  ,\\ \\ \ds
\xi_0^{-\nu} D_t^{\ceil{\alpha\nu}-1} \int_0^t{\cal B}_{\alpha,\nu}(\tau)_{ct}\,
D_tP^{\alpha,\nu}_{ij}(t-\tau)_{ct}{\rm d}\tau = \sum_{r=1}^N (H_{rj}-\delta_{ij}) P^{\alpha,\nu}_{ir}(t)_{ct} 
\end{array} \hspace{-1cm} t \in \mathbb{R}^{+}
\end{equation}
where $\delta_{ij} = P^{\alpha,\nu}_{ij}(t)_{ct}|_{t=0}$ is the initial condition and in the formulation of the last line we used on the right-hand side row-normalization of the transition matrix.
The continuous-time limit kernel ${\cal B}_{\alpha,\nu}(t)_{ct}$ was determined in (\ref{kernelwehave})
and $ {\cal K}^{(0)}_{\alpha,\nu}(t)_{ct}= \xi_0^{-\nu}D_t^{\ceil{\alpha\nu}-1} {\cal B}_{\alpha,\nu}(t)_{ct}$
in relation (\ref{thelimitfunction}).
Eq. (\ref{contitimeKFeqs}) is in accordance with the `{\it generalized fractional Kolmogorov-Feller equation}' derived for the Prabhakar CTRW on undirected networks \cite{TMM-APR-PhysicaA2020,MichelitschRiascosGFPP2019}.
\\ \\
{\it Case $\nu=1$ , $0<\alpha \leq 1$:}
\\ \\
Let us discuss the case $\nu=1$ with $0<\alpha < 1$ more closely, i.e.\ the walk subordinated to {\it fractional Bernoulli} (type B).
We refer this walk to as {\it Fractional Bernoulli Walk (FBW)}.
(\ref{equithis})
then takes with Eqs. (\ref{Kkernel})-(\ref{memoryfu}) the form of a fractional difference equation:
\begin{equation}
\label{purefractional-case}
\left(1-{\hat T}_{-1}\right)^{\alpha}P^{\alpha,1}_{ij}(t) - \delta_{ij}\,\left[ (-1)^t\left(\begin{array}{l}
\alpha-1 \\ \hspace{0.5cm} t\end{array}\right)  +\xi \delta_{t0}\right] =  -\xi P^{\alpha,1}_{ij}(t)+ \xi \sum_{r=1}^NH_{ir}P^{\alpha,1}_{rj}(t-1) ,\hspace{0.2cm} t \in \mathbb{N}_0 .
\end{equation}
This equation for $t=0$ yields the initial 
condition $P^{\alpha,1}_{ij}(t)\big|_{t=0}= \delta_{ij}$.
The memory term \\ \\    $$ \hspace{1cm} (-1)^t\left(\begin{array}{l}
\alpha-1 \\ \hspace{0.5cm} t\end{array}\right) \sim \frac{t^{-\alpha}}{\Gamma(1-\alpha)}
,\hspace{1cm} (t \,\, {\rm large}) ,\hspace{0.5cm} \alpha \in (0,1)   $$ 
\\ reflects the long-time memory and non-Markovianity of the process where this memory for $\alpha \rightarrow 0+0 $ becomes extremely long with $\frac{t^{-\alpha}}{\Gamma(1-\alpha)} \rightarrow \Theta(t)= 1 $; whereas $\alpha\rightarrow 1-0 $ represents the memoryless limit with $\frac{t^{-\alpha}}{\Gamma(1-\alpha)} \rightarrow \delta(t) =0$.
In order to analyze the continuous-time limit we account for the generating function (\ref{genermontrollB}) of the transition matrix
\begin{equation}
\label{Coxgen-fract-bernoulli}
{\bar {\mathbf P}}^{\alpha,1}(u) = (\xi+(1-u)^{\alpha-1})\left[\xi({\mathbf 1}-u{\mathbf H})+{\mathbf 1}(1-u)^{\alpha}\right]^{-1} .
\end{equation}
The transition matrix of the FBW has the well-scaled operator representation
\begin{equation}
\label{well-scaled-frac-ber}
{\mathbf P}^{\alpha,1}(t)_h  = (\xi +(1-{\hat T}_{-h})^{\alpha-1})\left[\xi({\mathbf 1}-{\mathbf H}{\hat T}_{-h})+
{\mathbf 1}(1-{\hat T}_{-h})^{\alpha}\right]^{-1} h\delta_h(t) ,\hspace{0.5cm} t \in h\mathbb{Z}.
\end{equation}
We utilize here notation ${\mathbf P}^{\alpha,1}(t)_h$ with subscript $(\ldots)_h$ 
when we refer to the time scaled walk with $t \in h\mathbb{Z}$.
It follows from representation (\ref{well-scaled-frac-ber}) that ${\mathbf P}^{\alpha,1}(t)_h$ solves the Cauchy problem
\begin{equation}
\label{PFBWsolves}
\begin{array}{l} \ds 
(1-{\hat T}_{-h})^{\alpha}{\mathbf P}^{\alpha,1}(t)_h - {\mathbf 1} (\xi +(1-{\hat T}_{-h})^{\alpha-1}) \, h\delta_h(t) =
\xi({\mathbf H}{\hat T}_{-h} - {\mathbf 1}){\mathbf P}^{\alpha,1}(t)_h     \\ \\ \ds
{\mathbf P}^{\alpha,1}(t)_h\big|_{t=0} = {\mathbf 1}
\end{array}
\hspace{0.5cm} t \in h\mathbb{Z}
\end{equation}
where this equation is the scaled version of (\ref{purefractional-case}). This equation is also consistent with
the fractional difference equations (\ref{fracKolmFeller}), (\ref{we-arrive-at-frac}) for the fractional Bernoulli state probabilities.
Consider generating function (\ref{intro-new-gen-fu}) with (\ref{operatortimescalegenPDTP}) for $\nu=1$, namely
\begin{equation}
\label{generalpha1}
{\cal E}_{\alpha}\left(\frac{t}{h},\xi,v\right)=  {\cal G}_{\alpha,1}\left(\frac{t}{h},v\right) =\sum_{n=0}^{\infty} \Phi_{\alpha,1}^{(n)}\left(\frac{t}{h},\xi\right) v^n =  
\frac{\xi +(1-{\hat T}_{-h})^{\alpha-1}}{\xi(1-v{\hat T}_{-h})+ (1-{\hat T}_{-h})^{\alpha}} h\delta_h(t) ,\hspace{0.2cm} t \in h\mathbb{Z}
\end{equation}
which also is obtained by using Eq. (\ref{oprepML}).
The scaled state distribution of the fractional Bernoulli process (\ref{state-prob-disceet-Laskin}) is then obtained from this generating function by
\begin{equation}
\label{generfustaefraBer}
\Phi_{\alpha,1}^{(n)}\left(\frac{t}{h},\xi\right) = \frac{1}{n!} \frac{d^n}{dv^n} {\cal E}_{\alpha}\left(\frac{t}{h},\xi,v\right)\Big|_{v=0} .
\end{equation}
This is the analogue equation as for the fractional Poisson distribution in the continuous-time limit which is shown a little later. We also show that generating function (\ref{generalpha1})
is a discrete-time approximation of the standard Mittag-Leffler function.
The transition matrix (\ref{well-scaled-frac-ber}) can then be written
in the form of the matrix function
\begin{equation}
\label{thesulution}
{\mathbf P}^{\alpha,1}(t)_h = {\cal E}_{\alpha}\left(\frac{t}{h},\xi,{\mathbf H} \right) ,\hspace{1cm} t \in h\mathbb{Z}.
\end{equation}
Let us first consider the continuous-time limit of generating function (\ref{generalpha1})
\begin{equation}
\label{genercontilimvfractber}
\begin{array}{l} \ds
{\cal E}_{\alpha}\left(t,\xi_0,v\right)_{ct} =
\lim_{h\rightarrow 0} {\cal E}_{\alpha}\left(\frac{t}{h},\xi_0h^{\alpha},v\right)  =
\lim_{h\rightarrow 0} \frac{\xi_0h +h^{1-\alpha}(1-{\hat T}_{-h})^{\alpha-1}}{\xi_0(1-v{\hat T}_{-h})+ h^{-\alpha}(1-{\hat T}_{-h})^{\alpha}} \delta_h(t) \\ \\ \ds \hspace{2.1cm}  = \frac{D_t^{\alpha-1}}{\xi_0(1-v)+D_t^{\alpha}} \delta(t)
= E_{\alpha}(-\xi_0(1-v)t^{\alpha}) ,\hspace{1cm} t \in \mathbb{R}^{+}
\end{array}
\end{equation}
retaining the initial condition ${\cal E}_{\alpha}(t,\xi_0,v)_{ct}\big|_{t=0}=1$ of (\ref{generalpha1}).
In this relation the standard Mittag-Leffler function 
\begin{equation}
\label{MLfunction}
E_{\alpha}(z) = \sum_{m=0}^{\infty}\frac{z^m}{\Gamma{(\alpha m +1)}}
\end{equation}
comes into play. This result is easily confirmed in view of 
Laplace transform ${\cal L}\{E_{\alpha}(-at^{\alpha})\}(s) = \frac{s^{\alpha-1}}{a+s^{\alpha}}$ of the Mittag-Leffler function. 
In this way we have shown that generating function (\ref{generalpha1}) also is a discrete-time approximation of the Mittag-Leffler function (\ref{genercontilimvfractber}).
Indeed Laskin's fractional Poisson distribution (\ref{fractionalPoissondistribution}) is obtained from the Mittag-Leffler generating function by \cite{Laskin2003,michel-riascos-springer2020} (and many others)
\begin{equation}
\label{frac-Poisson}
\Phi_{\alpha,1}(t)_{ct}= \frac{1}{n!} \frac{d^n}{dv^n}
E_{\alpha}(-\xi_0(1-v)t^{\alpha})\big|_{v=0} 
\end{equation}
which is also the continuous-time limit of Eq. (\ref{generfustaefraBer}).
The result (\ref{genercontilimvfractber}) allows us to get the continuous-time limit of 
the transition matrix (\ref{thesulution}) in the form of the Mittag-Leffler matrix function
\begin{equation}
\label{contlimfracbernoullitransmat}
\begin{array}{l} 
\ds  {\mathbf P}^{\alpha,1}(t)_{ct} =\lim_{h\rightarrow 0} 
{\cal E}_{\alpha}\left(\frac{t}{h},\xi_0h^{\alpha},{\mathbf H}\right) =
E_{\alpha}\left(-\xi_0({\mathbf 1}-{\mathbf H})t^{\alpha}\right) \\ \\ \ds \hspace{1.6cm} = |v_1\rangle\langle {\bar v}_1| 
+\sum_{m=2}^N |v_m\rangle\langle {\bar v}_m| E_{\alpha}(-\xi_0(1-\lambda_m)t^{\alpha}) 
,\hspace{0.5cm} t \in \mathbb{R}^{+}
\end{array}
\end{equation}
retaining initial condition $\it {\mathbf P}^{\alpha,1}(t)_{ct}\big|{t=0}={\mathbf 1}$. Accounting for the Mittag-Leffler asymptotic relation (holding for the eigenvalues $|\lambda|<1$)
$\it E_{\alpha}(-\xi_0(1-\lambda)t^{\alpha}) \sim \frac{1}{(1-\lambda)\xi_0}\frac{t^{-\alpha}}{\Gamma(1-\alpha)}$ for large $t$, we observe \eqref{contlimfracbernoullitransmat} agrees with (\ref{asform}) for $\nu=1$.
\\[1ex]
The continuous-time limit of Eq. (\ref{purefractional-case}) yields 
\begin{equation}
\label{contilim-frac}
\begin{array}{l} \ds 
D_t^{\alpha} P^{\alpha,1}_{ij}(t)_{ct} - \delta_{ij}\, \frac{t^{-\alpha}}{\Gamma(1-\alpha)} = \xi_0 
\sum_{r=1}^N(H_{ir} - \delta_{ir}) P^{\alpha,1}_{rj}(t)_{ct}  ,\\ \\    \ds    
\frac{1}{\Gamma(1-\alpha)}\int_0^t(t-\tau)^{-\alpha}D_{\tau}P^{\alpha,1}_{ij}(\tau)_{ct} = -\frac{\xi_0}{K_i}\sum_{r=1}^N L_{ir}P^{\alpha,1}_{rj}(t)_{ct}
\end{array} \hspace{1cm} t \in \mathbb{R}^{+} ,\hspace{0.5cm} 0<\alpha<1
\end{equation}
where the first line contains the Riemann-Liouville fractional derivative $D_t^{\alpha}$ of order $\alpha$ whereas
in the second line we utilize the Caputo-fractional derivative of 
order $\alpha$.
The fractional evolution equation (\ref{contilim-frac}) with initial condition $P^{\alpha,1}_{ij}(t)_{ct}\big|_{t=0}=\delta_{ij}$ indeed is solved by the Mittag-Leffler transition matrix (\ref{contlimfracbernoullitransmat}).
\\[1ex]
Eq.\ (\ref{contilim-frac}) is the fractional Kolmogorov-Feller equation governing fractional diffusion in the network,
i.e.\ a random walk in the network subordinated to a fractional Poisson process.
The fractional differential continuous-limit equation (\ref{contilim-frac})
refers to the class of equations governing semi-Markov processes related to $\alpha$-stable subordinators in \cite{PachonPolitoRicciuti2018} (Eqs. (14), (15) in that paper). 
Fractional differential equations of this type (mostly for infinite continuous spaces) with Mittag-Leffler solutions 
occur in a wide range of problems in fractional dynamics and anomalous diffusion
(see e.g.\ \cite{Mainardi-et-al2000,Zaslavsky2002,MainardiGorenfloScalas2004,TMM-APR-PhysicaA2020,MichelitschRiascosGFPP2019,MetzlerKlafter2000,MetzlerKlafter2004}). 
\\[1ex]
It remains us to consider the case $\nu=1$, $\alpha=1$ which is a walk subordinated to a {\it standard Bernoulli process}. We refer this walk to as {\it Bernoulli Walk (BW)}.
We then get for (\ref{equithis}) the difference equation
\begin{equation}
\label{alphanuone}
\left(1-{\hat T}_{-1}\right)P^{1,1}_{ij}(t) - \delta_{ij} \delta_{t0} (\xi + 1) = 
  \xi \sum_{r=1}^N (H_{ir}{\hat T}_{-1}-\delta_{ir}) P^{1,1}_{rj}(t) 
 ,\hspace{0.5cm} t \in \mathbb{N}_0 
\end{equation}
where with ${\hat T}_{-1}P^{1,1}_{rj}(t)=P^{1,1}_{rj}(t-1)$ and the memory term $\sim \delta_{t0}$ is null for $t>0$, i.e. the process is memory-less and Markovian reflecting these properties of standard Bernoulli process.
Eq. (\ref{alphanuone}) is in accordance with the difference equation given in \cite{PachonPolitoRicciuti2018} holding for Markov chains (See Eq. (65) in that paper).
Let us rewrite Eq. (\ref{alphanuone}) also in matricial representation
\begin{equation}
\label{matrical}
(1+\xi){\mathbf P}^{1,1}(t)  -\delta_{t0}{\mathbf 1}(\xi+1) = \left({\mathbf 1}+\xi{\mathbf H}\right){\mathbf P}^{1,1}(t-1) ,\hspace{0.5cm} t \in \mathbb{N}_0 .
\end{equation}
By using causality, i.e. ${\mathbf P}^{1,1}(-1)=0$ this equation recovers for $t=0$ the initial condition ${\mathbf P}^{1,1}(0)= {\mathbf 1}$. For $t\geq 1$ this equation yields the recursion
\begin{equation}
\label{recursionalphanuone}
{\mathbf P}^{1,1}(t) =  \frac{1}{1+\xi}
\left({\mathbf 1}+\xi{\mathbf H}\right){\mathbf P}^{1,1}(t-1) ,\hspace{0.5cm} {\mathbf P}^{1,1}(0)= {\mathbf 1}.
\end{equation}
Thus by iterating this recursion yields for the transition matrix ${\mathbf P}^{1,1}(t)= \frac{1}{(1+\xi)^t}\left({\mathbf 1}+\xi {\mathbf H}\right)^t$.
On the other hand in view of the Binomial distribution (\ref{issolved}) we can obtain this result also by employing the Cox series (\ref{Coxseries}), namely
\begin{equation}
\label{transmat}
{\mathbf P}^{1,1}(t) = {\mathbf P}^{1,1}(t,\xi) = \frac{1}{(1+\xi)^t}\sum_{n=0}^{t} 
\left(\begin{array}{l} t \\ n \end{array}\right)\xi^n{\mathbf H}^n = 
\frac{1}{(1+\xi)^t}\left({\mathbf 1}+\xi {\mathbf H}\right)^t =(p{\mathbf H}+q{\mathbf 1})^t , \hspace{0.5cm} t \in \mathbb{N}_0 
\end{equation}
where $p=\frac{\xi}{1+\xi}$ and $q=\frac{1}{1+\xi}$. In view of (\ref{transmat}) and also with Binomial distribution (\ref{issolved}) it follows that $p$ is the probability that the walker makes a step and $q$ that it does not make a step in a time unit $1$.
For $p=1$, $q=0$ the walker hence makes (almost surely) $t$ steps up to time $t$ and (\ref{transmat}) recovers
${\mathbf P}^{1,1}(t)= {\mathbf H}^t$ the definition of the $t$-step transition matrix. On the other hand
for $p=0$, $q=1$ the walker (almost surely) does not move thus in this case the walker remains on his departure node with ${\mathbf P}^{1,1}(t)= {\mathbf 1}$.
\\[1ex]
Now we can directly derive the continuous-time limit of (\ref{transmat}) by considering the process on $t \in h\mathbb{N}_0$ to arrive at (See Eq. (\ref{contilimim-cmul}))
\begin{equation}
\label{contilim-poisson}
\begin{array}{l} \ds 
{\mathbf P}^{1,1}(t)_{ct} = \lim_{h\rightarrow 0} {\mathbf P}^{1,1}\left(\frac{t}{h},\xi_0h\right)  =  \lim_{h\rightarrow 0} 
\frac{1}{(1+\xi_0h)^{\frac{t}{h}}}
\left({\mathbf 1}+\xi_0h {\mathbf H}\right)^{\frac{t}{h}} 
\\ \\ \ds  \hspace{1.5cm}
= e^{-\xi_0 t} e^{\xi_0 {\mathbf H}t} = e^{-\xi_0 t ({\mathbf 1}-{\mathbf H})},\hspace{0.5cm} t \in \mathbb{R}^{+}.
\end{array}
\end{equation}
This result is also consistent with the Cox-series generated with the Poisson-distribution state probabilities
\begin{align}
\nonumber
{\mathbf P}^{1,1}(t)_{ct} &= \Theta(t) e^{-\xi_0 t} \sum_{n=0}^{\infty} {\mathbf H}^n\frac{(\xi_0 t)^n}{n!} = \Theta(t) e^{-\xi_0 t} e^{\xi_0 {\mathbf H}t} \\ \label{cox-poisson}
&= \Theta(t)  ( |v_1\rangle\langle {\bar v}_1|
+ \sum_{m=2}^N |v_m\rangle\langle {\bar v}_m| e^{\xi_0(\lambda_m-1)t} ),\hspace{1cm} t \in \mathbb{R}
\end{align}
where we added here the Heaviside function to emphasize causality.
By accounting for the spectral structure of $\mathbf{H}$, namely $\lambda_1=1$
and $m-1$ eigenvalues with $\lambda_m-1 <0$ we have $\lim_{t\rightarrow \infty} e^{\xi_0(\lambda_m-1)t} \rightarrow 0$ for $m=2,\ldots N$ thus
the transition matrix ${\mathbf P}^{1,1}(t)_{ct} \rightarrow |v_1\rangle\langle {\bar v}_1|$ approaches asymptotically the stationary distribution. We also notice that the matrix exponential (\ref{contilim-poisson}), (\ref{cox-poisson}) is recovered for $\alpha=1$ by the Mittag-Leffler transition matrix (\ref{contlimfracbernoullitransmat}).
\\[1ex]
Finally the continuous-time limit of Eq. (\ref{alphanuone}) is obtained as 
\begin{equation}
\label{obtainedasstandardpoisson}
\begin{array}{l} \ds 
D_t P^{1,1}_{ij}(t)_{ct} - \delta_{ij}\delta(t) =
 -\xi_0 P^{1,1}_{ij}(t)_{ct} + \xi_0 \sum_{r=1}^NH_{ir}P^{1,1}_{rj}(t)_{ct} = -\frac{\xi_0}{K_i} \sum_{r=1}^N L_{ir} P^{1,1}_{rj}(t)_{ct}\\ \\ \ds
 D_t{\mathbf P}^{1,1}(t)_{ct} -{\mathbf 1}\delta(t) =  \xi_0\left({\mathbf H}-{\mathbf 1}\right){\mathbf P}^{1,1}(t)_{ct} 
 \end{array}
 \hspace{0.1cm} t \in \mathbb{R}
\end{equation}
with initial condition ${\mathbf P}^{1,1}(t)_{ct}|_{t=0} = {\mathbf 1}$. 
The continuous-time limit (\ref{obtainedasstandardpoisson}) is the {\it Kolmogorov-Feller equation} governing the transition probabilities in continuous-time Markov chains, 
i.e. CTRWs on graphs with underlying standard Poisson process. It is straightforward to see that the Cauchy problem  (\ref{obtainedasstandardpoisson}) is solved by the continuous-time limit exponential transition matrix (\ref{contilim-poisson}), (\ref{cox-poisson}). The Eq. (\ref{obtainedasstandardpoisson}) can be also be recovered from the fractional case (\ref{contilim-frac}) 
in the limit $\alpha\rightarrow 1-0$.
{\it Kolmogorov-Feller equations} are widely used as master equations to model
{\it Markovian walks} on graphs (Markov chains) and normal diffusion in multidimensional infinite spaces \cite{PachonPolitoRicciuti2018,TMM-APR-ISTE2019} (and see the references therein).
\\[1ex]
In this section we focused on walks subordinated to PDTPs which we refer to as Prabhakar DTRWs and their special cases
such as fractional Bernoulli and standard Bernoulli as well as their classical continuous-time limits. The advantage of the Prabhakar DTRW with the PDTP counting process of three parameters $\alpha \in (0,1]$, $\nu>0$, $\xi>0$ is the great flexibility to adapt to real-world stochastic processes.
\\[1ex]
Another interesting case (though not being Prabhakar) is the DTRW with Sibuya counting process. A detailed analysis of this walk is beyond the scope of the present paper. The essential aspects are analyzed in \cite{PachonPolitoRicciuti2018}. We confine us here to a brief outline in the spirit of a Montroll Weiss Sibuya DTRW 
model in Appendix \ref{Sibuya-DTRW}.
\section{\small INFLUENCE OF WAITING TIME INITIAL CONDITIONS ON DTRW FEATURES}
\label{influence-ini}
This section is devoted to briefly analyze the effect of the initial condition of the discrete-time waiting-time density
in a DTRW.
In this paper we have constructed discrete-time counting processes with discrete-time waiting time densities $w(t)$ on $t\in \{0,1,2,\ldots \}$ with $w(t)\big|_{t=0}=0$. 
Here our aim is to consider the effect of a discrete-time counting process which has $w(0)=1-\epsilon$ with 
$0<\epsilon \leq 1$.
The waiting time generating function of such a process then is given by
\begin{equation}
\label{genfu}
\mathbb{E} u^Z = {\bar w}(u) =1-\epsilon +u \sum_{k=1}^{\infty}u^{k-1}w(k).
\end{equation}
The state probabilities (probabilities for $n$ arrivals within $[0,t]$) then have the generating function
\begin{equation}
\label{generfustate}
{\bar \Phi}^{(n)}(u)= \frac{1-{\bar w}(u)}{1-u} ({\bar w}(u))^n ,\hspace{0.5cm} n=0,1,2,\ldots 
\end{equation}
Now let us consider the initial condition of the state probabilities
\begin{equation}
\label{inicons}
\Phi^{(n)}(t)\big|_{t=0} = {\bar \Phi}^{(n)}(u)\big|_{u=0} = \epsilon (1-\epsilon)^n ,\hspace{0.5cm} n=0,1,2,\ldots 
\end{equation}
with normalization $\ds \sum_{n=0}^{\infty}\Phi^{(n)}(t)\big|_{t=0} =\frac{\epsilon}{\epsilon} =1 $. 
\\[1ex]
\noindent (i)
The first observation is that the initial condition of the survival probability is
$\Phi^{(0)}(t)\big|_{t=0} = \epsilon $. The good initial condition $\Phi^{(n)}(0)=\delta_{n0}$
is {\it only} fulfilled for $\epsilon=1$ and as a consequence $w(0)=1-\epsilon =0$.
\\[1ex]
\noindent(ii) The second observation is that for 
$\epsilon \rightarrow 0$ we have $w(0)=1-\epsilon \rightarrow 1$ and hence $\Phi^{(0)}(t)\big|_{t=0} \rightarrow 0 $ for the survival probability. Hence at $t=0$ at least one event already has arrived (almost surely).
\\[1ex]
Consider now a Montroll-Weiss DTRW with this discrete-time counting process.
Then with (\ref{Coxseries}) the generating function of the transition matrix by accounting for Eq. 
(\ref{generfustate}) writes
\begin{equation}
\label{genermontroll}
{\bar {\mathbf P}}(u) = \sum_{n=0}^{\infty}{\bar \Phi}^{(n)}(u){\mathbf H}^n =  
 \frac{1-{\bar w}(u)}{1-u}\left({\mathbf 1}-{\bar w}(u){\mathbf H}\right)^{-1} .
\end{equation}
The `natural' initial condition is ${\mathbf P}(t=0)= (\delta_{ij})$ when the walker at $t=0$ is sitting on a given initial node $i$.
However, this natural initial condition is fulfilled if and only if the state probabilities fulfill the initial condition $\Phi^{(n)}(t=0)=\delta_{n0}$ which is true only for $\epsilon=1$ (See (\ref{inicons})).
Let us now consider the effect of $\epsilon < 1$ on the initial condition of the transition matrix, namely
\begin{align}
\label{Montroll-Weiss-walk-timezero}
\ds 
&{\mathbf P}(t)\big|_{t=0} ={\bar {\mathbf P}}(u)|_{u=0}  = \sum_{n=0}^{\infty}\epsilon(1-\epsilon)^n{\mathbf H}^n \notag \\
& \hspace{1.3cm}= | v_1\rangle\langle {\bar v}_1| \sum_{n=0}^{\infty}\epsilon(1-\epsilon)^n+ 
\epsilon \sum_{m=2}^N  \sum_{n=0}^{\infty} \left((1-\epsilon)\lambda_m\right)^n | v_m\rangle\langle {\bar v}_m| \\
& \ds 
 {\mathbf P}(0) =  | v_1\rangle\langle {\bar v}_1| + \epsilon \sum_{m=2}^N \frac{| v_m\rangle\langle {\bar v}_m|}{1-\lambda_m(1-\epsilon)} 
 = | v_1\rangle\langle {\bar v}_1| + \left(1-w(0)\right)\sum_{m=2}^N \frac{| v_m\rangle\langle {\bar v}_m|}{1-\lambda_m w(0)} . \notag
\end{align}
We observe that for $\epsilon \rightarrow 0$ (i.e. $w(0) \rightarrow 1-0$) the transition matrix already at $t=0$ takes the stationary distribution ${\mathbf P}(0) \rightarrow | v_1\rangle\langle {\bar v}_1| $ (and not as it is in a `good' DTRW for $t\rightarrow \infty$). Then
the departure node of the walker becomes maximally uncertain in the sense of a random initial condition where the walker at $t=0$ (almost surely) makes a huge number of steps to reach `immediately' the stationary equilibrium distribution which remains unchanged
${\mathbf P}(t) = |v_1\rangle\langle {\bar v}_1|$ $\forall t$. 
\\[1ex]
In view of this consideration we can formulate an `uncertainty principle' as follows:
{\it The more ${\mathbb P}(Z=t)=w(t)$ is `localized versus $t=0$', the more variable is the (random) initial condition of the transition matrix}. As a consequence for $w(0) > 0$ the transition matrix does not solve a Cauchy initial value problem and the departure node is `uncertain'.
On the other hand for $\epsilon \rightarrow 1$ (i.e. $w(0) \rightarrow 0$) the initial condition becomes
${\mathbf P}(t=0)= {\mathbf 1}$ thus the DTRW then is a `good walk' with a well-defined (`certain') departure node.
\\[1ex]
We notice that
relation (\ref{Montroll-Weiss-walk-timezero}) is consistent with recent results obtained by a different approach \cite{RiascosBoyerHerringerMateos2019}.
\section{\small CONCLUSIONS}
\label{Conclusions}
In this paper we analyzed discrete-time renewal processes and their continuous-time limits. We focused especially on counting processes which are `Prabhakar'. These processes are discrete-time approximations of the Prabhakar continuous-time renewal process (GFPP). 
Among the discrete-time Prabhakar processes one process (the PDTP) stands out 
as it contains for special choices of parameters (namely for $\nu=1$ and $0<\alpha<1$) the fractional
Bernoulli counting process and (for $\nu=1$ and $\alpha=1$) the standard Bernoulli process. 
The PDTP and the class of discrete-time `Prabhakar' counting processes converge in well-scaled continuous-time limits to the continuous-time Prabhakar process. 
\\[1ex]
The PDTP is constructed such 
that zero waiting times between events are forbidden leading to strictly positive interarrival times. This `good' initial condition makes the PDTP useful to define a Montroll-Weiss DTRW that solves a Cauchy problem with a well-defined departure node of the walker. We called this walk Prabhakar DTRW.
We derived for the Prabhakar DTRW generalized fractional discrete-time Kolmogorov-Feller equations that govern the resulting stochastic dynamics on undirected graphs. The Prabhakar DTRW (unless for $\alpha=1$ and $\nu=1$) is non-Markovian where in the range $\alpha\in (0,1)$ and $\nu>0$ for long observation times universal (Mittag-Leffler) power-law long-time memory effects emerge with fat-tailed waiting time densities.
\\[1ex]
We demonstrate explicitly that for certain choices of parameters the generalized fractional difference Kolmogorov-Feller equations of the Prabhakar DTRW turn into their classical counterparts of fractional Bernoulli and standard Bernoulli, and
the same is true for their continuous-time limits recovering fractional Poisson and Poisson, respectively.
\\[1ex]
In the present paper we analyzed stochastic motions that are defined as normal random walks subordinated to a PDTP and introduced the Prabhakar DTRW.
Generally the subordination approach is a powerful tool to define new stochastic processes. Anomalous diffusive motions where for instance the operational time 
(i.e. the number of jumps performed up to a physical time $t$) grows either slower than the physical time (subdiffusion) or faster than the physical time (superdiffusion) are  of great interest. A general approach to construct such stochastic motions with special emphasis on biased and strictly increasing walks has been developed in a recent article \cite{tmm-fp-apr-fractalfract}.
\\[1ex]
Although we focused in the present paper on undirected graphs, a promising field of applications of the PDTP and Prabhakar DTRW arises from biased walks on directed graphs (See \cite{RiascosMichelitschPizarro2020} for a recent model of fractional dynamics on directed graphs).
Among these cases Prabhakar DTRWs come along as strictly increasing walks with interesting applications such as `aging in complex systems'. These problems exhibit strictly increasing random quantities (`damage-misrepair accumulation') \cite{RiascosWang-Mi-Mi2019}. 
Generally the Prabhakar DTRW approach opens a huge potential of new interdisciplinary applications to topical problems as various as the time evolution of pandemic spread, communication in complex networks, dynamics in public transport networks,
anomalous relaxation, collapse of financial markets, just to denominate a few examples.
\subsection*{Acknowledgments}

F.~Polito has been partially supported by the project ``Memory in Evolving Graphs'' (Compagnia di San Paolo/Università degli Studi di Torino).
%
%
%

\begin{appendix}
\section{\small \it APPENDICES - SUPPLEMENTARY MATERIALS}
\subsection{\small DISCRETE CONVOLUTIONS AND GENERATING FUNCTIONS}
\label{generfuncts}
In this appendix we review some basic properties of generating functions (corresponding to the discrete Laplace transform in the sense soon after explained) which are powerful analytical tools.
The generating function is defined by
\begin{equation}
\label{discerte-Laplace}
{\bar g}(u) = \mathbb{E} u^Z = \sum_{k=0}^{\infty} u^k g(k), \hspace{1cm} |u |\leq 1 
\end{equation}
where we assume $\mathbb{P}(Z=k)= g(k)$ is
a discrete-time waiting time probability density in a discrete-time counting process.
Normalization is expressed by ${\bar g}(u)|_{u=1}=1$, $\mathbb{E} u^Z$ stands for expectation value of $u^{Z}$ and $Z \in \{0,1,2,..\} = \mathbb{N}_0$ a.s.\ indicates a discrete random variable which takes the values 
$Z=k$ with probability $\mathbb{P}(Z=k) = g(k)$. 
We notice the correspondence of generating functions with Laplace transforms by introducing the density (PDF) 
\begin{equation}
\label{densitycorr}
\chi(t) = \sum_{k=0}^{\infty} \delta(t-k) g(k) ,\hspace{1cm} t \in \mathbb{R}
\end{equation}
defined on continuous-time where $\delta(\tau)$ stands for Dirac's $\delta$-distribution. Density (\ref{densitycorr})
has Laplace transform
\begin{equation}\label{Lapfor}
{\cal L}\{\chi(t)\}(s)= \int_{0_{-}}^{\infty}e^{-st} \chi(t) {\rm d}t = \mathbb{E} e^{-sZ} = 
{\bar g}(e^{-s}) ,\hspace{1cm} \Re{s} \geq 0
\end{equation} 
which is the generating function (\ref{discerte-Laplace}) for $u=e^{-s}$ and
${\bar g}(e^{-s})\big|_{s=0}=\sum_{k=0}^{\infty}g(k)=1$ reflects normalization.
The discrete convolution operator of two causal functions $g(t), h(t)$ ($t\in \mathbb{N}_0$) is defined as
\begin{equation}
\label{convol-discrete}
(g \star h)(t) = (h \star g)(t) =\sum_{j=0}^tg(j)h(t-j) = \frac{1}{t!}\frac{d^t}{d u^t} \{{\bar g}(u){\bar h}(u)\}|_{u=0} ,\hspace{1cm} t \in \mathbb{N}_0
\end{equation}
where
$ {\bar g}(u){\bar h}(u) = \sum_{t=0}^{\infty} (g \star h)(t) \, u^t $ is the 
generating function of the discrete convolution. 
We use sometimes also the synonymous notation $g(t) \star h(t)$ for (\ref{convol-discrete}).
In the following we denote the $n$th convolution power with 
\begin{equation} 
\label{conv-n-power}
(g(t) \star)^n =  \frac{1}{t!}\frac{d^t}{d u^t} ({\bar g}(u))^n|_{u=0} ,\hspace{0.5cm} t, n \in \mathbb{N}_0
\end{equation}
which includes $n=0$ and yields unity $(g(t) \star)^0 =\delta_{t0}$ (where $\delta_{ij}$ denotes the Kronecker symbol, which in the following will also be denoted by the equivalent notation $\delta_{i,j}$).
Relation (\ref{conv-n-power}) can be extended to non-integer (especially fractional) convolution powers $n\in \mathbb{R}$. 
\\[1ex]
Now let us see the connection between generating functions and the shift operator ${\hat T}_{-1}$ which is defined by 
${\hat T}_{-1}g(t)=g(t-1)$. Consider again generating function (\ref{discerte-Laplace})
and replace $u$ 
by the shift operator ${\hat T}_{-1}$ and use $({\hat T}_{-1})^n\delta_{t,0} ={\hat T}_{-n}\delta_{t0} =\delta_{(t-n),0}=\delta_{tn}$
to arrive at
\begin{equation}
\label{shift-gen}
{\bar g}({\hat T}_{-1}) \delta_{t0} = 
\sum_{k=0}^{\infty}g(k) \delta_{tk} = g(t) = \frac{1}{t!} \frac{d^t}{du^t}{\bar g}(u)|_{u=0} ,\hspace{0.5cm} t\in \mathbb{N}_0 .
\end{equation}
In the present paper we extensively make use of the correspondence of generating functions and shift-operator counterparts. Such procedures turn out to be especially useful in the definition and determination of suitably scaled continuous-time limits and governing discrete-time evolution equations.
\subsection{\small SHIFT OPERATORS AND CAUSAL DISTRIBUTIONS}
\label{crucialdefinitions}
We consider a {\it causal} distribution $\Theta(t)g(t)$ having all its non-zero values on $t\in \{0,h,2h,\ldots\}=h\mathbb{N}_0$ ($h>0$) and zero values for $t<0$. For its definition we make use of the {\it discrete-time
Heaviside function} defined by
\begin{equation}
\label{heavisideintegers}
\ds
\Theta_h(t) = \Theta(t)= \left\{\begin{array}{l} 1 , \hspace{1cm} t\in \{ 0,h,2h,\ldots \hspace{0.2cm}\}\\ \\ 
                    0 ,\hspace{1cm}    t \in \{-h,-2h,\ldots \}
                   \end{array}\right. \hspace{1cm} t \in h\mathbb{Z}, \hspace{1cm} h>0.
                   \end{equation}
We especially emphasize that with this definition $\Theta_h(0)=1$. In the entire paper for Heaviside functions we write $\Theta_h(t)=\Theta(t)$, i.e.\ we may skip subscript $h$. 
Without loss of generality we may identify the discrete-time Heaviside function with its continuous-time counterpart, i.e. the conventional Heaviside- unit-step function with $\Theta(t)=1$ for $t\geq 0$ (especially $\Theta(0)=1$) and $\Theta(t)=0$
for $t<0$ defined on $t\in \mathbb{R}$. This is necessary when we use
${\hat T}_{-h} = e^{-hD_t}$ leading to the `distributional representation' 
$e^{-hD_t}\Theta(t)=\Theta(t-h)$ which is defined on $t\in \mathbb{R}$.
This allows to define Laplace transformation of $\delta_h(t)$ as in 
(\ref{Lapladeltah}) which has `good properties' in the limit for $h\rightarrow 0$.
\\[1ex]
Let $\delta_{k,l}$ (we also use notation $\delta_{kl}$) be the Kronecker symbol defined for {\it any positive and negative integer including zero} by
\begin{equation}
\label{Kroneckerdelta}
\ds 
\delta_{i,j} = \delta_{i-j,0} =  \left\{\begin{array}{l} 1 , \hspace{1cm} i=j  \\ \\ 
                    0,     \hspace{1cm}   i \neq j
                   \end{array}\right. \hspace{2cm} i,j \in \mathbb{Z}
                   \end{equation}
where we have the translational invariance property $\delta_{i+m,j+m}=\delta_{i,j}$ ($m\in \mathbb{Z}$).
Then we define the `{\it discrete-time $\delta$-distribution}' as follows
\begin{equation}
\label{discrete-time-delta}
\ds \delta_h(t) = \frac{\Theta(t)-\Theta(t-h)}{h} = 
\frac{1-{\hat T}_{-h}}{h}\Theta(t) = \,\, \frac{1}{h}
\delta_{\frac{t}{h},0} \,\,  = \ds \left\{\begin{array}{l}\ds  \frac{1}{h} ,\hspace{0.5cm} t=0 \\ \\ \ds
0 ,\hspace{0.6cm} t\neq 0  \end{array}\right.
\hspace{1cm} t \in h\mathbb{Z}
\end{equation}
where $\delta_{\frac{t}{h},0}$ in this relation indicates the Kronecker symbol (\ref{Kroneckerdelta}) where we notice that $\frac{t}{h} \in \mathbb{Z}$. We observe that for $h=1$ we have $\delta_1(t)=\delta_{t0}$.
Formula (\ref{discrete-time-delta}) defines a discrete-time density on
$t\in h\mathbb{Z}$. However, it makes sense to extend this definition to $t \in \mathbb{R}$. With this extended definition $\delta_h(t) = \frac{\Theta(t)-\Theta(t-h)}{h}$ becomes an integrable distribution in the Gelfand-Shilov sense \cite{GelfandShilov1968} with $\delta_h(t)= \frac{1}{h}$ for $t\in [0,h)$ and $\delta_h(t)=0$ else.
Then we have $\int_{-\infty}^{\infty}\delta_h(\tau){\rm d}\tau=\int_0^h\delta_h(\tau){\rm d}\tau=1$ and $\delta_h(t)$ has a well-defined Laplace transform (subsequent relation (\ref{Lapladeltah})).
Of great importance is 
the continuous-time limit of the discrete-time $\delta$-distribution 
\begin{equation}
\label{dirac-delta}
\lim_{h\rightarrow 0} \delta_h(t) = \lim_{h\rightarrow 0}\frac{1-e^{-hD_t}}{h} \Theta(t) = D_t\Theta(t) = \delta(t) ,\hspace{1cm} t \in \mathbb{R}
\end{equation}
recovering Dirac's continuous-time $\delta$-distribution.
Further, it appears also useful to briefly consider the Laplace transform of $\delta_h(t)$, namely
\begin{equation}
\label{Lapladeltah}
\begin{array}{l} \ds
{\cal L}\{\delta_h(t)\}(s) = \int_{0_{-}}^{\infty} \delta_h(t) e^{-st}{\rm d}t = \int_{-\infty}^{\infty} e^{-st} \frac{(1-e^{-hD_t})}{h}\Theta(t){\rm d}t \\ \\ \ds
\hspace{1cm} = \int_{-\infty}^{\infty} \Theta(t) \frac{(1-e^{+hD_t})}{h} e^{-st} {\rm d}t =\frac{1-e^{-hs}}{h} \int_0^{\infty}e^{-st} {\rm d}t \hspace{1cm} \Re\{s\} >0\\ \\ \ds \hspace{1cm}
 = \frac{1-e^{-hs}}{h s} 
\end{array}
\end{equation}
where indeed $\lim_{h\rightarrow 0} {\cal L}\{\delta_h(t)\}(s) = {\cal L}\{\delta(t)\}(s) = 1$ recovers the Laplace transform of Dirac's $\delta$-distribution 
reflecting (\ref{dirac-delta}).
We applied here partial integration and we used ${\hat T}_{-h}= e^{-hD_t}$ having adjoint operator $ ({\hat T}_{-h})^{\dagger}={\hat T}_{+h}=e^{+hD_t}$ (with $e^{hD_t}e^{-st} =e^{-hs}e^{-st}$).
\\ \\ \\
\noindent {\it Continuous-time limit of discrete-time densities}
\\ \\
Let $g(k)$ be a ‘discrete-time density’ such as for instance
a discrete-time waiting time distribution $\mathbb{P}(Z=k)=g(k,C)$ ($k \in \mathbb{N}_0$) and $C$ stands for an internal parameter (or a set of parameters). We often write $g(t)$ instead of $g(t,C)$ if only the time dependence is of interest. We assume $g(t)=\Theta(t) g(t)$ ($t\in \mathbb{Z}$) to be a causal distribution.
Then we can define the `{\it scaled causal discrete-time density}' by accounting for (\ref{discrete-time-delta}) 
as follows
\begin{equation}
\label{thepropertyogdelta}
\begin{array}{l} \ds
g(t)_h = {\bar g}({\hat T}_{-h}) \delta_h(t)= 
 \sum_{k=0}^{\infty} g(k,C){\hat T}_{-hk} \delta_h(t) = 
 \sum_{k=0}^{\infty} g(k,C)\delta_h(t-kh)  \\ \\ \ds \hspace{1cm} = \Theta(t) h^{-1} g\left(\frac{t}{h},C\right) ,
\hspace{1cm} t \in h \mathbb{Z} \\ \\ \ds
\hspace{1cm} = h^{-1} g\left(\frac{t}{h},C\right) ,\hspace{1cm} t \in h \mathbb{N}_0
\end{array}
\end{equation}
where due to the summation over non-negative $k$ (\ref{thepropertyogdelta}) is non-zero only for $t\geq 0$ (causality).
Here the process takes place on $t\in h \mathbb{N}_0$ and the scaled discrete-time density $g(t)_h=h^{-1} g(\frac{t}{h},C)$ ($t\in h\mathbb{N}_0$) as well as $g(t)_{h=1} = g(t,C)$ ($t\in \mathbb{N}_0$) both
have physical dimension $sec^{-1}$ and merit therefore to be called (discrete-time) {\it densities}. 
Note that $g(\frac{t}{h},C)$ is dimensionless since $\frac{t}{h}$ is dimensionless as $t$ and $h$ have units $\sec$.
When $h$ is finite the parameter $C$ does not necessarily need to be scaled with $h$.
However, we will see subsequently that generally we need to introduce a scaling law $C=C(h)=C_0h^{\rho}$ such that the continuous-time limit $h\rightarrow 0$ of (\ref{thepropertyogdelta}) exists. The constant $C_0>0$ is independent of $h$ and is a free parameter. $C_0$ has units $\sec^{-\rho}$ (as $C=C_0h^{\rho}$ is dimensionless).
\\ [1ex]
Now let us define the continuous-time limit density of (\ref{thepropertyogdelta}) and introduce the variable $\tau_k= hk \in h\mathbb{N}_0$ and rescaled parameter $C(h)=C_0h^{\rho}$ with scaling dimension $\rho$ such that the limit
\begin{equation}
\label{thepropertyogdelta-contilimit}
\begin{array}{l} \ds 
 g(t,C_0)_{ct} =
\lim_{h\rightarrow 0} {\bar g}({\hat T}_{-h}) \delta_h(t)= 
\lim_{h\rightarrow 0} \sum_{k=0}^{\infty} \frac{1}{h}g\left(\frac{\tau_k}{h},C_0h^{\rho}\right) \delta_h(t-hk) h  = \int_0^{\infty}  g(\tau,C_0)_{ct}\delta(t-\tau){\rm d}\tau  \\ \\ \ds
\hspace{0.5cm} =\lim_{h\rightarrow 0} \sum_{k=0}^{\infty} g\left(k,C_0h^{\rho}\right) \frac{1}{h}\delta_{\frac{t}{h},k} =  \lim_{h\rightarrow 0} \frac{1}{h} g\left(\frac{t}{h},C_0h^{\rho}\right) 
\end{array}
 \end{equation}
exists. In the second line we used $\delta_h(t-kh)=\frac{1}{h}\delta_{\frac{t}{h},k}$ of (\ref{discrete-time-delta}).
We call $g(t)_{ct}=g(t,C_0)_{ct}$ ($ t \in \mathbb{R}^{+}$) the continuous-time limit density of the discrete-time density $g(t,C)$ ($ t \in \mathbb{N}_0$). We refer $\frac{1}{h} g(\frac{t}{h},C)$ ($t =\{0,h,2h,\ldots, \}\in h\mathbb{N}_0 $) to as the {\it `scaled' discrete-time density} and $\frac{1}{h} g(\frac{t}{h},C_0h^{\rho})$ the {\it `well-scaled' discrete-time density}, the latter with suitably chosen scaling dimension $\rho$ such that the continuous-time limit exists. We further notice that $C_0^{-\frac{1}{\rho}}$ defines a characteristic time scale in the continuous-time density $g(t,C_0)_{ct}$.
\\[1ex]
Limiting relation (\ref{thepropertyogdelta-contilimit}) indeed is useful when we consider 
processes on discrete times $\in h\mathbb{N}_0$ and especially for the transition
to the continuous-time limit $h\rightarrow 0$ where the asymptotic behavior of $g(k)$ with $k=\frac{t}{h}$ large comes into play.
To see this we consider the important case of a fat-tailed behavior $g(k) = g(k,C) \sim C k^{-\alpha-1}$ for large $k$ with $\alpha \in (0,1)$ such as for instance in relation (\ref{appro}).
Then we need to rescale the internal parameter $C(h)=C_0h^{\rho}$ such that this limit 
\begin{equation}
\label{limitexist}
g(t,C_0)_{ct}=\lim_{h\rightarrow 0}= h^{-1} g\left(\frac{t}{h},C_0h^{\rho}\right) \lim_{h\rightarrow 0} h^{-1} C_0h^{\rho} \left(\frac{t}{h}\right)^{-\alpha-1} 
= \lim_{h\rightarrow 0} C_0 h^{\rho+\alpha} t^{-\alpha-1}   ,\hspace{0.5cm} t \in \mathbb{R}
\end{equation}  
exists. Clearly this limit exists if and only if $\rho+\alpha=0$ and hence $\rho=-\alpha$ thus the required scaling is $C(h)=C_0h^{-\alpha}$. Then the continuous-time limit yields
\begin{equation}
\label{contilimitexistsscaling}
g(t,C_0)_{ct} = \lim_{h\rightarrow 0} h^{-1} g\left(\frac{t}{h},C_0h^{-\alpha}\right) = C_0 t^{-\alpha-1} .
\end{equation}
We notice that a fat-tailed continuous-time limit 
of the form (\ref{contilimitexistsscaling}) originates from a discrete-time density having
generating function ${\bar g}(u)=1-C'(1-u)^{\alpha} +0(1-u)^{\alpha}$ (with $C'=C\frac{\Gamma(1-\alpha)}{\alpha} >0$)
which corresponds to a rescaled version of Sibuya$(\alpha)$. 
This underlines the utmost universal importance of Sibuya$(\alpha)$ which is recalled in more details in Appendices \ref{Sib}-\ref{Sibuyarenewal}.
Furthermore, to establish the connection to relation 
(\ref{appro}): When we identify the constants $C(h)=\frac{\alpha\nu}{\xi(h)\Gamma(1-\alpha)}$ (and $C_0=\frac{\nu}{\xi_0}\frac{\alpha}{\Gamma(1-\alpha)}$) then relations
(\ref{limitexist}), (\ref{contilimitexistsscaling}) recover (\ref{appro}). 
\\[1ex]
Well-scaled limiting procedures in the sense of 
relations (\ref{thepropertyogdelta-contilimit})-(\ref{contilimitexistsscaling}) are recurrently performed in the present paper.
\\ \\ \\
{\it Continuous-time limit of `cumulative discrete-time distributions'}
\\ \\
Apart from densities we also deal with `cumulative discrete-time distributions' which we define as
\begin{align}
\label{cumul} 
G(t) & = G(t)_{h=1} = \frac{1}{t!}\frac{d^t}{du^t}\left(\frac{{\bar g}(u)}{1-u}\right)|_{u=0}  = \sum_{k=0}^{t} g(t-k)\Theta(k) \notag\\
& = \Theta(t) \star g(t) = \frac{1}{1-{\hat T}_{-1}} g(t) = \sum_{k=0}^t g(k) , \hspace{1cm} t \in \mathbb{N}_0
\end{align}
where $g(t)=g(t,C)$ is a causal {\it discrete-time density} as defined above in (\ref{thepropertyogdelta}). Causality makes the upper limit $k=t$
in the sum $\frac{1}{1-{\hat T}_{-1}} g(t) = \sum_{k=0}^{\infty} {\hat T}_{-k}  \Theta(t)g(t) =
\sum_{k=0}^{\infty} \Theta(t-k)g(t-k)$.
Cumulative distributions in the sense of (\ref{cumul}) are distributions generated by summation of discrete-time 
densities and are dimensionless distributions.
Examples of cumulative discrete-time distributions include the (dimensionless)
state probabilities $\Phi_{\alpha,\nu}^{(n)}(t)$ (See (\ref{state-prob})).
As cumulative distributions are dimensionless, we have to think the sum (\ref{cumul}) multiplied with a time increments
$h=1$ (having physical dimension $sec$).
With these observations we define the `{\it cumulative scaled discrete-time distribution}' as
\begin{equation}
\label{cumuldiscrete-time}
G(t)_h = \frac{h}{1-{\hat T}_{-h}} 
{\bar g}({\hat T}_{-h}) \delta_h(t)  = {\bar g}({\hat T}_{-h})\frac{h}{1-{\hat T}_{-h}}\delta_h(t) = {\bar g}({\hat T}_{-h})\Theta(t),  \hspace{0.5cm} t \in h\mathbb{Z}
\end{equation}
where
$\frac{h}{1-{\hat T}_{-h}}\delta_h(t) =\Theta(t)$ (See (\ref{discrete-time-delta})) and $G(t)_{h=1}$ recovers (\ref{cumul}). 
Now in order to consider the continuous-time limit we rewrite (\ref{cumuldiscrete-time}) as
\begin{equation}
\label{cumultimelimit}
\begin{array}{l} \ds 
G(t)_h= \frac{h}{1-{\hat T}_{-h}} g(t)_h = \sum_{k=0}^{\infty} h g(t-kh)_h 
= \sum_{k=0}^{\infty} g(\frac{t-kh}{h},C) \\ \\ \ds \hspace{0.5cm}
=\sum_{k=0}^{\frac{t}{h}}g(k,C) = G\left(\frac{t}{h},C\right)
 \end{array}  \hspace{0.5cm} t \in h\mathbb{N}_0
\end{equation}
where $G(k,C)=G(k,C)_1$ ($k=\frac{t}{h}\in \mathbb{N}_0$) is the cumulative distribution (\ref{cumul}) and $g(t)_h = {\bar g}({\hat T}_{-h}) \delta_h(t)$ indicates a causal
scaled density as in (\ref{thepropertyogdelta}).
The continuous-time limit then is obtained as
\begin{equation}
\label{contilimim-cmul}
\begin{array}{l} 
\ds G(t,C_0)_{ct} = \lim_{h\rightarrow 0} G(t)_h  = \lim_{h\rightarrow 0} G\left(\frac{t}{h},C_0h^{\rho}\right)_1 \\ \\ \ds \hspace{0.5cm} = 
D_{t}^{-1} \,\, g(t,C_0)_{ct} \,\, =
\int_0^t  g(\tau,C_0)_{ct}{\rm d}\tau 
\end{array}
\hspace{0.5cm} t \in \mathbb{R}^+ .
\end{equation}
The scaling of the constant(s) $C(h)=C_0h^{\rho}$ has to be suitably chosen such that the limit (\ref{contilimim-cmul}) exists where the same scaling dimension $\rho$ occurs as for the density in (\ref{thepropertyogdelta-contilimit}).
This equation is constructive for the explicit evaluation of continuous-time limits 
for instance in the state probabilities
since it is the limit of the integer time expression $G(\frac{t}{h},C_0h^{\rho})_1$ ($\frac{t}{h} \in \mathbb{N}$ large) of Eq. (\ref{cumul}) for $h\rightarrow 0$.
\\ \\ \\
\noindent {\it Further properties of causal discrete-time distributions}
\\ \\
Now let us consider discrete-time convolutions of causal distributions $G(t) = \Theta(t) g(t)$. By accounting for the shift operator ${\hat T}_{-1}$ such that ${\hat T}_{-n}g(t)=g(t-n)$, especially we define
\begin{equation}
\label{defineshifts}
\begin{array}{l}
{\hat T}_{-n} \delta_{i,j} = \delta_{i-n,j} = \delta_{i,j+n} = \delta_{i,j} {\hat T}_{+n}  \\ \\
{\hat T}_{n} = ({\hat T}_{-n})^{\dagger} =({\hat T}_{-n})^{-1} 
\end{array}
\end{equation}
where with ${\hat A}^{\dagger}$ we indicate the adjoint (hermitian conjugate) operator to ${\hat A}$. It is well-known that
shift operators are unitary and represent a commutative (Abelian) group. From observation (\ref{defineshifts})
we introduce the notation
\begin{equation}
\label{defineshift}
\begin{array}{l}
{\hat T}_a \delta_{i,j} = \delta_{i+a,j} = \delta_{i,j-a} = \delta_{i,j}{\hat T}_{-a}   , \\ \\
\delta_{i,j} {\hat T}_b = \delta_{i,j+b} = \delta_{i-b,j} = {\hat T}_{-b}  \delta_{i,j} .
\end{array}
\end{equation}
We introduce the convention that if we write the shift operator ${\hat T}$ to the left of the Kronecker symbol $\delta_{i,j}$ 
it acts on the left index $(i,..)$
and if we write it to the right of the Kronecker symbol, then it acts on right index 
$(..,j)$.
Let us consider now a chain of discrete-time convolutions
\begin{align}
\label{disrep}
 (G_1 & \star G_2 \star \ldots \star G_n )(t) = {\bar g}_1({\hat T}_{-1}){\bar g}_2({\hat T}_{-1})\ldots {\bar g}_n({\hat T}_{-1}) \, \delta_{t,0} \notag\\
 &= \sum_{t_1=-\infty}^{\infty} \sum_{t_2=-\infty}^{\infty}\ldots \sum_{t_n=-\infty}^{\infty} 
\delta_{t,t_1+t_2+\ldots +t_n} \Theta(t_1)\Theta(t_2)\ldots \Theta(t_n) g_1(t_1)g_2(t_2) \ldots g_n(t_n) ,\hspace{0.5cm} t \in \mathbb{Z}
\end{align}
where the first line establishes the shift-operator representation. \\ We notice in view of
${\bar g}_i({\hat T}_{-1}){\bar g}_j({\hat T}_{-1}) = {\bar g}_j({\hat T}_{-1}){\bar g}_i({\hat T}_{-1})$ that discrete-time convolutions (as their continuous-time counterparts) indeed are commuting.
The convolution power (\ref{conv-n-power}) is a special case of (\ref{disrep}). The generating function of (\ref{disrep}) can be recovered by
\begin{equation}
\label{generforeobt}
\begin{array}{l}
\ds 
\sum_{t=0}^{\infty} (G_1 \star G_2 \star \ldots \star G_n )(t)u^t  \\ \\\ds  =  \sum_{t_1=-\infty}^{\infty} \sum_{t_2=-\infty}^{\infty}\ldots \sum_{t_n=-\infty}^{\infty}
\Theta(t_1)\Theta(t_2)\ldots \Theta(t_n) g_1(t_1)g_2(t_2) \ldots g_n(t_n) \times \sum_{t=0}^{\infty} \delta_{t,t_1+t_2+\ldots +t_n} u^t \\ \\ \ds
= \prod_{j=1}^n \sum_{t_j=0}^{\infty}  g_j(t_j) u^{t_j} = {\bar g}_1(u){\bar g}_2(u)\ldots {\bar g}_n(u) 
\end{array}
\end{equation}
reflecting also shift-operator representation in (\ref{disrep}).
Then define operator functions $A({\hat T}_{-1})$, $B({\hat T}_{+1})$ which commute and act on a convolution of the form
(\ref{generforeobt}) 
\begin{equation}
\label{convol}
\begin{array}{l}
\ds 
A({\hat T}_{-1}) B ({\hat T}_{+1}) (G_1 \star G_2 \star \ldots \star G_n )(t) 
\\ \\ \ds
=  \sum_{t_2=-\infty}^{\infty}\ldots \sum_{t_n=-\infty}^{\infty}
\Theta(t_1)\Theta(t_2)\ldots \Theta(t_n) g_1(t_1)g_2(t_2) \ldots g_n(t_n) A({\hat T}_{-1})  B ({\hat T}_{+1})  \delta_{t,t_1+t_2+\ldots +t_n}  \\ \\
\ds 
\sum_{t_2=-\infty}^{\infty}\ldots \sum_{t_n=-\infty}^{\infty}
\Theta(t_1)\Theta(t_2)\ldots \Theta(t_n) g_1(t_1)g_2(t_2) \ldots g_n(t_n) \delta_{t,t_1+t_2+\ldots +t_n} A({\hat T}_{+1})  B ({\hat T}_{-1})
\\ \\ \ds 
=  (G_1 \star G_2 \star \ldots \star G_n )(t) B ({\hat T}_{-1}) A({\hat T}_{+1})
\end{array}
\end{equation}
where it is crucial here that $t \in \mathbb{Z}$.
All the representations in (\ref{convol}) are equivalent where we utilized (\ref{defineshift}) with the adjoint 
operator function $\left\{A({\hat T}_{-1}) B({\hat T}_{+1})\right\}^{\dagger}
= B({\hat T}_{-1})  A({\hat T}_{+1})$ (commuting). For instance consider $A({\hat T}_{-1}) ={\hat T}_{-n}$, then
we have
\begin{equation}
\label{shiftex}
\begin{array}{l}
\ds
{\hat T}_{-n} (G_1 \star G_2 \star \ldots \star G_n )(t)= (G_1 \star G_2 \star \ldots \star G_n )(t-n) \\ \\ \ds =
\sum_{t_2=-\infty}^{\infty}\ldots \sum_{t_n=-\infty}^{\infty}
\Theta(t_1)\Theta(t_2)\ldots \Theta(t_n) g_1(t_1)g_2(t_2) \ldots g_n(t_n) \delta_{t-n,t_1+t_2+\ldots +t_n} \\ \\ 
 \ds =
(G_1 \star G_2 \star \ldots \star G_n )(t) {\hat T}_{+n} .
\end{array}
\end{equation}
Application of any shift-operator functions on discrete-time convolutions of causal distributions is hence well-defined in the sense 
(\ref{convol}) by reducing the shift operations to shift operations in the argument of the Kronecker symbol.
In the distributional sense
all operations involving operator functions of 
shift operators
act on convolutions (for instance the discrete-time fractional Kolmogorov-Feller equations 
(\ref{recursive-diff}), (\ref{explicitbackward}) and (\ref{alsointheform})) and 
are well-defined as they can be reduced to shift operations acting on Kronecker-$\delta$s defined on $\mathbb{Z}$.
It is straightforward to rewrite all these relations for their scaled counterparts defined on $h\mathbb{Z}$ involving discrete-time $\delta$-distributions (\ref{discrete-time-delta}) generalizing the Kronecker-$\delta$s. We can conceive the causal discrete-time distributions in a wider sense as discrete-time counterparts of 
Gelfand-Shilov generalized functions and distributions \cite{GelfandShilov1968}.

\subsection{\small GENERALIZED FRACTIONAL INTEGRALS AND DERIVATIVES}
\label{appendix1}
Here we consider some properties of discrete-time convolution operators and their continuous-time limits in order to establish the connection with fractional calculus. To this end let us introduce the notion of `generalized discrete-time fractional integrals and derivatives' in the spirit of Ref. \cite{michelCFM2011}.
Let $\Theta(t)f(t)$ be a causal discrete-time distribution defined on $t\in h\mathbb{Z}$ ($h>0$). Then we define the well-scaled `{\it discrete-time integral}' as
\begin{equation}
\label{discrte-time-convol}
 \frac{h}{1-{\hat T}_{-h}}\Theta(t)f(t) = h\sum_{k=0}^{\infty} {\hat T}_{-hk} \Theta(t)f(t)
= \sum_{k=0}^{\frac{t}{h}} h f(t-kh)= \sum_{k=0}^{\frac{t}{h}} hf(kh) ,\hspace{1cm} t\in h \mathbb{Z}
\end{equation}
which is well-scaled thus its continuous-time limit for $h\rightarrow 0$ exists for sufficiently `good', i.e. absolutely
integrable distributions $f(t)$ defined on $t\in \mathbb{R}^{+}$. 
The continuous-time limit then yields the standard integral

\begin{equation}
\label{casesimple}
\lim_{h\rightarrow 0} \frac{h}{(1-e^{-hD_t})} \Theta(t)f(t) = D_t^{-1}\Theta(t)f(t) = \int_0^tf(\tau){\rm d}\tau ,\hspace{1cm} t\in \mathbb{R}^{+} .
\end{equation}
Now we introduce `generalized derivatives and integrals' which include traditional and fractional derivatives and integrals. 
We generalize (\ref{discrte-time-convol}) by introducing the well-scaled
`{\it discrete-time fractional integral}' as
\begin{equation}
\label{renormalization}
\begin{array}{l}
\ds
 \left(\frac{h}{(1-{\hat T}_{-h})}\right)^{\mu} \cdot \Theta(t) f(t) =
h^{\mu} \sum_{k=0}^{\infty} (-1)^k
\left(\begin{array}{l} -\mu \\ \hspace{0.2cm} k\end{array}\right)f(t-hk)\Theta(t-kh) \\ \\ \ds \hspace{4.2cm} =
h^{\mu} \sum_{k=0}^{\frac{t}{h}} (-1)^k\left(\begin{array}{l} -\mu \\ \hspace{0.2cm} k\end{array}\right)f(t-hk)
\end{array} \hspace{0.1cm} \mu>0, \hspace{0.1cm} t \in h\mathbb{N}_0 .
\end{equation}
Taking into account that the shift operator can be represented as
\begin{equation}
\label{shiftasexponential}
{\hat T}_{-h}= e^{-hD_t} ,\hspace{1cm}  D_t=\frac{d}{dt}
\end{equation}
where $\lim_{h\rightarrow 0}h^{-1}(1-{\hat T}_{-h})\cdot f(t)$ 
denoting the traditional derivative of first order. Then we consider the limit 
\begin{equation}
\label{renormalization2}
\lim_{h\rightarrow 0} h^{\mu} (1-e^{-hD_t})^{-\mu}  f(t) =  D_t^{-\mu}\cdot f(t) 
\end{equation}
which defines first of all an operator representation of the integral of order $\mu >0$.
Now we evaluate explicitly the continuous-time limit of (\ref{renormalization}) to arrive at
\begin{equation}
\label{limitprocess}
\begin{array}{l}
 \ds \hspace{1cm} 
D_t^{-\mu} \cdot f(t) := \lim_{h\rightarrow 0} h^{\mu}(1-e^{-hD_t})^{-\mu} \cdot f(t) 
= \lim_{h\rightarrow 0} \sum_{k=0}^{\frac{t}{h}}
 \frac{h^{\mu}}{\Gamma(\mu)} \frac{\Gamma(\mu+k)}{\Gamma(k+1)}  f(t-kh)  \\ \\
\hspace{2.5cm} \ds = \sum_{k=0}^{\ceil{th^{-\delta}}} \frac{h^{\mu}}{\Gamma(\mu)} \frac{\Gamma(\mu+k)}{\Gamma(k+1)} f(t-kh) +   \sum_{k= 1+\ceil{th^{-\delta}}}^{th^{-1}} \frac{h^{\mu}}{\Gamma(\mu)} \frac{\Gamma(\mu+k)}{\Gamma(k+1)} f(t-kh)  \\ \\
\hspace{2.5cm} \ds    
=\lim_{h\rightarrow 0}  \sum_{k=1+\ceil{th^{-\delta}}}^{\frac{t}{h}} h 
\frac{(kh)^{\mu-1}}{\Gamma(\mu)}f(t-kh) + O(h^{(1-\delta)\mu+2\delta}) \\ \\
\hspace{2.5cm} \ds  = \lim_{h\rightarrow 0} \int_{h(1+\ceil{th^{-\delta}})}^t  \frac{\tau^{\mu-1}}{\Gamma(\mu)}f(t-\tau) {\rm d}\tau +O(h^{(1-\delta)\mu+2\delta}) \\ \\ \ds \hspace{2.5cm} = \int_0^t 
\frac{\tau^{\mu-1}}{\Gamma(\mu)}f(t-\tau) {\rm d}\tau .
\end{array}\hspace{-2cm} \delta \in (0,1)
\end{equation}
We indeed identify the last line with the {\it Riemann-Liouville fractional integral} of order $\mu$ ($\mu >0$) \cite{SamkoKilbasMarichev1993}  where 
$\mu \in \mathbb{N}$ recovers traditional integer order 
integration and for $\mu\rightarrow 0+$
the kernel $\frac{\tau^{\mu-1}}{\Gamma(\mu)} \rightarrow \delta(\tau)$ thus reproducing the identity $D_t^{0+}f(t)=f(t)$. 
For $\mu>0$ the kernel $\frac{\tau^{\mu-1}}{\Gamma(\mu)}$ is on $\mathbb{R}^{+}$ a locally integrable function.
\\[1ex]
In the deduction (\ref{limitprocess}) we employed the ceiling function $\ceil{..}$ (see the text below Eq. (\ref{further-important}) for a definition) namely $k_0=\ceil{th^{-\delta}}$ where we have $ \lim_{h\rightarrow 0} h(1+\ceil{th^{-\delta}}) \sim t h^{1-\delta} \rightarrow 0$.
We choose $\delta \in (0,1)$ such that the first sum $S_1(h)$ in the second line 
tends to zero for $h\rightarrow 0$ at least as $S_1(h) \sim h^{\mu(1-\delta)+2\delta}$ which is shown below. The crucial point for this choice is that
$hk_0(h) \sim t h^{1-\delta} \rightarrow 0$, however $k_0(h) \sim t h^{-\delta} \rightarrow \infty$ for $h\rightarrow 0$.
Hence in the second sum of the second line in all terms $k$ is large; so we can employ the asymptotic relation $\frac{\Gamma(k+\mu)}{\Gamma(\mu)\Gamma(k+1)}
\sim \frac{k^{\mu-1}}{\Gamma(\mu)}$.
The scaling behavior of $S_1(h)$ for $h\rightarrow 0$ is then obtained as follows
\begin{equation}
\label{sumestimate}
\begin{array}{l} \ds 
S_1(h) = \sum_{k=0}^{k_0} \frac{h^{\mu}}{\Gamma(\mu)} \frac{\Gamma(\mu+k)}{\Gamma(k+1)} f(t-kh) \sim f(t-\eta k_0h) h^{\mu} \sum_{k=0}^{k_0}\left(\begin{array}{l} -\mu \\ \hspace{0.2cm} k\end{array}\right)(-1)^k  ,\hspace{0.5cm} k_0(h)=\ceil{th^{-\delta}}  \\ \\ 
\ds \hspace{0.1cm} =f(t-\eta h k_0) h^{\mu}\frac{1}{k_0!}\frac{d^{k_0}}{du^{k_0}}(1-u)^{1-\mu}\big|_{u=0} =
f(t-\eta h k_0) h^{\mu}(-1)^{k_0} \left(\begin{array}{l}1 -\mu \\ \hspace{0.2cm} k_0\end{array}\right) ,
\hspace{0.5cm} \eta \in (0,1)
\\ \\ \ds \hspace{0.1cm}
= f(t-\eta h k_0) \frac{h^{\mu}\Gamma(\mu+k_0-1)}{\Gamma(k_0+1)\Gamma(\mu-1)}  = f(t-\eta h k_0) h^{\mu}\frac{(\mu-1)_{k_0}}{k_0!}  \\ \\ \ds S_1(h)
\sim f(t-\eta h k_0) \frac{h^{\mu}k_0^{\mu-2}}{\Gamma(\mu-1)}  \sim 
f(t) \frac{t^{\mu-2}}{\Gamma(\mu-1)} h^{\mu} h^{-\delta(\mu-2)} 
\sim h^{\mu(1-\delta)+2\delta}  \rightarrow 0 ,\hspace{0.5cm} (h\rightarrow 0)
\end{array}
\end{equation}
where in (\ref{sumestimate}) $\mu \neq 1$. 
The case $\mu=1$ of standard integration is considered above (See (\ref{discrte-time-convol})-(\ref{casesimple})).
In the last line we used that $k_0=\ceil{th^{-\delta}} \sim th^{-\delta}$ (large) with 
$\delta\in (0,1)$ and $\mu>0$ and the interim value $f(t-\eta h k_0) \sim f(t) $ is approximately constant over the (small) interval $[t-k_0h,t]$ especially for $h\rightarrow 0$.
\\[1ex]
A more sophisticated generalized integral operator which we refer to as `Mittag-Leffler integral' is of interest \cite{KilbasSaigoSaxena2004}.
We can define this convolution operator by the following well-scaled continuous-time limiting procedure
\begin{equation}
\label{limiting-process}
{\cal I}_{\alpha}\cdot f(t) = 
\lim_{h\rightarrow 0}\frac{\xi(h)}{\xi(h)+(1-{\hat T}_{-h})^{\mu}} \cdot f(t) ,\hspace{0.5cm} \mu >0 ,\hspace{0.5cm} t \in h\mathbb{N}_0
\end{equation}
where $\xi(h) =\xi_0h^{\mu}$ is chosen such that (\ref{limiting-process}) exists. Then (\ref{limiting-process}) has the Laplace transform
\begin{equation}
\label{Lapla}
\begin{array}{l} \ds 
\lim_{h\rightarrow 0} {\cal L}\{I_{\mu}\cdot \delta_h(t)\}(s) = \lim_{h\rightarrow 0}
\int_{0_{-}}^{\infty} e^{-st} \frac{\xi_0h^{\mu}}{\xi_0h^{\mu} +(1-e^{-hD_t})^{\mu}} \delta_h(t) {\rm d}t    \\ \\ \ds \hspace{1cm}
=\lim_{h\rightarrow 0} \frac{\xi_0}{\xi_0+h^{-\mu}(1-e^{-hs})^{\mu}} {\cal L}\{\delta_h(t)\}(s)   = \frac{\xi_0}{\xi_0+s^{\mu}}
\end{array} \hspace{1cm} \mu >0
\end{equation}
where we made use of the limiting property $\lim_{h\rightarrow 0}{\cal L}\{\delta_h(t)\}(s)={\cal L}\{\delta(t)\}(s)=1$
of (\ref{Lapladeltah}).
We identify 
(\ref{Lapla}) with the Laplace transform of the Mittag-Leffler density. 
Now rewrite (\ref{limiting-process}) as a series of Riemann-Liouville fractional integral operators
\begin{equation}
\label{rewritemL}
\begin{array}{l} \ds 
{\cal I}_{\mu}\cdot \delta(t)  = \lim_{h\rightarrow 0} 
\sum_{n=1}^{\infty} (-1)^{n-1}\left(\xi_0 h^{\mu}(1-e^{-hD_t})^{-\mu}\right)^n \, \cdot \delta_h(t) \\ \\
\ds {\cal I}_{\mu}\cdot \delta(t) = \sum_{n=1}^{\infty} (-1)^{n-1}\xi_0^n D_t^{- n \mu} \,\cdot \delta(t) =\frac{\xi_0}{\xi_0+D_t^{\mu}} \cdot \delta(t) 
\end{array} \hspace{0.25cm} \mu >0, \hspace{0.25cm} t\in \mathbb{R}^{+}
\end{equation}
where the Riemann-Liouville integrals $D_t^{- n \mu} \,\cdot \delta(t)$ are defined by (\ref{limitprocess}). Thus we obtain
the continuous-time limit convolution kernel
\begin{equation}
\label{ML-int}
\begin{array}{l} \ds 
{\cal I}_{\mu} \cdot \delta(t) = \frac{\xi_0}{\xi_0+D_t^{\mu}} \cdot \delta(t) \\ \\
\ds \hspace{0.5cm}  = \int_0^t
\sum_{n=1}^{\infty} (-1)^{n-1}\xi_0^n\frac{\tau^{n\mu-1}}{\Gamma(n\mu)}\delta(t-\tau) {\rm d}\tau = \xi_0 t^{\mu-1} 
E_{\mu,\mu}(-\xi_0t^{\mu}) = \frac{d}{d t}(1-E_{\mu}( \xi_0t^{\mu}))
\end{array} t \in \mathbb{R}^{+}
\end{equation}
which we identify with the Mittag-Leffler density
where here we have $\mu>0$ and $\mu=1$ recovers the exponential density ${\cal I}_{1} \cdot \delta(t)=\xi_0e^{-\xi_0 t}$.
\\[1ex]
So far in this appendix we generalized well-scaled discrete-time fractional {\it integrals} which involve uniquely non-negative powers of the well-scaled discrete-time integral operator $ \frac{h}{1-{\hat T}_{-h}}$. Now we analyze operators involving its negative 
powers leading to the notion of generalized fractional derivatives. 
To this end we define the well-scaled `{\it discrete-time fractional derivative}' as
\begin{equation}
\label{further-important}
\begin{array}{l}
\ds
\left(\frac{h}{1-{\hat T}_{-h}}\right)^{-\mu} \cdot f(t) = h^{-\mu} (1-{\hat T}_{-h})^{\mu} \cdot \Theta(t) f(t) \\ \\ \ds \hspace{3.2cm} =
h^{-\mu} \sum_{k=0}^{\frac{t}{h}} \left(\begin{array}{l} \mu \\  k\end{array}\right)(-1)^k f(t-kh)
\end{array} \hspace{1cm} \mu >0, \hspace{1cm} t \in h\mathbb{N}_0 .
\end{equation}
The limit $h\rightarrow 0$ of this relation can be identified with the Gr\"unwald-Letnikov derivative, e.g. \cite{Meerschert-et-al2018}.
Then, recall the {\it ceiling function} $\ceil{\mu}$ indicates the smallest integer larger or equal to $\mu$.
We have $-1 <\mu-\ceil{\mu} \leq 0$ with $-1 <\mu-\ceil{\mu} <0$ if $\mu \notin \mathbb{N}$ and zero if $\mu = \ceil{\mu}$ is integer.
Then (\ref{further-important}) can be decomposed as follows
\begin{equation}
\label{decomp}
h^{-\mu} (1-{\hat T}_{-h})^{\mu} = h^{-\ceil{\mu}} (1-{\hat T}_{-h})^{\ceil{\mu}}\,\, h^{\ceil{\mu}-\mu} 
(1-{\hat T}_{-h})^{\mu-\ceil{\mu}}
\end{equation}
which is a discrete-time derivative of integer-order $\ceil{\mu}$
applied to a (generalized fractional discrete-time) integral of order $\ceil{\mu}-\mu$. 
The continuous-time limit of (\ref{further-important}) is then obtained by accounting for 
\begin{equation}
\label{integer-order}
\lim_{h\rightarrow 0}  h^{-\ceil{\mu}} (1-{\hat T}_{-h})^{\ceil{\mu}}  = \frac{d^{\ceil{\mu}}}{dt^{\ceil{\mu}}} = D_t^{\ceil{\mu}}
\end{equation}
converges to the integer-order traditional derivative of order $\ceil{\mu}$ and

\begin{equation}
\label{fracintegral}
\lim_{h\rightarrow 0} h^{\ceil{\mu}-\mu} (1-{\hat T}_{-h})^{\mu-\ceil{\mu}} \cdot f(t) 
\end{equation}
converges by taking into account (\ref{limitprocess}) to the {\it Riemann-Liouville fractional integral} of order $\ceil{\mu}-\mu$.
Hence the continuous-time limit of (\ref{further-important}) yields
\begin{equation}
\label{RLint}
D_t^{\mu}\cdot f(t)= \lim_{h\rightarrow 0} h^{-\mu} (1-{\hat T}_{-h})^{\mu} \cdot f(t) = D_t^{\ceil{\mu}} \int_0^t f(t-\tau) 
\frac{\tau^{\ceil{\mu}-\mu-1}}{\Gamma(\ceil{\mu}-\mu)}{\rm d}\tau ,\hspace{0.5cm} t \in \mathbb{R}^{+}
\end{equation}
which we identify with the {\it Riemann-Liouville fractional derivative} or order $\mu$ \cite{SamkoKilbasMarichev1993} where $\tau^{\ceil{\mu}-\mu-1}$ due to $-1 < \ceil{\mu}-\mu-1 <0 $ is if $\mu \notin \mathbb{N}_0$ a weakly singular integrable function. 
For $\ceil{\mu}-\mu \rightarrow 0$, i.e. in the limit of $\mu$ approaching integer we have
the distributional relation $\frac{\tau^{\ceil{\mu}-\mu-1}}{\Gamma(\ceil{\mu}-\mu)} \rightarrow \delta(\tau)$ \cite{GelfandShilov1968} thus (\ref{RLint}) converges to the integer order derivative $D_t^{\mu}\cdot f(t) \rightarrow D_t^{\ceil{\mu}}\cdot f(t)$.

\subsection{\small SIBUYA WALK}
\label{Sib}
We consider now a strictly increasing random walk with $Z_j \in \mathbb{N}$ are IID steps 
occurring with the probabilities
\begin{align}
\nonumber
\mathbb{P}(Z_j=r) &= w_{\alpha}(r)= (-1)^{r-1}\left(\begin{array}{l} \alpha \\ r \end{array}\right)\\ &= (-1)^{r-1}\,\frac{\alpha(\alpha-1)..(\alpha-r+1)}{r!} , \hspace{0.5cm} r \in \mathbb{N} ,\hspace{0.5cm} w_{\alpha}(r)\big|_{r=0} = 0,  
\hspace{0.3cm} \alpha \in (0,1). \label{Sibuya}
\end{align}
This distribution is referred to as {\it Sibuya$(\alpha)$} (See \cite{PachonPolitoRicciuti2018} for a 
thorough analysis of properties).
The limit $\alpha \rightarrow 1-0$ 
gives a `trivial walk' where the walker in 
each step hops from $r$ to $r+1$ (almost surely) with transition probabilities 
$w_1(r)= \delta_{r1}$ ($r\in \mathbb{N}_0$). Despite $\alpha=1$ is admissible, per definition  
{\it Sibuya$(\alpha)$} covers only the fractional interval $\alpha\in (0,1)$.
\\[1ex]
{\it Sibuya$(\alpha)$} has the generating function
\begin{equation}
\label{wgen}
{\bar w}_{\alpha}(u)=  \sum_{k=1}^{\infty} (-1)^{k-1}\left(\begin{array}{l} \alpha \\ k \end{array}\right)  u^k =1-(1-u)^{\alpha} ,\hspace{0.5cm} 
|u| \leq 1
\end{equation}
where $w_{\alpha}(r)=\frac{1}{r!}\frac{d^r}{du^r}{\bar w}(u)|_{u=0}$ gives the probabilities (\ref{Sibuya}).
We notice that for large $r$ using the asymptotic formula $\frac{\Gamma(r+c)}{\Gamma(r)} \sim r^c$ we have
\begin{align}
\label{fat-tail}
w_{\alpha}(r)&= (-1)^{r-1}\left(\begin{array}{l} \alpha \\ r\end{array}\right) \notag \\
& = \frac{\alpha}{r} \frac{\Gamma(r-\alpha)}{\Gamma(1-\alpha)\Gamma(r)} \sim
\frac{\alpha r^{-\alpha-1}}{\Gamma(1-\alpha)} = -\frac{r^{-\alpha-1}}{\Gamma(-\alpha)} ,\hspace{0.5cm} (r\,\,\, {\rm large}) ,\hspace{0.5cm} \alpha \in (0,1)
\end{align}
which is fat-tailed with the same tail as the Mittag-Leffler density.
\\[1ex]
An interesting quantity in a strictly increasing walk is the {\it expected hitting number} $\tau(r)$ of a node $r$. For $\alpha=1$ this quantity
is trivial and yields $\tau_{\alpha=1}(r)=1$ which is the upper bound in a strictly increasing walk. Since (\ref{Sibuya}) allows long-range jumps (with overleaping of nodes) and in view that the walk is strictly increasing so that nodes can be visited either once or never, intuitively we infer that the expected hitting number fulfills $0<\tau_{\alpha}(r)<1$ in the fractional interval $\alpha \in (0,1)$. This will be indeed confirmed subsequently by means of explicit formulas.
The expected hitting number $\tau_{\alpha}(r)$ indicates the average number of visits of a node $r$
and generally is an important quantity to describe recurrence phenomena 
\cite{michel-riascos-recurrence-frac-2017,michel-riascos-springer2018,TMM-APR-ISTE2019}. It is clear that a strictly increasing walk 
always is transient, since the walker cannot return to a node. Generally the expected hitting number of a node in a strictly increasing walk fulfills $0<\tau(r) \leq 1$. By simple conditioning arguments we have that
\begin{equation}
\label{explnumvisits}
\tau_{\alpha}(r) = \sum_{n=0}^{\infty} (w_{\alpha}\star)^n(r) 
\end{equation}
where $(w_{\alpha}\star)^n(r)$ indicates convolution powers (See Appendix \ref{generfuncts}).
For an analysis of this quantity for space-fractional Markovian walks on undirected networks we refer to \cite{TMM-APR-ISTE2019}.
It follows that (\ref{explnumvisits}) for {\it Sibuya$(\alpha)$} has the generating function  
\begin{equation}
\label{SibuyaGenFu}
{\cal G}_{\alpha}(u)= \sum_{n=0}^{\infty}({\bar w}_{\alpha}(u))^n  =\frac{1}{1-{\bar w}_{\alpha}(u)} = (1-u)^{-\alpha} ,\hspace{1cm} |u|<1 
\end{equation}
where ${\bar w}_{\alpha}(u)$ is the generating function (\ref{wgen}) of the Sibuya-steps.
The expected hitting number of a node $r \in \mathbb{N}_0$ in a Sibuya walk is then obtained as
\begin{equation}
\label{expectedhiting}
\begin{array}{l}
\ds 
\tau_{\alpha}(r) = \frac{1}{r!} \frac{d^r}{du^r} \frac{1}{1-{\bar w}_{\alpha}(u)}|_{u=0} = 
\frac{1}{r!} \frac{d^r}{du^r}(1-u)^{-\alpha}|_{u=0} \\ \\ \ds 
\hspace{0.5cm} = \frac{\alpha(\alpha+1)\ldots(\alpha+r-1)}{r!}= \frac{1}{\Gamma(\alpha)} \frac{\Gamma(\alpha+r)}{\Gamma(r+1)} = 
\frac{(\alpha)_r}{r!} 
\end{array} \hspace{1cm} r\in \mathbb{N}_0
\end{equation}
where $(\alpha)_r$ stands for the Pochhammer-symbol (\ref{Pochhammer}).
For $\alpha=1$ we have $\tau_{1}(r) =\frac{r!}{r!\Gamma(1)}=1$  $\forall r \in \mathbb{N}_0$ which is the 
above anticipated result when the walker hops in each step
from $r$ to $r+1$, thus the walker hits each node once (almost surely).
Hence  $\tau_{1}(r)=1$ is the upper bound for $\tau_{\alpha}(r)$ 
which can be seen accounting for
\begin{equation}
\label{beseen}
\tau_{\alpha}(r) = \prod_{j=1}^r \left(1 + \frac{\alpha-1}{j}\right), \hspace{1cm} 0<\alpha \leq 1
\end{equation}
where each multiplier fulfills $0 < 1 + \frac{\alpha-1}{j} \leq 1$ thus $0< \tau_{\alpha}(r) \leq 1$. The lower bound $\tau_{0+}(r)\sim \delta_{t0}$ (i.e. null for $r>0$) is approached for $\alpha\rightarrow {0+}$ where also for $r$ large we then have $\tau_{\alpha}(r) \sim  \frac{r^{\alpha-1}}{\Gamma(\alpha)} \sim \delta(t) = 0$.
The physical interpretation is that the steps then
become infinitely long, thus
any finite node $r>0$ is over-leaped where only the departure node $r=0$ 
is occupied once.
\\[1ex]
On the other hand for $\alpha\rightarrow 1-0$ we have $\tau_{1-0}(k)=1$ thus each node is hit once (almost surely) since the walker moves in any time step from $r$ to its right neighbor node $r+1$ (with transition probabilities $w_1(r)=\delta_{r1}$).
\\[1ex]
It is instructive to recall briefly the limiting distribution of 
the rescaled integer random variable $\sigma_{\alpha}(n) \rightarrow \lambda_n \sigma_{\alpha}(n)$ 
for infinitely many steps $n \rightarrow \infty$ 
where $\sigma(n) = \sum_{j=1}^nZ_j$ is an integer random variable with IID Sibuya steps $Z_j$.
Then let us consider the Laplace transform 
\begin{equation}
\label{limiting}
\begin{array}{l}
\ds 
\lim_{n\rightarrow \infty} \mathbb{E} e^{-s\lambda_n\sigma_{\alpha}(n)}  = 
\prod_{j=1}^n\left(\mathbb{E}e^{-s \lambda_nZ_j}\right) = 
\left(\sum_{r=0}^{\infty}e^{-s \lambda_n r}w_{\alpha}(r) \right)^n \\ \\
\ds \hspace{2.9cm} = \left(1-(1-e^{-\lambda_n s})^{\alpha}\right)^n = \left({\bar w}_{\alpha}(u=e^{-\lambda_ns})\right)^n .
\end{array}
\end{equation}
In this limiting relation 
we rescale the Sibuya random variable $\sigma(n) =\sum_{j=1}^nZ_j \rightarrow \lambda_n \sigma(n) $ 
such that the limit (\ref{limiting}) exists. This requires the scaling factor $\lambda_n \rightarrow 0$ 
and $1-(1-e^{-\lambda_n s})^{\alpha} = 1-s^{\alpha}\lambda_n^{\alpha} +O(\lambda_n^{\alpha}s^{\alpha})$ 
thus (\ref{limiting}) remains finite (i.e. neither null nor infinite)
\begin{equation}
\label{finitelimit}
\mathbb{E} e^{-s\lambda_n\sigma_{\alpha}(n)}  \sim \left(1-\lambda_n^{\alpha}s^{\alpha}\right)^n, \qquad n \to \infty
\end{equation} 
if and only if we choose $\lambda_n = (\frac{\xi}{n})^{\frac{1}{\alpha}}$ where $\xi$ is an arbitrary positive constant independent on $n$.
With this choice we get
\begin{equation}
\label{finitelimitrescaled}
\mathbb{E} e^{-s\lambda_n\sigma_{\alpha}(n)}  \sim \left(1-\frac{\xi}{n}s^{\alpha}\right)^n \rightarrow
e^{-\xi s^{\alpha}} .
\end{equation}
Equation (\ref{finitelimitrescaled}) is the Laplace transform of the limiting distribution emerging for many Sibuya steps 
($n\rightarrow \infty$) of the rescaled random variable $\frac{\xi^{\frac{1}{\alpha}}}{n^{\frac{1}{\alpha}}}\sigma_{\alpha}(n)$ 
which hence is a stable subordinator (see e.g.\ \cite{PachonPolitoRicciuti2018} and the references therein).
\subsection{\small SIBUYA COUNTING PROCESS}
\label{Sibuyarenewal}
In this appendix we review some pertinent results for the `Sibuya counting process'
(i.e. the discrete-time counting 
process with waiting times following Sibuya($\alpha$)) to develop the Sibuya Montroll-Weiss DTRW as a walk subordinated to the Sibuya counting process in Appendix \ref{Sibuya-DTRW}. The interested reader is referred to \cite{PachonPolitoRicciuti2018} for a profound analysis, and see also the references therein.
The generating function of the Sibuya
survival probability is then with (\ref{wgen}) obtained as
\begin{equation}
\label{survivalSibuya}
{\bar \Phi}_{S}^{(0)}(u) = \frac{1}{1-u}(1-{\bar w}_{\alpha}(u)) = (1-u)^{\alpha-1}
\end{equation}
which yields the survival probability
\begin{equation}
\label{Sibuyasurvival}
\Phi_{S,\alpha}^{(0)}(t) = \frac{1}{t!}\frac{d^t}{du^t} {\bar \Phi}_{S}^{(0)}(u)|_{u=0} =
\frac{1}{t!}\frac{d^t}{du^t}(1-u)^{\alpha-1}|_{u=0} = 
\frac{\Gamma(t+1-\alpha)}{\Gamma(t+1)\Gamma(1-\alpha)} ,\hspace{1cm} t\in \mathbb{N}_0 
\end{equation} 
fulfilling initial condition $\Phi_{S,\alpha}^{(0)}(t)\big|_{t=0}=1$.
We observe that for $\alpha\rightarrow 1-0$ the survival probability $\Phi_{S,1}(t)=\delta_{t0}$ for $t>0$ is vanishing 
since for $\alpha=1$ the waiting time is equal to one (almost surely with $w_1(t)=\delta_{t1}$). Now consider the 
behavior for $t$ large which yields
\begin{equation}
\label{klarge}
\Phi_{S,\alpha}^{(0)}(t)= \frac{\Gamma(t+1-\alpha)}{\Gamma(t+1)\Gamma(1-\alpha)} \sim 
\frac{t^{-\alpha}}{\Gamma(1-\alpha)} ,\hspace{1cm} \alpha \in (0,1)
\end{equation} 
with the same tail as has the Mittag-Leffler waiting time distribution. Hence for large observation times the Sibuya counting process behaves as a fractional Poisson process and the same is true for the PDTP (See asymptotic relations (\ref{state-asymp}) for the PDTP state probabilities.).
We mention that the state-probabilities, i.e. probabilities for $n$ arrivals in the discrete time interval $[0,t]$ are given by
\begin{equation}
\label{steprob}
\Phi_{S,\alpha}^{(n)}(t) = \frac{1}{t!}\frac{d^t}{du^t} 
\left\{(1-u)^{\alpha-1}\left(1-(1-u)^{\alpha}\right)^n\right\}|_{u=0} ,\hspace{0.5cm} n,t \in \mathbb{N}_0 
\end{equation}
with generating function 
$ \frac{(1-{\bar w}_{\alpha}(u))}{(1-u)} ({\bar w}_{\alpha}(u))^n=  {\bar \Phi}_{S,\alpha}^{(n)}(u)= (1-u)^{\alpha-1}\left(1-(1-u)^{\alpha}\right)^n$ and refer to Ref. \cite{PachonPolitoRicciuti2018} for a thorough analysis. 
\subsection{\small SIBUYA DTRW}
\label{Sibuya-DTRW}
Due to the general importance of {\it Sibuya}$({\alpha})$ let us briefly outline
the {\it Sibuya Montroll-Weiss DTRW} where we subordinate a walk with one step-transition matrix ${\mathbf H}$ to a Sibuya counting process. 
Then we obtain for the generating function of the transition matrix the Cox-series (See (\ref{Coxseries}))
\begin{equation}
\label{Cox-ser-sib}
{\bar {\mathbf P}}_S(u) = 
\sum_{n=0}^{\infty}{\mathbf H}^n  {\bar \Phi}_{S,\alpha}^{(n)}(u) = (1-u)^{\alpha-1}
\left[{\mathbf 1} - {\mathbf H} + {\mathbf H}(1-u)^{\alpha}\right]^{-1} .
\end{equation}
We call the random walk with transition matrix generating function (\ref{Cox-ser-sib})
the `{\it Sibuya DTRW}'. 
Formula (\ref{Cox-ser-sib}) is related to the
generating function
\begin{equation}
\label{double-gen-fou-mod}
\begin{array}{l}
\ds 
{\cal G}_{S,\alpha}(u,v)= \sum_{n=0}^{\infty}\sum_{r=0}^{\infty} \Phi_{S,\alpha}^{(n)}(r)  u^r v^n  =  \frac{(1-{\bar w}_{\alpha}(u))}{1-u} \sum_{n=0}^{\infty}({\bar w}_{\alpha}(u))^nv^n  \\ \\ \ds = \frac{(1-{\bar w}_{\alpha}(u))}{1-u} \frac{1}{1-v{\bar w}_{\alpha}(u)} = \frac{(1-u)^{\alpha-1}}{1-v+v(1-u)^{\alpha}} ,\hspace{0.5cm}
|u| \leq 1, \hspace{0.5cm} |v| <1
\end{array}
\end{equation}
by $v \rightarrow {\mathbf H}$.
Now it is only a small step to write down the Kolmogorov-Feller type equations for
the Sibuya DTRW, namely
\begin{equation}
\label{SibyaDTRWkG}
(1-{\hat T}_{-1})^{\alpha} {\mathbf H} {\mathbf P}_S(t) = 
({\mathbf H}-{\mathbf 1}) {\mathbf P}_S(t) + {\mathbf 1}(1-{\hat T}_{-1})^{\alpha-1} \delta_{t0}
\end{equation}
or
\begin{equation}
\label{equisibkf}
(1-{\hat T}_{-1})^{\alpha} {\mathbf H} {\mathbf P}_S(t)  = 
({\mathbf H}-{\mathbf 1}) {\mathbf P}_S(t)  +  {\mathbf 1}(-1)^t \left(\begin{array}{l} \alpha-1 \\  \hspace{0.3cm} t\end{array}\right).
\end{equation}
This equation is the Kolmogorov-Feller fractional difference equation solved by the transition matrix $ {\mathbf P}_S(t)$ in a Montroll-Weiss {\it Sibuya DTRW}.
\\[1ex]
Let us also consider a walk not taking place on an undirected network, but the simplest variant of strictly increasing walk 
where the walker at each event makes a unit step into the positive integer direction. This walk simply counts the arrivals and has the transition probabilities $H_{ij}=\delta_{i+1,j}$ i.e. we have then in (\ref{equisibkf})
$({\mathbf H} {\mathbf P}_S(t))_{ij} =
\sum_{k=0}^{\infty} \delta_{i+1,k} P_{kj} = P_{i+1,j}$. For this walk Eq. (\ref{equisibkf}) then takes the form
\begin{equation}
\label{wegetSibuyKFaeqs}
(1-{\hat T}_{-1})^{\alpha} P_{i+1,j}(t) - (-1)^t \left(\begin{array}{l} \alpha-1 \\  \hspace{0.3cm} t\end{array}\right)\delta_{ij} =  P_{i+1,j}(t)-P_{ij}(t) 
\end{equation}
where this equation is also obtained in \cite{PachonPolitoRicciuti2018} (Remark 7 in that paper).
\end{appendix}

\end{document}